\documentclass[review]{elsarticle}
\usepackage{multirow}
\usepackage{booktabs}
\usepackage{ulem}
\usepackage{lineno,hyperref}
\usepackage{amsmath,bm}
\usepackage{epstopdf}
\usepackage{color,xcolor}
\modulolinenumbers[5]

\journal{Applied Mathematics and Computation}

\begin{document}

\begin{frontmatter}

\title{A high-order compact finite difference scheme and precise integration method based on modified Hopf-Cole transformation for numerical simulation of n-dimensional Burgers' system}

\author[mymainaddress]{Changkai Chen}

\author[mymainaddress,mysecondaryaddress]{Xiaohua Zhang\corref{mycorrespondingauthor}}
\cortext[mycorrespondingauthor]{Corresponding author.\newline
E-mail addresses: zhangxiaohua07@163.com;}

\author[mymainaddress]{Zhang Liu}

\address[mymainaddress]{College of Science, China Three Gorges University, Yichang, 443002, China}
\address[mysecondaryaddress]{Three Gorges Mathematical Research Center, China Three Gorges University, Yichang, 443002, China}

\begin{abstract}
This paper modifies a n-dimensional Hopf-Cole transformation to the n-dimensional Burgers' system.
We obtain the n-dimensional heat conduction equation through the modification of the Hopf-Cole transformation.
Then the fourth-order precise integration method (PIM) in combination with a spatially global sixth-order compact finite difference (CFD) scheme is presented to solve the equation with high accuracy.
Moreover, coupling with the Strang splitting method, the scheme is extended to multi-dimensional (two, three-dimensional) Burgers' system.
Numerical results show that the proposed method appreciably improves the computational accuracy compared with the existing numerical method.
Moreover, the two-dimensional and three-dimensional examples demonstrate excellent adaptability, and the numerical simulation results also have very high accuracy in medium Reynolds numbers.
\\

\end{abstract}

\begin{keyword}
n-dimensional Hopf-Cole transformation, n-dimensional Burges' system, compact finite difference, precise integration method, Strang splitting method
\end{keyword}

\end{frontmatter}

\section{Introduction}
Burgers' equation is a nonlinear partial differential {equation (PDE)} which was first introduced by Bateman \cite{BATEMAN}, and was later treated as the turbulence of the mathematical model \cite{book:968450,BURGERS1948171}.
Burgers' equation is an especially important PDEs in fluid mechanics,
which combines the characteristics of the first order wave equation and heat conduction equation.
Burgers' equation is used as a tool to describe the interaction between convection and diffusion.
Over the decades, Burgers' equation has a large variety of applications in the modeling of water in dynamic soil water, surface disturbances electromagnetic waves, density waves, statistics of flow problems, mixing and turbulent diffusion, cosmology and seismology \cite{Zabusky1965,Dubansky2017}, etc.
Hopf \cite{Hopf1950} and Cole \cite{Cole1951} showed independently that for any given initial conditions the Burgers' equation can be reduced to a linear homogeneous heat equation that can be solved analytically and the analytical solution of the old Burgers' equation can be expressed in the form of Fourier series.
Even though the analytical solution is available in the form of the Fourier series, accurate and efficient numerical schemes are still required to solve the Burgers' equation which consists of a multi-dimensional system or the complex initial condition.
In such situations, the Fourier series solutions for the practical applications are very limited, which converges slowly or diverges in many cases.
The analytical or numerical solutions are essential for the corresponding Burgers' equations.
Apart from the limited number of these problems, most of them do not have exact analytical solutions, so it is imperative to get a satisfactory solution of Burgers' equation.
Here, we first analyze one-dimensional coupled Burgers' equation.
Owing to the nonlinear convection term and viscous term, the coupled Burgers' equation can be studied as a simple example of the Navier-Stokes equation.
\begin{itemize}
\item The one-dimensional coupled nonlinear Burgers' equation \cite{Kumar2014}
\end{itemize}

\begin{equation}\label{equ.1}
  \left\{ \begin{array}{l}
\frac{{\partial u_1}}{{\partial t}} = \omega_1 \frac{{{\partial ^2}u_1}}{{\partial {x^2}}} - \kappa_1 u_1 \frac{{\partial u_1}}{{\partial x}} - \delta_1 \frac{{\partial ({u_1}{u_2})}}{{\partial x}}\\
\frac{{\partial u_2}}{{\partial t}} = \omega_2 \frac{{{\partial ^2}u_2}}{{\partial {x^2}}} - \kappa_2 u_2 \frac{{\partial u_2}}{{\partial x}} - \delta_2\frac{{\partial ({u_1}{u_2})}}{{\partial x}}
\end{array} \right. ,{\kern 1pt} {\kern 1pt} {\kern 1pt} {\kern 1pt} {x \in \Omega  = [a,b] },{\kern 1pt} {\kern 1pt} t \in [0,T]
\end{equation}
subject to the initial conditions:
\begin{equation}\label{}
\begin{array}{l}
u_1(x,0) = {g_1}(t){\kern 1pt} {\kern 1pt} ,{\kern 1pt} {\kern 1pt} {\kern 1pt}{ x \in \Omega  = [a,b]}\\
u_2(x,0) = {g_2}(t){\kern 1pt} {\kern 1pt} ,{\kern 1pt} {\kern 1pt} {\kern 1pt}{ x \in \Omega  = [a,b]}
\end{array}
\end{equation}
and boundary conditions

\begin{equation}\label{}
\begin{array}{l}
u_1(a,t) = {f_1}(x){\kern 1pt} {\kern 1pt} ,{\kern 1pt} {\kern 1pt} u_1(a,t) = {f_2}(x){\kern 1pt} {\kern 1pt} ,{\kern 1pt} {\kern 1pt} {\kern 1pt} { x \in \Omega  = [a,b]}{\kern 1pt} {\kern 1pt} {\kern 1pt} ,{\kern 1pt} {\kern 1pt} {\kern 1pt}{\kern 1pt} {\kern 1pt} t \in [0,T]\\
u_2(b,t) = {f_3}(x){\kern 1pt} {\kern 1pt} ,{\kern 1pt} {\kern 1pt} {\kern 1pt} u_2(b,t) = {f_4}(x){\kern 1pt} {\kern 1pt} ,{\kern 1pt} {\kern 1pt} {\kern 1pt} {\kern 1pt} {x \in \Omega  = [a,b]}{\kern 1pt} {\kern 1pt} ,{\kern 1pt} {\kern 1pt} {\kern 1pt}{\kern 1pt} {\kern 1pt} t\in [0,T]
\end{array}
\end{equation}
where $\omega_1,\omega_2$ is kinematic viscosity parameters of the fluid, which correspond to an inverse of Reynolds number $Re$ (if $\omega_1=\omega_2$, then $\omega_1=\omega_2=\frac{1}{Re}$ ) . $\kappa_1,\kappa_2$  are real constants and $\delta_1,\delta_2$ are arbitrary constants. $f(x)$ and $g(t)$ are given smooth functions.

It is pervasively acknowledged that the nonlinear coupled Burgers' equation (\ref{equ.1}) does not have precise analytic solutions.
Researchers are interested in using various numerical techniques to study the properties of Burgers' equation, because it has wide applicability in various fields of science and engineering.
Due to the existence of nonlinear terms and viscosity parameters, numerical approximation of nonlinear coupled Burgers' equation is a challenging task.
In Burgers' equation, discontinuities may appear in finite time, even if the initial condition is smooth. They give rise to the phenomenon of shock waves which have important applications in physics \cite{BREZIS199876,Frisch2001}.
Recently, some contributions related to the time-dependent coupled viscous Burgers' equation have been published, which analyze the theoretical and numerical aspects.
Several numerical experiments on non-coupled and coupled Burgers' equation were run to compare the accuracy of the proposed schemes with other existing methods \cite{Bhatt2016,Gao2016,Khater2008,Mittal2011,Jaradat2018,Jiwari2013}.
For the sake of clarity, a brief description of the comparison method is provided below.

Bhatt et al. \cite{Bhatt2016} proposed A-stable and L-stable Fourth-order exponential time difference Runge-Kutta schemes in combination with a global fourth-order CFD scheme for the numerical solution of the coupled Burgers' equations.

In Ref. \cite{Gao2016}, the analytical solutions of two-dimensional and three-dimensional Burgers' equations are derived.
For multi-dimensional problems, these solutions can describe the shock wave phenomenon in large Reynolds numbers ($R e\geq100$), { which can be used as a reference for testing numerical methods.}

In Ref. \cite{Khater2008}, the authors develop a Chebyshev spectral collocation method for solving approximate solutions of nonlinear PDEs.  Using Chebyshev spectral collocation method, this problem is reduced to a set of ordinary differential equations (ODEs), and then solved with Runge-Kutta fourth-order method.

In Ref. \cite{Mittal2011}, the author uses a cubic B-spline function to construct a collocation method for numerical simulation of coupled Burgers' equation.
The time derivative term is discretized by the conventional Crank-Nicolson (C-N) scheme, while space derivative term is discretized by the cubic B-spline method.
The results obtained by finite difference cubic B-spline show that the accuracy of the solution decreases with the increase of time due to the time truncation error of the time derivative term.

Jaradat et al. \cite{Jaradat2018} establish new two-mode coupled Burgers' equations which are introduced.
The authors find the necessary conditions in which the multiple {kinks} and multiple singular kink solutions exist and present the two-front solutions.

{Jiwari et al. \cite{Jiwari2013,Mittal2013} developed a differential quadrature method to solve time dependent Burgers' equation. And results are accurately produced by such two numerical schemes provided by the authors. These two schemes are also found quite easy to implement.}

{In Ref. \cite{Tamsir2016}, the authors proposed an algorithm based on exponential modified cubic-B-spline differential quadrature method for Burgers' equation. With some modifications, such method is flexible enough to solve model equations in multi-dimensional problems including mechanical, physical or biophysical effects.}

{In Ref. \cite{Olshanskii2007}, finite element analysis and approximation of Burgers'-Fisher equation with non-smooth initial data is presented, which provides results with high accuracy and efficiency.}

 In this paper, we mainly discuss the numerical scheme of n-dimensional Burgers' system. { The n-dimensional ($n\geq 1$) Burgers' system includes low-dimensional ($n=1$) Burgers' equations and multi-dimensional ($n\geq 2$) Burgers' equations.}

\begin{itemize}
\item The n-dimensional Burgers' system \cite{Chan2010,Chen2016,Wang2017}
\end{itemize}

\begin{equation}\label{equ.4}
  {\bm{U}_t} + {\kappa}(\bm{U} \cdot \nabla )\bm{U} = \omega \Delta \bm{U}
\end{equation}
where $\bm{U} = ( u_1, u_2,...,u_n )^T$  is the fluid velocity fields, $n$ is dimension of space, $\omega$ is the kinematic viscosity of the fluid and $\kappa$  are real constants.
$\Delta  = \frac{{{\partial ^2}}}{{\partial x_1^2}} + \frac{{{\partial ^2}}}{{\partial x_2^2}} +  \cdots  + \frac{{{\partial ^2}}}{{\partial x_n^2}}$ denotes the Laplace operator,
and $\nabla  = {(\frac{\partial }{{\partial {x_1}}} , \frac{\partial }{{\partial {x_2}}} ,  \cdots  , \frac{\partial }{{\partial {x_n}}})^T}$ is the Hamilton gradient operator.

The system of Eq. (\ref{equ.4}) is Burgers' equations which includes non-coupled and coupled problems.
When $n = 1,2,3$ and $\kappa=1$, the system of Eq.(\ref{equ.4}) respectively becomes one-dimensional, two-dimensional and three-dimensional Burgers' equation
 \begin{equation}\label{equ.5}
  {u_t} + u{u_x} = \omega {u_{xx}}
\end{equation}

 \begin{equation}\label{equ.6}
\begin{array}{l}
{u_t} + u{u_x} + v{u_y} = \omega ({u_{xx}} + {u_{yy}})\\
{v_t} + u{v_x} + v{v_y} = \omega ({v_{xx}} + {v_{yy}})
\end{array}
 \end{equation}
 \begin{equation}\label{equ.7}
   \begin{array}{l}
{u_t} + u{u_x} + v{u_y} + w{u_z} = \omega ({u_{xx}} + {u_{yy}} + {u_{zz}})\\
{v_t} + u{v_x} + v{v_y} + w{v_z} = \omega ({v_{xx}} + {v_{yy}} + {v_{zz}})\\
{w_t} + u{w_x} + v{w_y} + w{w_z} = \omega ({w_{xx}} + {w_{yy}} + {w_{zz}})
\end{array}
 \end{equation}

Chen et al. \cite{Chen2016} show an n-dimensional Hopf-Cole transformation between the n-dimensional Burgers' system and an n-dimensional heat equation under an irrotational condition. Motivated by this idea, the purpose of this paper is to intend to extend the Hopf-Cole transformation to linearize the n-dimensional Burgers' equation (\ref{equ.4});
After obtaining the n-dimensional heat conduction equation, the CFD scheme with high precision and high efficiency is used to solve it.

Currently, there are many numerical methods for heat conduction equation \cite{Andrea2001,Marsden1994}, such as finite difference method (FDM), finite element method (FEM), finite volume method (FVM) and spectrum method, etc.
The traditional FDM shows great limitations in accuracy. An important measure to improve the accuracy of the traditional FDM is to refine mesh, which in turn will increase the amount of storage and calculating speed, especially in high-dimensional cases.
Therefore, it is of great theoretical significance and practical value to construct a scheme with high accuracy and excellent stability in time and space.

The CFD scheme is one of the most studied FDM at present.
Experience proves that the compact scheme is much more accurate than the corresponding explicit scheme of the same order \cite{Li2008}.
Over the past three decades, the methods for developing high-order CFD scheme have made great progress.
Dennis et. al. proposed the fourth-order CFD scheme for convection-diffusion problems \cite{Dennis1989}, this scheme can get more accurate results with a thicker grid.
Lele \cite{Lele1992} developed CFD scheme with pseudo spectral resolution on the basis of summarizing the previous work and proposed a linear sixth-order central CFD scheme,
which can achieve the accuracy of the spectral method.
Subsequently, many scholars constructed different schemes of CFD scheme and solved many types of partial differential equations \cite{Zhao2013,Lai2007,Nihei2003,Sutmann2007},
such as integro-differential equations, three-dimensional Poisson equations, the shallow water equations, and the Helmholtz equations, they all achieved better numerical results.
Sengupta et. al. developed a class of upwind compact difference schemes, and such schemes could be applied to different fields \cite{Wang2005}.
In that same year, Kumar \cite{Kumar2009} discussed a high-order compact difference scheme for singularly perturbed reaction diffusion problems on a new Shish Kin mesh.
Sen \cite{Mehra2017,Sen2016} discussed the fourth-order exact compact difference scheme for mixed derivative parabolic problems with variable coefficients.

The CFD scheme is a widely used method for spatial discretization of heat conduction equations to obtain the ODEs, and then other methods of time discretization are used for discretizing the ODEs, such as Euler method, multistep methods and Runge-Kutta method.
The exact solution of heat conduction equation contains the calculation of exponential matrix.
How to accurately calculate the exponential matrices is an essential problem in solving PDEs.
Moler et al. \cite{Moler2003} summarized nineteen schemes for calculating the exponential matrices.
These nineteen schemes are aimed at different practical problems, and their numerical solutions also have corresponding advantages and disadvantages.
In 1994, Zhong \cite{Wan-Xie2004} proposed the precise integration method (PIM) of exponential matrices to solve the initial value problem of linear ODEs.
PIM is an approximated method to calculate the exponential matrices, which contains Taylor approximation and Pad\'{e} approximation.
The PIM avoided the computer error caused by fine division and improved the numerical solution of exponential matrices by the accuracy of the computation.

Alternating Direction Implicit (ADI) method is a classical numerical scheme for solving multi-dimensional heat conduction equation.
ADI, such as Peacemen-Rachford scheme, D'Yakonov scheme and Douglas scheme, are only the second-order accuracy schemes \cite{Zhang2016,Karaa2004,Journal2013}.
ADI often fail to meet the accuracy requirements of practical problems.
Strang splitting method (SSM) is a numerical method for solving differential equations that are decomposable into a sum of differential operators,
which is to solve multi-dimensional PDEs by reducing their dimensionality to a sum of one-dimensional problems \cite{Gilbert2013}.
This is a scheme of operator splitting method.
If the differential operators of the SSM commute, then it will lead to no loss of accuracy.
Therefore, the proposed schemes will extend to multi-dimensional heat conduction equation through SSM.

The remainder of the paper is arranged as follows.
The n-dimensional Hopf-Cole transformation between the n-dimensional Burgers' system and n-dimensional heat conduction equation are presented in Section 2; Moreover, we give the modification of the Hopf-Cole transformation.
The high-order exponential time differencing PIM in combination with a spatially global sixth-order CFD scheme for solving n-dimensional heat condution equations are presented in Section 3. In Section 4, the Strang splitting method is described and the proposed schemes are extended to multi-dimensional problems. In Section 5, numerical examples are carried out to test the accuracy and adaptability of the proposed schemes. The conclusions are drawn in Section 6.

\section{{The} n-dimensional {Hopf-Cole} transformation}
The purpose of n-dimensional Hopf-Cole transformation is to convert Eq. (\ref{equ.4}) into {the} n-dimensional heat equation
\begin{equation}\label{equ.8}
 { {{\phi} _t} - \omega {\triangle {\phi}} = 0}
\end{equation}
by the n-dimensional Hopf-Cole transformation
\begin{equation}\label{equ.9}
 {{u_i} =  - 2\frac{\omega }{\kappa }{\partial _{{x_i}}}\ln {\phi } =  - 2\frac{\omega }{\kappa }\frac{{{\phi _{{x_i}}}}}{\phi }}
\end{equation}
 {where $i=1,2,...,n$. Note: $x_1=x,x_2=y,x_3=z$; $u_1=u,u_2=v,u_3=w$.}

When $n = 1,2,3$ and $\kappa=1$, the system of Eq. (\ref{equ.8}) respectively becomes one-dimensional, two-dimensional and three-dimensional heat equations
 {
\begin{equation}\label{e.10}
  {{\phi} _t} - \omega {{\phi_{xx}}} = 0
\end{equation}
\begin{equation}\label{e.11}
  {{\phi} _t} - \omega {({\phi_{xx}}+{\phi_{yy}})} = 0
\end{equation}
\begin{equation}\label{e.12}
  {{\phi} _t} - \omega {({\phi_{xx}}+{\phi_{yy}+{\phi_{zz}}})} = 0
\end{equation}
}

The initial and boundary conditions are
\begin{equation}\label{equ.13}
  { {\phi}({x_i},0) = \exp(\int_0^{\bm{x_i}} { - \frac{{{{u_i}(\xi,0)}}}{{2\omega }}} d\xi),{\kern 1pt} {\kern 1pt} {\kern 1pt} {x_i} \in \Omega  = [a,b]{\kern 1pt} {\kern 1pt} {\kern 1pt}}
\end{equation}
\begin{equation}\label{equ.14}
 {{{\phi_{x_i}}}(a,t)={\phi_{x_i}}(b,t)=0,{\kern 1pt} {\kern 1pt} {\kern 1pt}{\kern 1pt} {\kern 1pt} t \in [0,T]}
\end{equation}

Based on this method, we intend to extend Hopf-Cole transformation to  n-dimensional Burgers' system. Assuming that the n-dimensional heat conduction equation has the irrotational condition
\begin{equation}\label{equ.10}
  \nabla  \times \bm{U} = \sum\limits_{i,j = 1}^n {(\frac{{\partial {u_j}}}{{\partial {x_i}}} - \frac{{\partial {u_i}}}{{\partial {x_j}}}){e_i} \wedge {e_j}}  = 0
\end{equation}
where ${e_i}{\kern 1pt} {\kern 1pt} ,{\kern 1pt} {\kern 1pt} {e_j}$ are the basis of n-dimensional Euclidean
space.

To facilitate readers to understand the derivation process, Eqs. (\ref{equ.4}) and (\ref{equ.10}) can be written as the following scalar forms
\begin{equation}\label{equ.11}
 { \frac{{\partial {u_i}}}{{\partial t}} + \kappa \sum\limits_{j = 1}^n {{u_j}\frac{{\partial {u_i}}}{{\partial {x_j}}}}  = \omega \Delta u_i ,{\kern 1pt} {\kern 1pt} i = 1,2, \ldots ,n}
\end{equation}
\begin{equation}\label{equ.17}
  \frac{{\partial {u_i}}}{{\partial {x_j}}} = \frac{{\partial {u_j}}}{{\partial {x_i}}},{\kern 1pt} {\kern 1pt} {\kern 1pt} {\kern 1pt} i,j = 1,2, \ldots ,n{\kern 1pt} {\kern 1pt} {\kern 1pt} ({\kern 1pt} i \ne j)
\end{equation}

Let $u_i={\frac{{\partial \varphi }}{{\partial {x_i}}}}$, substituting in Eq. (\ref{equ.11}), we obtain
\begin{equation}\label{equ.18}
 \frac{{{\partial ^2}\varphi }}{{\partial {x_j}\partial t}} + \kappa \sum\limits_{j = 1}^n {\frac{{\partial \varphi }}{{\partial {x_j}}}\frac{{{\partial ^2}\varphi }}{{\partial {x_i}\partial {x_j}}}} - \omega \frac{{\partial \Delta \varphi }}{{\partial {x_j}}}=0,{\kern 1pt} {\kern 1pt} i = 1,2, \ldots ,n
\end{equation}
then Eq. (\ref{equ.18}) can be written as
\begin{equation}\label{14}
{\frac{{\partial {\varphi}}}{{\partial t}} + \frac{\kappa}{2} \sum\limits_{j = 1}^n {{{(\frac{{\partial \varphi }}{{\partial {x_j}}})}^2} - \omega \Delta \varphi } }  = 0
\end{equation}

Applying Hopf-Cole transformation, Eq. (\ref{14}) will becomes the n-dimensional heat conduction equations. {The detailed derivation process is as follows:
\\ (1) Eq. (\ref{14}) can be written as
\begin{equation}\label{eeqq.20}
\frac{{\partial \varphi }}{{\partial t}} + \frac{\kappa }{2}{\sum\limits_{j = 1}^n {\left( {\frac{{\partial \varphi }}{{\partial {x_j}}}} \right)} ^2} - \omega \sum\limits_{j = 1}^n {\frac{{{\partial ^2}\varphi }}{{\partial x_j^2}}}  = 0
\end{equation}
where $\sum\limits_{j = 1}^n {\frac{{{\partial ^2}}}{{\partial x_j^2}}}  = \Delta$.
\\ (2) Introduce $\varphi=-2\frac{\omega }{\kappa }\ln{\bm{\phi}}$ for the system of Eq. (\ref{eeqq.20})
\begin{equation}\label{eeqq.21}
- \frac{{2\omega }}{\kappa }\frac{{{\phi _t}}}{\phi } + \frac{\kappa }{2}{\left( {\frac{{ - 2\omega }}{\kappa }\frac{1}{\phi }\sum\limits_{j = 1}^n {\frac{{\partial \phi }}{{\partial {x_j}}}} } \right)^2} + \frac{{2{\omega ^2}}}{\kappa }\frac{1}{{{\phi ^2}}}\left( {\sum\limits_{j = 1}^n {\left( {\frac{{{\partial ^2}\phi }}{{\partial x_j^2}}} \right)\phi  - \sum\limits_{j = 1}^n {{{\left( {\frac{{\partial \phi }}{{\partial {x_j}}}} \right)}^2}} } } \right) = 0
\end{equation}
\\ (3) The two sides of Eq.(\ref{eeqq.21}) are multiplied by $\kappa$ and then simplified.
\begin{equation}\label{eeqq.22}
  - 2\omega \frac{{{\phi _t}}}{\phi } + 2{\omega ^2}\frac{1}{{{\phi ^2}}}\sum\limits_{j = 1}^n {{{\left( {\frac{{\partial \phi }}{{\partial {x_j}}}} \right)}^2}}  - 2{\omega ^2}\frac{1}{\phi }\sum\limits_{j = 1}^n {\left( {\frac{{{\partial ^2}\phi }}{{\partial x_j^2}}} \right)}  - 2{\omega ^2}\frac{1}{{{\phi ^2}}}\sum\limits_{j = 1}^n {{{\left( {\frac{{\partial \phi }}{{\partial {x_j}}}} \right)}^2}}  = 0
\end{equation}
It is especially noted that $\kappa$ disappears in Eq. (\ref{eeqq.22}).
\\(4) And further simplify to obtain
\begin{equation}\label{eeqq.23}
  {\phi _t} + \omega \sum\limits_{j = 1}^n {\frac{{{\partial ^2}\phi }}{{\partial x_j^2}}}  = 0\Rightarrow {\phi _t} + \omega \Delta \phi  = 0
\end{equation}
}
{\subsection{The modification of Hopf-Cole transformation}}

 With the Development of Hopf-Cole transformation in the past decades,
 {Kadalbajoo et al. \cite{Kadalbajoo2006} proposed the C-N scheme based on the Hopf-Cole transformation for Eq. (\ref{equ.5}) }.
  They discretized the space twice with C-N scheme and central difference.
  Due to {the twice spatial dispersions} of {${\phi_x}$} of Eq. (\ref{equ.9}), the numerical solution results were in loss of accuracy.
  In 2015, Mukundan et al. \cite{Mukundan2015} presented numerical techniques for Burgers' equation, which use backward difference and central difference for {${\phi_x}$}.
  The accuracy of these numerical schemes will decline because of the twice discretizations of {${\phi_x}$}.
  We have improved Hopf-Cole transformation, which will only be dispersed once in space. {Hopf-Cole transformation is used again, but the object to be solved this time is the first derivative ${\phi_x}$ of the heat conduction equation.} Firstly, Eq. (\ref{equ.18}) can be written as

\begin{equation}\label{15}
{\frac{{\partial {\varphi_x}}}{{\partial t}} + \frac{\kappa}{2}\sum\limits_{j = 1}^N {{{(\frac{{\partial \varphi_x }}{{\partial {x_j}}})}^2} - \omega \Delta \varphi_x } }  = 0
\end{equation}

Substituting {$\varphi=-2\frac{\omega }{\kappa }\ln{{{\phi_x}}}$} into Eq. (\ref{15})
\begin{equation}\label{equ.16}
 { \frac{{\partial {\phi_x} }}{{\partial t}} - \omega {\triangle {\phi_x}} = 0}
\end{equation}

{Initial and boundary conditions of Eq. (\ref{equ.16}) can be obtained from Eqs. (\ref{equ.13}) and (\ref{equ.14}).}

Then the solution of Eq. (\ref{equ.16}) can be obtained by utilizing high precision numerical schemes such as CFD scheme.
Thus, Eq. (\ref{equ.16}) will get  {${\phi_x}$} after a spatial discretization for n-dimensional Burgers' equations.
In this way, the modification of n-dimensional Hopf-Cole transformation avoids the truncation error of twice spatial difference and can obtain the first derivative of  {${\phi}({x},t)$} with higher precision.

{
The modification of n-dimensional Hopf-Cole transformation design in this section lies in two points:
\\(1) The modification of Hopf-Cole transformation is more general and suitable for Burgers' system where $\kappa$ is a variable;
\\(2) Hopf-Cole transformation is used twice to solve the first derivative ${\phi_x}$ and solution ${\phi}$ of the heat conduction equation, thus avoiding the second truncation error.}

\subsection{The simplification of initial value problem}
For some initial value problems, Fourier series solutions of Hopf-Cole transformation will converge very slowly, which dramatically increases the complexity of the calculation.
In this subsection, we simplify the initial value condition of n-dimensional(one-dimensional, two-dimensional, three-dimensional) Burgers' equation.
{In 2016, Gao et al. \cite{Gao2016} gave numerical modification of analytical solution for two and three dimensional Burgers' equation. Their modification is similar to our simplification, but Gao et al. did not provide one-dimensional case. Therefore, the following two-dimensional and three-dimensional improvements refer to the ideas put forward by Gao et al.}

\subsubsection{One-dimensional modification}\label{2.2.1}
Researchers have proposed the one-dimensional Burgers' equation with the following initial and boundary condition \cite{Kadalbajoo2006,Jiwari2015,Jiwari2012,Seydaoglu2018,G.W.Wei1998}
\begin{equation}\label{e.60}
\begin{array}{l}
u(x,0) =u_0(x)= \sin \pi x{\kern 1pt} {\kern 1pt} {\kern 1pt} ,{\kern 1pt} {\kern 1pt} {\kern 1pt} x \in [0,{\rm{1}}]\\
u(0,t) = u(1,t) = 0{\kern 1pt} {\kern 1pt} ,{\kern 1pt} {\kern 1pt} t > 0
\end{array}
\end{equation}

It is widely noted that the analytical solution of {the} one-dimensional heat conduction equation can be written in the standard form of the Fourier series
\begin{equation}\label{eee.27}
\phi (x,t) = \sum\limits_{\alpha  = 0}^\infty  {{C_\alpha }\exp ( - {\alpha ^2}{\pi ^2}\omega t)\cos (\alpha \pi x)}
\end{equation}
where ${C_\alpha }$ is Fourier coefficient.

The initial conditions of {the} one-dimensional heat conduction equation are extracted from the  {Eq. (\ref{eee.27})}
\begin{equation}\label{ee.23}
 { \phi (x,0) = \sum\limits_{\alpha  = 0}^\infty  {{C_\alpha }\cos ({\alpha \pi x})} ,{\kern 1pt} {\kern 1pt} {\kern 1pt}{x} \in \Omega  = [0,1]}
\end{equation}
with boundary conditions
\begin{equation}\label{ee.25}
{  {\phi _x}(0,t) = {\phi _x}(1,t) = 0{\kern 1pt} {\kern 1pt} ,{\kern 1pt} {\kern 1pt} t \in [0,T]}
\end{equation}
by {the} one-dimensional Hopf-Cole transformation
\begin{equation}\label{ee.26}
u =- 2\omega\frac{{\phi}_x}{{\phi}}
\end{equation}

Applying the Fourier transformation to Eq. (\ref{ee.23}), we will obtain Fourier coefficient ${C_\alpha }$
\begin{equation}\label{}
  {C_\alpha }={A_\alpha }{B_\alpha }
\end{equation}
where

\begin{equation}\label{}
{A_\alpha } = \left\{ \begin{array}{l}
1,{\kern 1pt} {\kern 1pt} {\kern 1pt} {\kern 1pt} {\kern 1pt} {\kern 1pt} if{\kern 1pt} {\kern 1pt} {\kern 1pt} \alpha = 0\\
2,{\kern 1pt} {\kern 1pt} {\kern 1pt} {\kern 1pt} {\kern 1pt} {\kern 1pt} if{\kern 1pt} {\kern 1pt} {\kern 1pt} \alpha \ne 0
\end{array} \right.
\end{equation}

\begin{equation}\label{ee.29}
\begin{array}{l}
{B_\alpha } = \int_0^1 {\exp [ - \frac{1}{{2\omega }}\int_0^x {{u_0}(\zeta )d\zeta ]} } dx\\
{\kern 1pt} {\kern 1pt} {\kern 1pt} {\kern 1pt} {\kern 1pt} {\kern 1pt} {\kern 1pt} {\kern 1pt} {\kern 1pt} {\kern 1pt} {\kern 1pt} {\kern 1pt} {\kern 1pt} {\kern 1pt} {\kern 1pt} {\kern 1pt}  = \exp ( - \frac{1}{{2\omega \pi }})\int_0^1 {\exp (\frac{{\cos \pi x}}{{2\omega \pi }})} \cos (\alpha \pi x)dx
\end{array}
\end{equation}

The challenge of the initial value problem is to calculate the coefficient $B_\alpha$ of the Eq. (\ref{ee.29}), which is difficult for the single integral consisting of the exponential and trigonometric function.
To simplify {the} one-dimensional problem, the main work is to convert the calculations of $B_\alpha$ into more efficient kind. $B_\alpha$ can be written as
\begin{equation}\label{ee.30}
{B_\alpha } = \left\{ \begin{array}{l}
0,{\kern 1pt} {\kern 1pt} {\kern 1pt} {\kern 1pt} {\kern 1pt} {\kern 1pt} {\kern 1pt} {\kern 1pt} {\kern 1pt} {\kern 1pt} {\kern 1pt} {\kern 1pt} {\kern 1pt} {\kern 1pt} {\kern 1pt} {\kern 1pt} {\kern 1pt} {\kern 1pt} {\kern 1pt} {\kern 1pt} {\kern 1pt} {\kern 1pt} {\kern 1pt} {\kern 1pt} {\kern 1pt} {\kern 1pt} {\kern 1pt} {\kern 1pt} {\kern 1pt} {\kern 1pt} {\kern 1pt} {\kern 1pt} {\kern 1pt} {\kern 1pt} {\kern 1pt} {\kern 1pt} {\kern 1pt} {\kern 1pt} {\kern 1pt} {\kern 1pt} {\kern 1pt} {\kern 1pt} {\kern 1pt} {\kern 1pt} {\kern 1pt} {\kern 1pt} {\kern 1pt} {\kern 1pt} {\kern 1pt} if{\kern 1pt} {\kern 1pt} {\kern 1pt} \alpha {\kern 1pt} {\kern 1pt} {\kern 1pt} {\kern 1pt} {\kern 1pt} {\kern 1pt} is{\kern 1pt} {\kern 1pt} {\kern 1pt} {\kern 1pt} {\kern 1pt} odd\\
{I_{(\alpha )/2}}({\textstyle{1 \over {2\omega \pi }}}),{\kern 1pt} {\kern 1pt} {\kern 1pt} {\kern 1pt} {\kern 1pt} {\kern 1pt} if{\kern 1pt} {\kern 1pt} {\kern 1pt} \alpha {\kern 1pt} {\kern 1pt} {\kern 1pt} is{\kern 1pt} {\kern 1pt} {\kern 1pt} even
\end{array} \right.
\end{equation}
where $I_{(n)}(1/2\omega \pi )$ is the modified Bessel function of the first kind and of order $n$.

\subsubsection{Two-dimensional modification}\label{2.2.2}

For the two-dimensional Burgers' equation (\ref{equ.6}), we select the modified Bessel function of the first kind and of order $n$ to replace Fourier coefficient under the initial condition. Considering the following initial and boundary conditions \cite{Gao2016,Zhang2009,Zhang2010,Siraj-ul-Islam2012}
\begin{equation}\label{e.70}
\begin{array}{*{20}{l}}
\begin{array}{l}
u(x,y,0) = {u_0}(x,y) = \sin \pi x\cos \pi y\\
v(x,y,0) = {v_0}(x,y) = \cos \pi x\sin \pi y
\end{array}\\
{u(0,y,t) = u(1,y,t) = v(x,0,t) = v(x,1,t) = 0}
\end{array}
\end{equation}
where the space domain is $(x,{\kern 1pt} {\kern 1pt} y) \in \Omega  = [0,1] \times [0,1]$, and the time domain is $t > 0$.

It is widely noted that the analytical solution of {the} two-dimensional heat conduction equation can be written in the standard form of the Fourier series
\begin{equation}\label{ee.22}
\phi (x,y,t) = \sum\limits_{\alpha ,\beta  = 0}^\infty  {{C_{\alpha \beta }}\exp [ - ({\alpha ^2} + {\beta ^2}){\pi ^2}\omega t]\cos (\alpha \pi x)} \cos (\beta \pi y)
\end{equation}
where ${C_{\alpha\beta} }$ is Fourier coefficient.

The initial conditions of {the} two-dimensional heat conduction equation are extracted from the Eq. (\ref{ee.22})
\begin{equation}\label{ee.33}
\phi (x,y,0) = \sum\limits_{\alpha ,\beta  = 0}^\infty  {{C_{\alpha \beta }}\cos (\alpha \pi x)} \cos (\beta \pi y),{\kern 1pt} {\kern 1pt} {\kern 1pt}{x} \in \Omega  = [0,1]
\end{equation}
and the boundary conditions
\begin{equation}\label{e.72}
  {\phi _x}(0,y,t) = {\phi _x}(1,y,t) = {\phi _y}(x,0,t)= {\phi _y}(x,1,t)=0
\end{equation}
by {the} two-dimensional Hopf-Cole transformation
\begin{equation}\label{}
u =- 2\omega\frac{{\phi}_x}{{\phi}}{\kern 1pt} {\kern 1pt},{\kern 1pt} {\kern 1pt}v =- 2\omega\frac{{\phi}_y}{{\phi}}
\end{equation}

Applying the Fourier transformation to Eq. (\ref{ee.33}), we will obtain Fourier coefficient ${C_{\alpha\beta}}$
\begin{equation}\label{}
  {C_{\alpha\beta} }={A_{\alpha\beta} }{B_{\alpha\beta} }
\end{equation}
where
\begin{equation}\label{}
{A_{\alpha\beta} } = \left\{ \begin{array}{l}
1,{\kern 1pt} {\kern 1pt} {\kern 1pt} {\kern 1pt} {\kern 1pt} {\kern 1pt} if{\kern 1pt} {\kern 1pt} {\kern 1pt} \alpha = 0{\kern 1pt} {\kern 1pt} {\kern 1pt} {\kern 1pt} {\kern 1pt} {\kern 1pt} and{\kern 1pt} {\kern 1pt} {\kern 1pt} \beta = 0\\
2,{\kern 1pt} {\kern 1pt} {\kern 1pt} {\kern 1pt} {\kern 1pt} {\kern 1pt} if{\kern 1pt} {\kern 1pt} {\kern 1pt} \alpha = 0{\kern 1pt} {\kern 1pt} {\kern 1pt} {\kern 1pt} {\kern 1pt} {\kern 1pt} and{\kern 1pt} {\kern 1pt} {\kern 1pt}\beta \ne 0\\
2,{\kern 1pt} {\kern 1pt} {\kern 1pt} {\kern 1pt} {\kern 1pt} {\kern 1pt} if{\kern 1pt} {\kern 1pt} {\kern 1pt} \alpha \neq 0{\kern 1pt} {\kern 1pt} {\kern 1pt} {\kern 1pt} {\kern 1pt} {\kern 1pt} and{\kern 1pt} {\kern 1pt} {\kern 1pt} \beta = 0\\
4,{\kern 1pt} {\kern 1pt} {\kern 1pt} {\kern 1pt} {\kern 1pt} {\kern 1pt} if{\kern 1pt} {\kern 1pt} {\kern 1pt} \alpha \ne 0{\kern 1pt} {\kern 1pt} {\kern 1pt} {\kern 1pt} {\kern 1pt} {\kern 1pt} and{\kern 1pt} {\kern 1pt} {\kern 1pt} \beta  \neq 0\\

\end{array} \right.
\end{equation}

\begin{equation}\label{ee.39}
\begin{array}{l}
{B_{\alpha \beta }} = \int_0^1 {\int_0^1 {\exp [ - \frac{1}{{2\omega }}D(x,y)]} dxdy} \\
{\kern 1pt} {\kern 1pt} {\kern 1pt} {\kern 1pt} {\kern 1pt} {\kern 1pt} {\kern 1pt} {\kern 1pt} {\kern 1pt} {\kern 1pt} {\kern 1pt} {\kern 1pt} {\kern 1pt} {\kern 1pt} {\kern 1pt} {\kern 1pt} {\kern 1pt} {\kern 1pt} {\kern 1pt}  = \exp ( - \frac{1}{{2\omega \pi }})\int_0^1 {\int_0^1 {\exp (\frac{{\cos \pi x\cos \pi y}}{{2\omega \pi }})} } \cos (\alpha \pi x)\cos (\beta \pi y)dxdy
\end{array}
\end{equation}
where
\begin{equation}\label{}
  D(x,y) = \frac{1}{2}\{ \int_0^x {[{u_0}(\zeta ,y) + {u_0}(\zeta ,0)]} d\zeta  + \int_0^y {[{v_0}(x,\zeta ) + {v_0}(0,\zeta )]} d\zeta \}
\end{equation}

The challenge of the initial value problem is to calculate the coefficient $B_{\alpha\beta}$ of the Eq. (\ref{ee.39}), which is difficult for the double integral consisting of the exponential and trigonometric function.
To simplify {the} two-dimensional problem, the main work is to convert the calculations of $B_{\alpha\beta}$ into more efficient kind. $B_{\alpha\beta}$ can be written as

\begin{equation}\label{79}
{B_{\alpha \beta }} = \left\{ \begin{array}{l}
0,{\kern 1pt} {\kern 1pt} {\kern 1pt} {\kern 1pt} {\kern 1pt} {\kern 1pt} {\kern 1pt} {\kern 1pt} {\kern 1pt} {\kern 1pt} {\kern 1pt} {\kern 1pt} {\kern 1pt} {\kern 1pt} {\kern 1pt} {\kern 1pt} {\kern 1pt} {\kern 1pt} {\kern 1pt} {\kern 1pt} {\kern 1pt} {\kern 1pt}{\kern 1pt} {\kern 1pt} {\kern 1pt} {\kern 1pt} {\kern 1pt}  {\kern 1pt} {\kern 1pt} {\kern 1pt} {\kern 1pt} {\kern 1pt} {\kern 1pt} {\kern 1pt} {\kern 1pt} {\kern 1pt} {\kern 1pt} {\kern 1pt} {\kern 1pt} {\kern 1pt}{\kern 1pt} {\kern 1pt} {\kern 1pt} {\kern 1pt} {\kern 1pt} {\kern 1pt} {\kern 1pt} {\kern 1pt} {\kern 1pt} {\kern 1pt} {\kern 1pt} {\kern 1pt} {\kern 1pt} {\kern 1pt} {\kern 1pt} {\kern 1pt} {\kern 1pt} {\kern 1pt} {\kern 1pt} {\kern 1pt} {\kern 1pt} {\kern 1pt} {\kern 1pt} {\kern 1pt} {\kern 1pt} {\kern 1pt} {\kern 1pt} {\kern 1pt} {\kern 1pt} {\kern 1pt} {\kern 1pt} {\kern 1pt} {\kern 1pt} {\kern 1pt} {\kern 1pt} {\kern 1pt} {\kern 1pt} {\kern 1pt} {\kern 1pt} {\kern 1pt} {\kern 1pt} {\kern 1pt} {\kern 1pt} {\kern 1pt} {\kern 1pt} {\kern 1pt} {\kern 1pt} {\kern 1pt} {\kern 1pt} {\kern 1pt} {\kern 1pt} {\kern 1pt} {\kern 1pt} {\kern 1pt} {\kern 1pt} {\kern 1pt} {\kern 1pt} {\kern 1pt} {\kern 1pt} {\kern 1pt} {\kern 1pt} {\kern 1pt} {\kern 1pt} {\kern 1pt} {\kern 1pt} {\kern 1pt} {\kern 1pt} {\kern 1pt} {\kern 1pt} {\kern 1pt} {\kern 1pt} {\kern 1pt} {\kern 1pt} {\kern 1pt} {\kern 1pt} {\kern 1pt} {\kern 1pt} {\kern 1pt} {\kern 1pt} {\kern 1pt} {\kern 1pt} {\kern 1pt} {\kern 1pt} if{\kern 1pt} {\kern 1pt} {\kern 1pt} \alpha  + \beta {\kern 1pt} {\kern 1pt} {\kern 1pt} {\kern 1pt} {\kern 1pt} {\kern 1pt} is{\kern 1pt} {\kern 1pt} {\kern 1pt} {\kern 1pt} {\kern 1pt} odd\\
{I_{(\alpha  + \beta )/2}}({\textstyle{1 \over {4\omega \pi }}}){I_{(\alpha  - \beta )/2}}({\textstyle{1 \over {4\omega \pi }}}),{\kern 1pt} {\kern 1pt} {\kern 1pt} {\kern 1pt} {\kern 1pt} {\kern 1pt} if{\kern 1pt} {\kern 1pt} {\kern 1pt} \alpha  + \beta {\kern 1pt} {\kern 1pt} {\kern 1pt} is{\kern 1pt} {\kern 1pt} {\kern 1pt} even
\end{array} \right.
\end{equation}
where $I_{(n)}(\textstyle\frac{1}{4\omega \pi})$ is the modified Bessel function of the first kind and of order $n$. Gao \cite{Gao2016} proved that Eq. (\ref{ee.39}) and Eq. (\ref{79}) are equal, the scheme will get more precise and efficient analytical solution through the modification.

\subsubsection{Three-dimensional modification}\label{2.2.3}
Because the three-dimensional Burgers' equation is too complex, its improvement is quite different from the one-dimensional and two-dimensional ones. Considering the following initial and boundary conditions \cite{Gao2016}
\begin{equation}\label{80}
\begin{array}{l}
u(x,y,z,0) = \sin \pi x{\cos \pi y}{\cos \pi z}\\
v(x,y,z,0) = \sin \pi x{\cos \pi y}{\cos \pi z}\\
{w(x,y,z,0) = \sin \pi x{\cos \pi y}{\cos \pi z}}

\end{array}
\end{equation}

\begin{equation}\label{81}
\begin{array}{l}
u(0,y,z,t) = u(1,y,t)=0 \\
v(x,0,z,t) = v(x,1,t) = 0\\
w(x,y,0,t) = w(x,1,t) = 0

\end{array}
\end{equation}
where the space domain is $(x,y,z) \in \Omega  = [0,1] \times [0,1]\times [0,1]$,and the time domain is $t > 0$.

It is widely noted that the analytical solution of {the} three-dimensional heat conduction equation can be written in the standard form of the Fourier series
\begin{equation}\label{ee.43}
\begin{array}{l}
\phi (x,y,z,t) = \\
\sum\limits_{\alpha ,\beta ,\gamma  = 0}^\infty  {{C_{\alpha \beta \gamma }}\exp ( - ({\alpha ^2} + {\beta ^2} + {\gamma ^2}){\pi ^2}\omega t)\cos (\alpha \pi x)} \cos (\beta \pi y)\cos (\gamma \pi z)
\end{array}
\end{equation}
where ${C_{\alpha\beta\gamma} }$ is Fourier coefficient.

The initial conditions of {the} three-dimensional heat conduction equation are extracted from the Eq. (\ref{ee.43})
\begin{equation}\label{ee.44}
\phi (x,y,z,0) = \sum\limits_{\alpha ,\beta ,\gamma  = 0}^\infty  {{C_{\alpha \beta \gamma }}\cos (\alpha \pi x)} \cos (\beta \pi y)\cos (\gamma \pi z)
\end{equation}
and the boundary conditions
\begin{equation}\label{83}
\begin{array}{l}
  {\phi _x}(0,y,z,t) = {\phi _x}(1,y,z,t) =0\\
   {\phi _y}(x,0,z,t)= {\phi _y}(x,1,z,t)=0\\
   {\phi _z}(x,y,0,t)= {\phi _z}(x,y,1,t)=0
  \end{array}
\end{equation}
by {the} three-dimensional Hopf-Cole transformation
\begin{equation}\label{e.84}
u =- 2\omega\frac{{\phi}_x}{{\phi}}{\kern 1pt} {\kern 1pt},{\kern 1pt} {\kern 1pt}v =- 2\omega\frac{{\phi}_y}{{\phi}},{\kern 1pt} {\kern 1pt} {w =- 2\omega\frac{{\phi}_z}{{\phi}}}
\end{equation}

Applying the Fourier transformation to Eq. (\ref{ee.44}), we will obtain Fourier coefficient ${C_{\alpha\beta\gamma}}$
\begin{equation}\label{}
  {C_{\alpha\beta\gamma} }={A_{\alpha\beta\gamma} }{B_{\alpha\beta\gamma} }
\end{equation}
where

\begin{equation}\label{}
{A_{\alpha\beta\gamma} } = \left\{ \begin{array}{l}
1,{\kern 1pt} {\kern 1pt} {\kern 1pt} {\kern 1pt} {\kern 1pt} {\kern 1pt} if{\kern 1pt} {\kern 1pt} {\kern 1pt} \alpha = 0{\kern 1pt} {\kern 1pt} {\kern 1pt} {\kern 1pt} {\kern 1pt} {\kern 1pt} and{\kern 1pt} {\kern 1pt} {\kern 1pt} \beta = 0{\kern 1pt} {\kern 1pt} {\kern 1pt} {\kern 1pt} {\kern 1pt} {\kern 1pt} and{\kern 1pt} {\kern 1pt} {\kern 1pt} \gamma = 0\\
2,{\kern 1pt} {\kern 1pt} {\kern 1pt} {\kern 1pt} {\kern 1pt} {\kern 1pt} if{\kern 1pt} {\kern 1pt} {\kern 1pt} \alpha \neq 0{\kern 1pt} {\kern 1pt} {\kern 1pt} {\kern 1pt} {\kern 1pt} {\kern 1pt} and{\kern 1pt} {\kern 1pt} {\kern 1pt}\beta = 0{\kern 1pt} {\kern 1pt} {\kern 1pt} {\kern 1pt} {\kern 1pt} {\kern 1pt} and{\kern 1pt} {\kern 1pt} {\kern 1pt} \gamma = 0\\
2,{\kern 1pt} {\kern 1pt} {\kern 1pt} {\kern 1pt} {\kern 1pt} {\kern 1pt} if{\kern 1pt} {\kern 1pt} {\kern 1pt} \alpha = 0{\kern 1pt} {\kern 1pt} {\kern 1pt} {\kern 1pt} {\kern 1pt} {\kern 1pt} and{\kern 1pt} {\kern 1pt} {\kern 1pt} \beta \neq 0{\kern 1pt} {\kern 1pt} {\kern 1pt} {\kern 1pt} {\kern 1pt} {\kern 1pt} and{\kern 1pt} {\kern 1pt} {\kern 1pt} \gamma = 0\\
2,{\kern 1pt} {\kern 1pt} {\kern 1pt} {\kern 1pt} {\kern 1pt} {\kern 1pt} if{\kern 1pt} {\kern 1pt} {\kern 1pt} \alpha = 0{\kern 1pt} {\kern 1pt} {\kern 1pt} {\kern 1pt} {\kern 1pt} {\kern 1pt} and{\kern 1pt} {\kern 1pt} {\kern 1pt} \beta = 0{\kern 1pt} {\kern 1pt} {\kern 1pt} {\kern 1pt} {\kern 1pt} {\kern 1pt} and{\kern 1pt} {\kern 1pt} {\kern 1pt} \gamma \neq 0\\
4,{\kern 1pt} {\kern 1pt} {\kern 1pt} {\kern 1pt} {\kern 1pt} {\kern 1pt} if{\kern 1pt} {\kern 1pt} {\kern 1pt} \alpha \ne 0{\kern 1pt} {\kern 1pt} {\kern 1pt} {\kern 1pt} {\kern 1pt} {\kern 1pt} and{\kern 1pt} {\kern 1pt} {\kern 1pt} \beta\neq 0{\kern 1pt} {\kern 1pt} {\kern 1pt} {\kern 1pt} {\kern 1pt} {\kern 1pt} and{\kern 1pt} {\kern 1pt} {\kern 1pt} \gamma = 0\\
4,{\kern 1pt} {\kern 1pt} {\kern 1pt} {\kern 1pt} {\kern 1pt} {\kern 1pt} if{\kern 1pt} {\kern 1pt} {\kern 1pt} \alpha = 0{\kern 1pt} {\kern 1pt} {\kern 1pt} {\kern 1pt} {\kern 1pt} {\kern 1pt} and{\kern 1pt} {\kern 1pt} {\kern 1pt} \beta\neq 0{\kern 1pt} {\kern 1pt} {\kern 1pt} {\kern 1pt} {\kern 1pt} {\kern 1pt} and{\kern 1pt} {\kern 1pt} {\kern 1pt} \gamma\neq 0\\
4,{\kern 1pt} {\kern 1pt} {\kern 1pt} {\kern 1pt} {\kern 1pt} {\kern 1pt} if{\kern 1pt} {\kern 1pt} {\kern 1pt} \alpha \ne 0{\kern 1pt} {\kern 1pt} {\kern 1pt} {\kern 1pt} {\kern 1pt} {\kern 1pt} and{\kern 1pt} {\kern 1pt} {\kern 1pt} \beta= 0{\kern 1pt} {\kern 1pt} {\kern 1pt} {\kern 1pt} {\kern 1pt} {\kern 1pt} and{\kern 1pt} {\kern 1pt} {\kern 1pt} \gamma \neq 0\\
8,{\kern 1pt} {\kern 1pt} {\kern 1pt} {\kern 1pt} {\kern 1pt} {\kern 1pt} if{\kern 1pt} {\kern 1pt} {\kern 1pt} \alpha \ne 0{\kern 1pt} {\kern 1pt} {\kern 1pt} {\kern 1pt} {\kern 1pt} {\kern 1pt} and{\kern 1pt} {\kern 1pt} {\kern 1pt} \beta\neq 0{\kern 1pt} {\kern 1pt} {\kern 1pt} {\kern 1pt} {\kern 1pt} {\kern 1pt} and{\kern 1pt} {\kern 1pt} {\kern 1pt} \gamma \neq 0\\
\end{array} \right.
\end{equation}

\begin{equation}\label{ee.49}
\begin{array}{*{20}{l}}
{{B_{\alpha \beta \gamma }} = \int_0^1 {\int_0^1 {\int_0^1 {\exp [ - \frac{1}{{2\omega }}D(x,y,z)]} dxdydz} } }\\
{ = \exp ( - \frac{1}{{2\omega \pi }})\int_0^1 {\int_0^1 {\int_0^1 {\exp (\frac{{\cos \pi x\cos \pi y\cos \pi z}}{{2\omega \pi }})} } } \cos (\alpha \pi x)\cos (\beta \pi y)\cos (\gamma \pi z)dxdydz}
\end{array}
\end{equation}
where
\begin{equation}\label{}
\begin{array}{l}
D(x,y,z) = \\
\frac{1}{3}\left( \begin{array}{l}
\int_0^x {[{u_0}(\zeta ,y,z) + {u_0}(\zeta ,0,z) + {u_0}(\zeta ,0,0)]} d\zeta \\
 + \int_0^y {[{v_0}(x,\zeta ,z) + {v_0}(x,\zeta ,0) + {v_0}(0,\zeta ,0)]} d\zeta \\
 + \int_0^z {[{w_0}(x,y,\zeta ) + {w_0}(0,y,\zeta ) + {w_0}(0,0,\zeta )]} d\zeta
\end{array} \right)
\end{array}
\end{equation}

The challenge of the initial value problem is to calculate the coefficient $B_{\alpha\beta\gamma}$ of the Eq. (\ref{ee.49}), which is difficult for that the triple integral consisting of the exponential and trigonometric function.
To simplify {the} three-dimensional problem, the main job is to convert the calculations of $B_{\alpha\beta\gamma}$ into more efficient kind. $B_{\alpha\beta\gamma}$ can be written as
\begin{equation}\label{90}
 \begin{array}{l}
{B_{\alpha \beta \gamma }} = \frac{{{{(1/8\omega \pi )}^\beta }}}{{[(\alpha  + \beta )/2]![(\alpha  - \beta )/2]!}}\\
 \times \left\{ \begin{array}{l}
\sum\limits_{j = 1}^{(\gamma  + 1)/2} {{\mu _j}\frac{{(\beta  + 2j - 2)!!}}{{(\beta  + 2j - 1)!!}}G,{\kern 1pt} {\kern 1pt} {\kern 1pt} {\kern 1pt} {\kern 1pt} {\rm{if}}{\kern 1pt} {\kern 1pt} \alpha ,\beta ,\gamma {\kern 1pt} {\kern 1pt} {\kern 1pt} {\rm{are}}{\kern 1pt} {\kern 1pt} {\kern 1pt} {\kern 1pt} {\rm{all}}{\kern 1pt} {\kern 1pt} {\kern 1pt} {\kern 1pt} {\rm{odd}}} \\
\sum\limits_{j = 0}^{(\gamma )/2} {{\mu _j}\frac{{(\beta  + 2j - 2)!!}}{{(\beta  + 2j - 1)!!}}G,{\kern 1pt} {\kern 1pt} {\kern 1pt} {\kern 1pt} {\kern 1pt} {\kern 1pt} {\kern 1pt} {\kern 1pt} {\kern 1pt} {\kern 1pt} {\rm{if}}{\kern 1pt} {\kern 1pt} {\kern 1pt} \alpha ,\beta ,\gamma {\kern 1pt} {\kern 1pt} {\kern 1pt} {\rm{are}}{\kern 1pt} {\kern 1pt} {\kern 1pt} {\kern 1pt} {\rm{all}}{\kern 1pt} {\kern 1pt} {\kern 1pt} {\kern 1pt} {\rm{even}}} \\
0,{\kern 1pt} {\kern 1pt} {\kern 1pt} {\kern 1pt} {\kern 1pt} otherwise
\end{array} \right.
\end{array}
\end{equation}
where $G={}_3{F_4}({\textstyle{{\beta  + 1} \over 2}},{\textstyle{\beta  \over 2}} + 1,{\textstyle{\beta  \over 2}} + j;{\textstyle{{\beta  + \alpha } \over 2}} + 1,{\textstyle{{\beta  - \alpha } \over 2}} + 1,m + 1,{\textstyle{{\beta  - \alpha } \over 2}} + j + 1;{({\textstyle{1 \over {4\omega \pi }}})^2})$ is the generalized hypergeometric series. Ref. \cite{Perger1993} defines the generalized hypergeometric series
\begin{equation}\label{}
{}_p{F_q}({\alpha _1},{\alpha _2}, \cdots ,{\alpha _p},{\beta _1},{\beta _2}, \cdots ,{\beta _p};\varpi ) = \sum\limits_{s = 0}^\infty  {\frac{{{{({\alpha _1})}_s}{{({\alpha _2})}_s} \cdots {{({\alpha _p})}_s}}}{{{{({\alpha _s})}_s}{{({\alpha _2})}_s} \cdots {{({\alpha _p})}_s}}}} \frac{{{\varpi ^s}}}{{s!}}
\end{equation}
in which $( a )_k$ is the Pochhammer symbol and is defined as:
\begin{equation}\label{}
\begin{array}{l}
{\alpha _0} = 1\\
{(\alpha )_k} = \alpha (\alpha  + 1)(\alpha  + 2) \cdots (\alpha  + l - 1),k \ge 1
\end{array}
\end{equation}

The coefficient $\mu_j$ in Eq. (\ref{90}) is defined by the following equation:
\begin{equation}\label{}
\cos (\gamma \varpi ) = \left\{ \begin{array}{l}
\sum\limits_{j = 1}^{\gamma /2} {{\mu _j}co{s^{2j - 1}}(\varpi ),{\kern 8pt} {\kern 1pt} {\kern 1pt} {\kern 1pt} {\kern 1pt} {\rm{if}}{\kern 1pt} {\kern 1pt} \gamma {\kern 1pt} {\kern 1pt} {\kern 1pt} {\rm{are}}{\kern 1pt} {\kern 1pt} {\kern 1pt} {\kern 1pt} {\rm{even}}} \\
\sum\limits_{j = 1}^{(\gamma  + 1)/2} {{\mu _j}co{s^{2j - 1}}(\varpi ),{\kern 1pt} {\kern 1pt} {\kern 1pt} {\kern 1pt} {\kern 1pt} {\rm{if}}{\kern 1pt} {\kern 1pt} \gamma {\kern 1pt} {\kern 1pt} {\kern 1pt} {\rm{are}}{\kern 1pt} {\kern 1pt} {\kern 1pt} odd}
\end{array} \right.
\end{equation}

\section{High-order numerical scheme}
 In this section, to solve the n-dimensional heat conduction equation obtained by Hopf-Cole transformation, we will present the sixth-order CFD scheme and the precise integration method (PIM).
 For simplicity, we consider one-dimensional heat conduction equation (\ref{e.10}) with mesh size ${h=x_{i+1}-x_{i}}$, in which $x_{i}=ih,i=1,2,...,N$, where $h$ is spatial step size.
 We firstly apply the sixth-order CFD scheme to discretization in space. If $\phi=\phi(x_{i})$ and $\phi{}''_{i}$ represent the second derivative of $\phi(x)$ at $x_{i}$, then an approximation of the second derivatives at interior nodes may be expressed as
\begin{equation}\label{ee.56}
\textstyle{\frac{2}{11}\phi{}''_{i-1}+\phi{}''_{i}+\frac{2}{11}\phi{}''_{i}}=\textstyle{\frac{3}{44h^{2}}({\phi_{i + 2}} - {\phi_i} + {\phi_{i - 2}})+\frac{12}{11h^{2}}({\phi_{i + 1}} - 2{\phi_i} + {\phi_{i - 1}})}
\end{equation}

In order to make those near-boundary points have the same order accuracy as interior nodes, they should be obtained by matching Taylor series expansions to the order of $O(h^6)$ at boundary points $1, 2, N-1$ and $N$, hence we get the following formulae \cite{Li2008}
\begin{equation}\label{ee.57}
\begin{array}{l}
{\kern 10pt}{\phi''_{1}} + {\textstyle{{126} \over {11}}}{\phi''_{2}}\\
 = {\textstyle{1 \over {{h^2}}}}\left( {{\textstyle{{2077} \over {157}}}{\phi_1} - {\textstyle{{2943} \over {110}}}{\phi_2} + {\textstyle{{573} \over {44}}}{\phi_3} + {\textstyle{{167} \over {99}}}{\phi_4} - {\textstyle{{18} \over {11}}}{\phi_5} + {\textstyle{{57} \over {110}}}{\phi_6} - {\textstyle{{131} \over {1980}}}{\phi_7}} \right)
\end{array}
\end{equation}
\begin{equation}\label{ee.58}
\begin{array}{l}
{\kern 10pt}{\textstyle{{11} \over {128}}}{\phi''_{1}} + {\phi''_{2}} + {\textstyle{{11} \over {128}}}{\phi''_{3}}\\
 = \frac{1}{{{h^2}}}\left( {{\textstyle{{585} \over {512}}}{\phi_1} - {\textstyle{{141} \over {64}}}{\phi_2} + {\textstyle{{459} \over {512}}}{\phi_3} + {\textstyle{9 \over {32}}}{\phi_4} - {\textstyle{{81} \over {512}}}{\phi_5} + {\textstyle{3 \over {64}}}{\phi_6} - {\textstyle{3 \over {512}}}{\phi_7}} \right)
\end{array}
\end{equation}

\begin{equation}\label{ee.59}
 \begin{array}{l}
{\kern 10pt}{\textstyle{{11} \over {128}}}{\phi''_{N}} + {\phi''_{N-1}} + {\textstyle{{11} \over {128}}}{\phi''_{N-2}}\\
 = {\textstyle{1 \over {{h^2}}}}\left( \begin{array}{l}
{\textstyle{{585} \over {512}}}{\phi_N} - {\textstyle{{141} \over {64}}}{\phi_{N - 1}} + {\textstyle{{459} \over {512}}}{\phi_{N - 2}} + \\
{\textstyle{9 \over {32}}}{\phi_{N - 3}} - {\textstyle{{81} \over {512}}}{\phi_{N - 4}} + {\textstyle{3 \over {64}}}{\phi_{N - 5}} - {\textstyle{3 \over {512}}}{\phi_{N - 6}}
\end{array} \right)
\end{array}
\end{equation}

\begin{equation}\label{ee.60}
\begin{array}{l}
{\kern 10pt} {\textstyle{{126} \over {11}}}{\phi''_{N-1}} +{\phi''_{N}}\\
 = \frac{1}{{{h^2}}}\left( \begin{array}{l}
{\textstyle{{2077} \over {157}}}{\phi_N} - {\textstyle{{2943} \over {110}}}{\phi_{N - 1}} + {\textstyle{{573} \over {44}}}{\phi_{N - 2}} + \\
{\textstyle{{167} \over {99}}}{\phi_{N - 3}} - {\textstyle{{18} \over {11}}}{\phi_{N - 4}} + {\textstyle{{57} \over {110}}}{\phi_{N - 5}} - {\textstyle{{131} \over {1980}}}{\phi_{N - 6}}
\end{array} \right)
\end{array}
\end{equation}

Writing Eqs. (\ref{ee.56},\ref{ee.57},\ref{ee.58},\ref{ee.59},\ref{ee.60}) in matrix form as
\begin{equation}\label{}
{\bm{A\Phi''}} = {\bm{B\Phi}}
\end{equation}
where
\begin{equation}\label{}
  {\bm{A}} = {\left[ {\begin{array}{*{20}{c}}
1&{{\textstyle{{126} \over {11}}}}&{}&{}&{}&{}&{}\\
{{\textstyle{{11} \over {128}}}}&1&{{\textstyle{{11} \over {128}}}}&{}&{}&{}&{}\\
{}&{{\textstyle{2 \over {11}}}}&1&{{\textstyle{2 \over {11}}}}&{}&{}&{}\\
{}&{}& \ddots & \ddots & \ddots &{}&{}\\
{}&{}&{}&{{\textstyle{2 \over {11}}}}&1&{{\textstyle{2 \over {11}}}}&{}\\
{}&{}&{}&{}&{{\textstyle{{11} \over {128}}}}&1&{{\textstyle{{11} \over {128}}}}\\
{}&{}&{}&{}&{}&{{\textstyle{{126} \over {11}}}}&1
\end{array}} \right]_{N \times N}}
\end{equation}
\begin{equation}\label{}
  {\bm{B}} = \frac{1}{{{h^2}}}{\left[ {\begin{array}{*{20}{c}}
{{\textstyle{{2077} \over {157}}}}&{{\textstyle{{ - 2943} \over {110}}}}&{{\textstyle{{574} \over {44}}}}&{{\textstyle{{167} \over {99}}}}&{{\textstyle{{ - 18} \over {11}}}}&{{\textstyle{{57} \over {110}}}}&{{\textstyle{{ - 131} \over {1980}}}}&{}\\
{{\textstyle{{585} \over {512}}}}&{{\textstyle{{ - 141} \over {64}}}}&{{\textstyle{{459} \over {512}}}}&{{\textstyle{9 \over {32}}}}&{{\textstyle{{ - 81} \over {512}}}}&{{\textstyle{3 \over {64}}}}&{{\textstyle{{ - 3} \over {512}}}}&{}\\
{{\textstyle{3 \over {44}}}}&{{\textstyle{{12} \over {11}}}}&{{\textstyle{{ - 51} \over {22}}}}&{{\textstyle{{12} \over {11}}}}&{{\textstyle{3 \over {44}}}}&{}&{}&{}\\
{}& \ddots & \ddots & \ddots & \ddots & \ddots &{}&{}\\
{}&{}&{{\textstyle{3 \over {44}}}}&{{\textstyle{{12} \over {11}}}}&{{\textstyle{{ - 51} \over {22}}}}&{{\textstyle{{12} \over {11}}}}&{{\textstyle{3 \over {44}}}}&{}\\
{}&{{\textstyle{{ - 3} \over {512}}}}&{{\textstyle{3 \over {64}}}}&{{\textstyle{{ - 81} \over {512}}}}&{{\textstyle{9 \over {32}}}}&{{\textstyle{{459} \over {512}}}}&{{\textstyle{{ - 141} \over {64}}}}&{{\textstyle{{585} \over {512}}}}\\
{}&{{\textstyle{{ - 131} \over {1980}}}}&{{\textstyle{{57} \over {110}}}}&{{\textstyle{{ - 18} \over {11}}}}&{{\textstyle{{167} \over {99}}}}&{{\textstyle{{574} \over {44}}}}&{{\textstyle{{ - 2943} \over {110}}}}&{{\textstyle{{2077} \over {157}}}}
\end{array}} \right]_{N \times N}}
\end{equation}
\begin{equation}\label{}
  {\bm{\Phi}} = {\left( {{\phi_1},{\phi_2}, \cdots ,{\phi_{N - 1}},{\phi_N}} \right)^T}
\end{equation}
\par Therefore the sixth-order compact finite difference approximation of second derivatives ${\bm{\Phi''}}$ is given by
\begin{equation}\label{equ.64}
 {\bm{\Phi''}}={\bm{A^{-1}B\Phi}}={\bm{H\Phi}}
\end{equation}
where ${\bm{H}}={\bm{A^{-1}B}}$.
\subsection{ Precise integration method}\label{3.1}
After the spatial discretization, the governing PDEs become the following ODEs
\begin{equation}\label{ee.66}
  \frac{\textcolor{green}{\rm{d}}{\bm{\Phi}}}{\textcolor{green}{\rm{d}}t}={\bm{H\Phi}}
\end{equation}
\par Giving $\tau {\rm{ = }}{t_{k + 1}} - {t_k}$ as the temporal step size, then integrating Eq. (\ref{ee.66}) directly, the following recurrence formula is obtained
\begin{equation}\label{}
  {{\bm{\Phi}}^{k + 1}} = {e^{{\bm{H}}\tau }}{{\bm{\Phi}}^k}{\rm{ = }}{\bm{T}}\left( \tau  \right){{\bm{\Phi}}^k}
\end{equation}
where ${\bm{T}}\left( \tau  \right) = {e^{{\bm{H}}\tau }}$ is an exponential matrix.
\par The present work will focus on how to compute the exponential matrix ${\bm{T}}$ very precisely.
Moler et al. \cite{Moler2003} had discussed nineteen dubious ways to compute the exponential matrix,
they pointed out that the problem of calculating exponential matrix had not been fully solved.
In this paper, we apply the PIM to calculate the exponential matrix, which was proposed by Zhong et al. \cite{Wan-Xie2004}.
The PIM is a algorithm of high precision for computing exponential matrix,
which avoids the computer truncation error caused by the fine division and improves the numerical accuracy of the exponential matrix.
In short, PIM is a series of matrix or vector multiplication calculations.
Therefore, the main problem is how to calculate the exponential matrix ${e^{{\bm{H}}\tau }}$.
The precise computation of exponential matrix has two core contents \cite{Zhong2001}:
\\(1) The incremental part of the exponential matrix is calculated separately, rather than as a whole.
\\(2) The addition theorem of exponent is achieved by $2^n$ algorithm.

Using the addition theorem for the exponential matrix ${e^{{\bm{H}}\tau }}$, the following equation is obtained

\begin{equation}\label{e.34}
  {e^{{\bm{H}}\tau }} = {\left( {{e^{{\bm{H}}\Delta t}}} \right)^m}
\end{equation}
{where $m$ is a relatively large positive integer and $\Delta t = \frac{\tau }{m}$. Thus $\Delta t$ is an extremely short time.} In order to ensure computational accuracy, Ref. \cite{Wan-Xie2004} suggested $ m = {2^n},n = 20,m = 1048576$
, $n$ is defined as bisection order.
\subsubsection{Taylor approximation methods}
The accurate computation of ${e^{{\bm{H}}\tau }}$ is a challenging problem for the numerical analysis \cite{Bhatt2016, Moler2003}.
The major issue is the cancellation error arising during the direct computation of ${e^{{\bm{H}}\tau }}$ for eigenvalues of ${e^{{\bm{H}} }}$ close to 0. To overcome this problem and other numerical issues associated with it, many researchers have proposed different methods.
This study focuses on the new technique and develops CFD scheme based on Taylor approximation of ${e^{{\bm{H}}\tau }}$ in order to alleviate computational difficulties associated with them.
The Taylor expansion formula to the exponential matrix ${e^{{\bm{H}}\Delta t}}$ is defined as
\begin{equation}\label{}
  {e^{{\bm{H}}\Delta t}}{\rm{ = }}\sum\limits_{j = 0}^\infty  {\frac{{{{({\bm{H}}\Delta t)}^j}}}{{j!}}}
\end{equation}
{ where $\Delta t$ is extremely short,} so the truncation error in time can be ignored, the fourth-order Taylor expansion can have high precision. Hence
\begin{equation}\label{}
{ {\bm{T}}\left( {\Delta t} \right)\approx \sum\limits_{j = 0}^4 {\frac{{{{({\bm{H}}\Delta t)}^j}}}{{j!}}{\rm{ = }}{\kern 1pt} {\kern 1pt} {\kern 1pt} {\bm{I}} + {\bm{H}}\Delta t + \frac{{{{({\bm{H}}\Delta t)}^2}}}{{2!}} + \frac{{{{({\bm{H}}\Delta t)}^3}}}{{3!}} + \frac{{{{({\bm{H}}\Delta t)}^4}}}{{4!}}}}
\end{equation}

Because $\tau$ is very small, it is enough to expand only the the first five terms of the series. The exponential matrix ${\bm{T}}\left( {\Delta t} \right)$ departs from the unit matrix ${\bm{I}}$ to a very small extent. Hence it should be distinguished as
\begin{equation}\label{e.37}
  {e^{{\bm{H}}\Delta t}} \cong {\bm{I}} + {{\bm{T}}_a}{\kern 1pt}  = {\bm{I}} + {\bm{H}}\Delta t + \frac{{{{({\bm{H}}\Delta t)}^2}}}{{2!}} + \frac{{{{({\bm{H}}\Delta t)}^3}}}{{3!}} + \frac{{{{({\bm{H}}\Delta t)}^4}}}{{4!}}
\end{equation}

In order to obtain exponential matrix ${\bm{T}}( \tau )$, we need to use $2^n$ algorithm for the matrix ${\bm{T}}( {\Delta t} )$.
\subsubsection{$2^n$ algorithm of the exponential matrix}
PIM has the problem of complete loss of precision in the exponential additional theorem \cite{Zhong2001,Zhang2012,Wang2009}.
One of the core contents of PIM is the identity matrix ${\bm{I}}$ cannot be directly added to the incremental matrix ${\bm{T_a}}$ for Eq. (\ref{e.37}).
Because ${\bm{T_a}}$ is a miniature matrix.
If they add up directly, ${\bm{T_a}}$ becomes the mantissa of ${\bm{I+T_a}}$ in the process of computer operation. Thus, ${\bm{T_a}}$ will become an appended part and its precision will seriously drop during the round-off operation in computer arithmetic. As a matter of fact, ${\bm{T_a}}$ is an incremental part, which must be calculated and stored separately.
Therefore, we will apply $2^n$  algorithm to calculate ${\bm{T_a}}$.

For computing the matrix ${\bm{T}}( \tau  ) = {e^{{\bm{H}}\tau }}$, Eq. (\ref{e.34}) should be factored as
\begin{equation}\label{e.38}
  {\bm{T}}( \tau  ) = {({\bm{I}} + {{\bm{T}}_a})^{{2^n}}} = {( {{\bm{I}} + {{\bm{T}}_a}} )^{{2^{n - 1}}}} \times {( {{\bm{I}} + {{\bm{T}}_a}} )^{{2^{n - 1}}}}
\end{equation}

Because ${\bm{T}}(\alpha )$ has the following equation relation of factorization

\begin{equation}\label{e.39}
  {{\bm{T}}_{a + 1}} \times {{\bm{T}}_{a + 1}} = {\bm{I}} + (2{{\bm{T}}_a} + {{\bm{{\rm T}}}_a} \times {{\bm{T}}_a})
\end{equation}

Thus, Eq. (\ref{e.38}) can be written as
\begin{equation}\label{}
  {\bm{T}}(\tau ) = {\bm{I}} + {(2{{\bm{T}}_a} + {{\bm{{\bm T}}}_a} \times {{\bm{T}}_a})^n}
\end{equation}

The factorization (\ref{e.39}) should be iterated $n$ times for ${\bm{T}}( \tau  )$. Then, ${\bm{T_a}}$ no longer has a small value after such an iteration circulated $n$ times according to the following computer cycle language
\begin{equation}\label{}
  for{\kern 1pt} {\kern 1pt} (a = 1:n){\kern 1pt} {\kern 1pt} {{\bm{T}}_a} = 2{{\bm{T}}_a} + {{\bm{T}}_a} \times {{\bm{T}}_a}
\end{equation}

At the end of the $n$ cycles, the computer stores ${\bm{T_a}}$. At this point, ${\bm{T_a}}$ can be directly added to the identity matrix ${\bm{I}}$ to obtain the exponential matrix ${\bm{T}}( \tau)$
\begin{equation}\label{}
 {\bm{T}}(\tau ) = {\bm{I}} + {{\bm{T}}_a} = {e^{{\bm{H}}\tau }}
\end{equation}

Therefore, the Taylor approximation of the PIM can be combined with $2^n$ algorithm to calculate the exponential matrix ${\bm{T}}(t)$ to obtain a high precision numerical solution $\bm{\Phi}(x,t)$.
In the same way, we can use the CFD based on PIM (CFD-PIM) to obtain the solution $\bm{\Phi}_x(x,t)$ of the first derivative of the one-dimensional heat conduction equation.
According to one-dimensional Hopf-Cole transformation (\ref{ee.26}), $\bm{\Phi}(x,t)$ and $\bm{\Phi}_x(x,t)$ are submitted into Eq. (\ref{ee.26}) to get the solution $u(x,t)$ of one-dimensional Burgers' equation.

\subsection{Stability analysis}
\subsubsection{Stability of Periodic boundary condition}
To study the stability of our scheme, we only consider the periodic boundary condition for simplicity.

In Eq. (\ref{equ.64}), if ${\lambda _i}(i = 1,2,...,N - 1)$ is the eigenvalue of matrix ${\bm{H}}$, then ${e^{{\lambda _i}\tau }}$ is the eigenvalue of exponential matrix ${e^{{\bm{H}}\tau }}$ with the same corresponding eigenvector $\bm{x} = {({x_1},{x_2}, \ldots ,{x_{N - 1}})}$.
To prensent that CFD-PIM scheme is unconditionally stability, we need to prove the spectral radius of matrix ${e^{{\bm{H}}\tau }}$ is less than 1.
To this end, the following two lemmas are needed.
\\ \textbf{Lemma 1.} If ${\lambda _i}$ is an eigenvalue of matrix ${\bm{H}} = {{\bm{A}}^{ - 1}}{\bm{B}}$ with its corresponding eigenvector $\bm{x}$, then the eigenvalue ${\lambda _i}$ is real number and ${\lambda _i}\leq0$  .
\\ \textbf{Proof.} By the definitions of eigenvalue and eigenvector, we may write  , implying that \cite{Zhou2018}. This gives
\begin{equation}\label{equ.76}
  {{\bm{x}}^T}{\bm{Bx}} = {\lambda _i}{{\bm{x}}^T}{\bm{Ax}}
\end{equation}
\par Here, for periodic boundary condition the matrix ${\bm{A}}$ and ${\bm{B}}$ are as follows

{
\begin{equation}\label{}
{\bf{A}} = \left[ \begin{array}{*{20}{c}}
1&{\frac{2}{{11}}}&{}&{}&{}\\
{\frac{2}{{11}}}&1&{\frac{2}{{11}}}&{}&{}\\
{}& \ddots & \ddots & \ddots &{}\\
{}&{}&{\frac{2}{{11}}}&1&{\frac{2}{{11}}}\\
{}&{}&{}&{\frac{2}{{11}}}&1
\end{array} \right]
\end{equation}
}
\begin{equation}\label{}
{\bm{B}} = \frac{1}{{{h^2}}}\left[ {\begin{array}{*{20}{c}}
{{\textstyle{{ - 51} \over {22}}}}&{{\textstyle{{12} \over {11}}}}&{{\textstyle{3 \over {44}}}}&{}&{}&{}&{}\\
{{\textstyle{{12} \over {11}}}}&{{\textstyle{{ - 51} \over {22}}}}&{{\textstyle{{12} \over {11}}}}&{{\textstyle{3 \over {44}}}}&{}&{}&{}\\
{{\textstyle{3 \over {44}}}}&{{\textstyle{{12} \over {11}}}}&{{\textstyle{{ - 51} \over {22}}}}&{{\textstyle{{12} \over {11}}}}&{{\textstyle{3 \over {44}}}}&{}&{}\\
{}& \ddots & \ddots & \ddots & \ddots & \ddots &{}\\
{}&{}&{{\textstyle{3 \over {44}}}}&{{\textstyle{{12} \over {11}}}}&{{\textstyle{{ - 51} \over {22}}}}&{{\textstyle{{12} \over {11}}}}&{{\textstyle{3 \over {44}}}}\\
{}&{}&{}&{{\textstyle{3 \over {44}}}}&{{\textstyle{{12} \over {11}}}}&{{\textstyle{{ - 51} \over {22}}}}&{{\textstyle{{12} \over {11}}}}\\
{}&{}&{}&{}&{{\textstyle{3 \over {44}}}}&{{\textstyle{{12} \over {11}}}}&{{\textstyle{{ - 51} \over {22}}}}
\end{array}} \right]
\end{equation}

\par Obviously, the matrix ,${\bm{A}}$ and ${\bm{B}}$ are really symmetrical, so the eigenvalue ${\lambda _i}$ is real number. Meanwhile, for arbitrary ${\bm{x}} \ne {\bm{0}}$, the right-hand side of the Eq. (\ref{equ.76}) is
\begin{equation}\label{}
  {{\bm{x}}^T}{\bm{Ax}} = x_1^2 + {\textstyle{4 \over {11}}}{x_1}{x_2} + x_2^2 + {\textstyle{4 \over {11}}}{x_2}{x_3} +  \cdots  + {\textstyle{4 \over {11}}}{x_{N - 2}}{x_{N - 1}} + x_{N - 1}^2
\end{equation}
\par Using the inequality {$2xy < {x^2} + {y^2}$}, we obtain
\begin{equation}\label{equ.80}
  \begin{array}{l}
{{\bm{x}}^T}{\bm{Ax}} > x_1^2 - {\textstyle{2 \over {11}}}(x_1^2 + x_2^2) + x_2^2 - {\textstyle{2 \over {11}}}(x_2^2 + x_3^2) +  \cdots  \\
{\kern 1pt} {\kern 1pt} {\kern 1pt} {\kern 1pt} {\kern 1pt} {\kern 1pt} {\kern 1pt} {\kern 1pt} {\kern 1pt} {\kern 1pt} {\kern 1pt} {\kern 1pt} {\kern 1pt} {\kern 1pt} {\kern 1pt} {\kern 1pt} {\kern 1pt} {\kern 1pt} {\kern 1pt} {\kern 1pt} {\kern 1pt} {\kern 1pt} {\kern 1pt} {\kern 1pt} {\kern 1pt} {\kern 1pt} {\kern 1pt} {\kern 1pt} {\kern 1pt} {\kern 1pt} {\kern 1pt} {\kern 1pt} {\kern 1pt}
- {\textstyle{2 \over {11}}}(x_{N - 2}^2 + x_{N - 1}^2) + x_{N - 1}^2\\
{\kern 1pt} {\kern 1pt} {\kern 1pt} {\kern 1pt} {\kern 1pt} {\kern 1pt} {\kern 1pt} {\kern 1pt} {\kern 1pt} {\kern 1pt} {\kern 1pt} {\kern 1pt} {\kern 1pt} {\kern 1pt} {\kern 1pt} {\kern 1pt} {\kern 1pt} {\kern 1pt} {\kern 1pt} {\kern 1pt} {\kern 1pt} {\kern 1pt} {\kern 1pt} {\kern 1pt} {\kern 1pt} {\kern 1pt} {\kern 1pt} {\kern 1pt} {\kern 1pt} {\kern 1pt}  > {\textstyle{9 \over {11}}}x_1^2 + {\textstyle{7 \over {11}}}\sum\limits_{i = 2}^{N - 2} {x_i^2}  + {\textstyle{9 \over {11}}}x_{N - 1}^2 > 0
\end{array}
\end{equation}
and the left-hand side of the Eq. (\ref{equ.76}) is
\begin{equation}\label{}
\begin{array}{l}
{{\bm{x}}^T}{\bm{Bx}}{\kern 1pt} \\
{\kern 1pt} {\kern 1pt} {\kern 1pt} {\kern 1pt} {\kern 1pt} {\kern 1pt} {\kern 1pt} {\kern 1pt} {\kern 1pt} {\kern 1pt} {\kern 1pt} {\kern 1pt} {\kern 1pt} {\kern 1pt} {\kern 1pt} {\kern 1pt}  = {\textstyle{{ - 51} \over {22}}}x_1^2 + {\textstyle{{12} \over {11}}}{x_1}{x_2} + {\textstyle{3 \over {44}}}{x_1}{x_3} + {\textstyle{{12} \over {11}}}{x_2}{x_1} - {\textstyle{{51} \over {22}}}x_2^2 + {\textstyle{{12} \over {11}}}{x_2}{x_3} + {\textstyle{3 \over {44}}}{x_{\rm{3}}}{x_{\rm{4}}}\\
{\kern 1pt} {\kern 1pt} {\kern 1pt} {\kern 1pt} {\kern 1pt} {\kern 1pt} {\kern 1pt} {\kern 1pt} {\kern 1pt} {\kern 1pt} {\kern 1pt} {\kern 1pt} {\kern 1pt} {\kern 1pt} {\kern 1pt} {\kern 1pt} {\kern 1pt} {\kern 1pt} {\kern 1pt} {\kern 1pt} {\kern 1pt} {\kern 1pt} {\kern 1pt} {\kern 1pt}  + {\textstyle{3 \over {44}}}{x_{\rm{1}}}{x_3}{\rm{ + }}{\textstyle{{12} \over {11}}}{x_{\rm{2}}}{x_3} - {\textstyle{{51} \over {22}}}x_3^2 + {\textstyle{{12} \over {11}}}{x_3}{x_4} + {\textstyle{3 \over {44}}}{x_4}{x_5} +  \cdots \\
{\kern 1pt} {\kern 1pt} {\kern 1pt} {\kern 1pt} {\kern 1pt} {\kern 1pt} {\kern 1pt} {\kern 1pt} {\kern 1pt} {\kern 1pt} {\kern 1pt} {\kern 1pt} {\kern 1pt} {\kern 1pt} {\kern 1pt} {\kern 1pt} {\kern 1pt} {\kern 1pt} {\kern 1pt} {\kern 1pt} {\kern 1pt} {\kern 1pt} {\kern 1pt}  + {\textstyle{3 \over {44}}}{x_{N - 1}}{x_{N - 3}} + {\textstyle{{12} \over {11}}}{x_{N - 2}}{x_{N - 3}}{\kern 1pt}  - {\textstyle{{51} \over {22}}}x_{N - 3}^2 + {\textstyle{{12} \over {11}}}{x_{N - 3}}{x_{N - 4}}\\
{\kern 1pt} {\kern 1pt} {\kern 1pt} {\kern 1pt} {\kern 1pt} {\kern 1pt} {\kern 1pt} {\kern 1pt} {\kern 1pt} {\kern 1pt} {\kern 1pt} {\kern 1pt} {\kern 1pt} {\kern 1pt} {\kern 1pt} {\kern 1pt} {\kern 1pt} {\kern 1pt} {\kern 1pt} {\kern 1pt} {\kern 1pt} {\kern 1pt}  + {\textstyle{3 \over {44}}}{x_{N - 3}}{x_{N - 5}} + {\textstyle{{12} \over {11}}}{x_{N - 1}}{x_{N - 2}} - {\textstyle{{51} \over {22}}}x_{N - 2}^2 + {\textstyle{{12} \over {11}}}{x_{N - 2}}{x_{N - 3}}\\
{\kern 1pt} {\kern 1pt} {\kern 1pt} {\kern 1pt} {\kern 1pt} {\kern 1pt} {\kern 1pt} {\kern 1pt} {\kern 1pt} {\kern 1pt} {\kern 1pt} {\kern 1pt} {\kern 1pt} {\kern 1pt} {\kern 1pt} {\kern 1pt} {\kern 1pt} {\kern 1pt} {\kern 1pt} {\kern 1pt} {\kern 1pt} {\kern 1pt} {\kern 1pt}  + {\textstyle{3 \over {44}}}{x_{N - 2}}{x_{N - 4}}{\kern 1pt}  - {\textstyle{{51} \over {22}}}x_{N - 1}^2 + {\textstyle{{12} \over {11}}}{x_{N - 1}}{x_{N - 2}} + {\textstyle{3 \over {44}}}{x_{N - 1}}{x_{N - 3}}
\end{array}
\end{equation}
\par Similarly, using the inequality {$2xy < {x^2} + {y^2}$}, we obtain
\begin{equation}\label{equ.82}
  {{\bm{x}}^T}{\bm{Bx}} < { - \frac{{12}}{{11}}\left( {x_1^2 + x_{N - 1}^2} \right) - \frac{3}{{44}}\left( {x_1^2 + x_2^2 + x_{N - 2}^2 + x_{N - 1}^2} \right)} < 0
\end{equation}
\par From Eqs. (\ref{equ.80}) and (\ref{equ.82}), we can see that the right-hand side of the Eq. (\ref{equ.76}) is ${{\bm{x}}^T}{\bm{Ax}} > 0$ and the left-hand side of the Eq. (\ref{equ.76}) is ${{\bm{x}}^T}{\bm{Bx}} < 0$. Thus, ${\lambda _i} \le 0$.
\\ \textbf{Lemma 2.} Let ${\bm{W}}$ be an arbitrary square matrix. Then for any operator matrix norm $\left\| \cdot \right\|$, we obtain ${\lambda _i}(\bm{W}) \le \left\| \bm{W} \right\|$, where ${\lambda _i}(\bm{W})$ is the spectral radius of matrix ${\bm{W}}$ \cite{Han2013}.
\\ \textbf{Theorem 1.} The Taylor approximation of CFDS-PIM is unconditionally stability.
\\ \textbf{Proof.} the Taylor series approximation of ${e^{{\lambda _i}\tau }}$ is defined as	
\begin{equation}\label{}
  {e^{{\lambda _i}\tau }}{\rm{ = }}\sum\limits_{j = 0}^\infty  {\frac{{{{({\lambda _i}\tau )}^j}}}{{j!}}}
\end{equation}
\par Using Lemma 1, we obtain ${\lambda _i}\tau  \le 0$, thus ${e^{{\lambda _i}\tau }} \le 1$. For fourth-order Taylor series approximation of the PIM, ${e^{{\lambda _i}\tau }} < 1$.
We use Lemma 2 to get the spectral radius of matrix ${e^{{\bm{H}}\tau }}$ is less than 1 in fourth-order Taylor approximation of the PIM. Thus, the Taylor approximation of CFD-PIM scheme is unconditionally stability.
\textcolor{blue}{\subsubsection{Amplification symbol}
\textbf{Definition 1.} The rational approximation $R_{r,s}(z)$ to the exponential $e^{-z}$ is called A-acceptable when $\left| {{R_{r,s}}( - z)} \right| < 1$ holds for all $-z$ with negative real part. The approximation is called L-acceptable when it is A-acceptable and it also satisfies $\left| {{R_{r,s}}( - z)} \right| \to 0$ as $\Re ( - z) \to  - \infty $.}
\begin{figure}[h!]
  \includegraphics[width=4.5in]{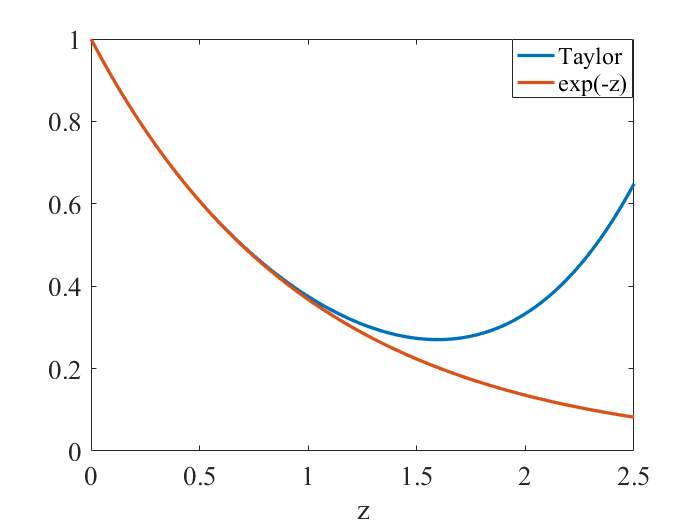}

  \textcolor{blue}{\caption{The behavior of Taylor approximation with the real number field $z \in [0,2.5]$.}\label{fn1}}
  \end{figure}
\begin{figure}[h!]
  \includegraphics[width=2.35in]{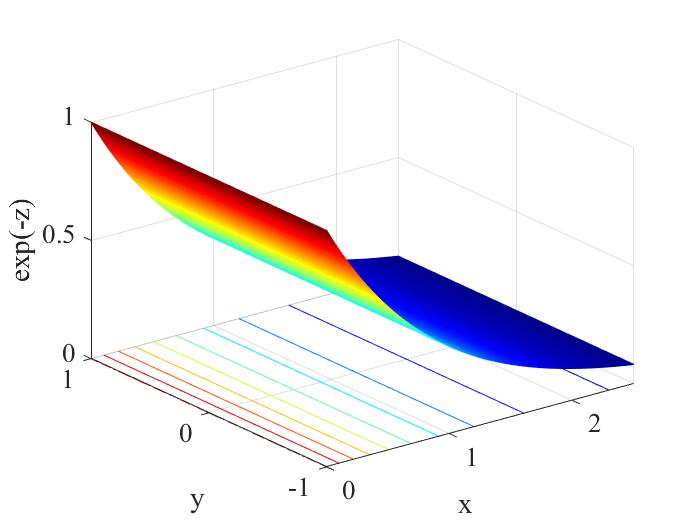}
  \includegraphics[width=2.35in]{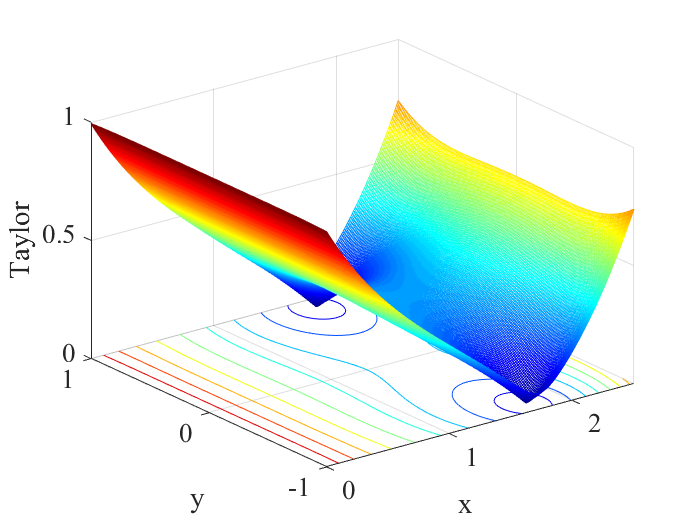}
  \textcolor{blue}{\caption{The behavior of Taylor approximation with the complex plane $z=x+iy$$\in [0,2.5]\times[-1,1]$}\label{fn2}}
  \end{figure}
\par \textcolor{blue}{In Fig. \ref{fn1}, we compare the behavior of $e^{-z}$ and Taylor approximation. It can be observed from the traces that Taylor approximation are A-acceptable.}
\par \textcolor{blue}{Fig. \ref{fn2} illustrates the traces of $e^{-z}$ and Taylor approximation for the different complex planes. Since the results of the functions are complex, we plot their real parts. It can be seen from the plots that Taylor approximation conforms to A-acceptable Definition 1.
}
\textcolor{blue}{\subsubsection{Stability region}
The stability of the CFDS-PIM can be observed from the plots of their stability region \cite{Bhatt2016,Liang2018}.The linear ordinary differential equation (\ref{ee.66}) can be rewrite as
\begin{equation}\label{eeee}
  {\bm{u}}_t=c{\bm{u}}
\end{equation}
We assume that a fixed point $u_0$ satisfying $cu_0= 0$ exists, and $u$ is the perturbation of $u_0$. If $Re(c)< 0$, then we can say the fixed point $u_0$ is stable. We denote $x = c\tau$, with $\tau$ being a single time step, and then apply Taylor approximation to Eq. (\ref{eeee}). The amplification factors $r(x)$ of Taylor approximations can be calculated in the following way:
\begin{equation}
\frac{u_{n+1}}{u_{n}}=r(x,y)=1+x+\frac{1}{2}x^{2}+\frac{1}{6}x^{3}+\frac{1}{24}x^{4}
\end{equation}
\begin{figure}[h!]
  \includegraphics[width=4.5in]{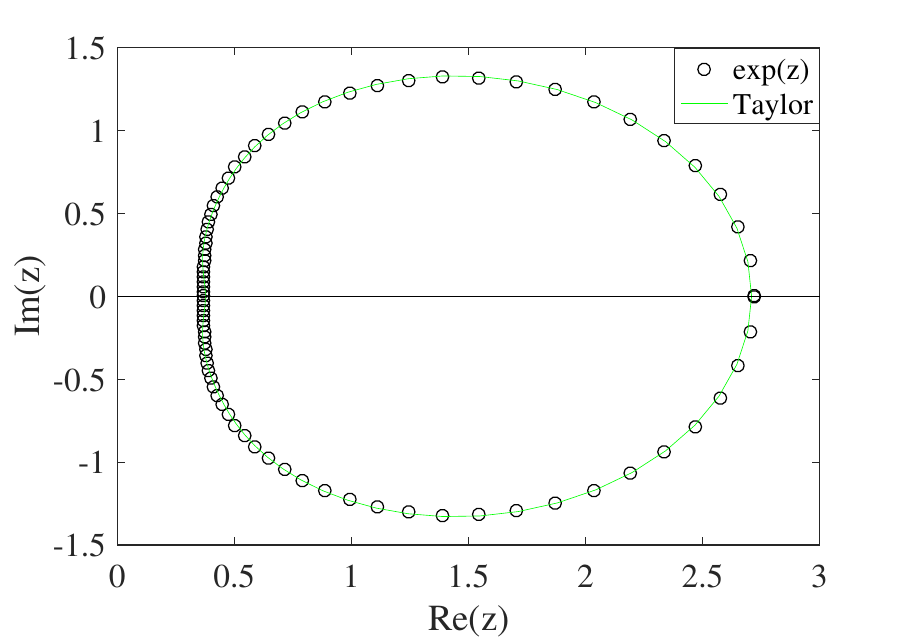}
  \textcolor{blue}{\caption{The stability regions of Taylor approximation}\label{fn4}}
  \end{figure}}
\par \textcolor{blue}{Notice that we assumed $r(x) < 1$ to obtain the stability region. Suppose that $x$ is complex. As can be seen in Fig. \ref{fn4}, the stability regions of proposed scheme is plotted. The axes of Fig. \ref{fn4} are real and imaginary parts of $z$. It can be observed from Fig. \ref{fn4} that the stability region of Taylor approximation are in good agreement with exponential approximation $e^z$.
}

\section{The n-dimensional numerical method}

Strang splitting method (SSM) is a numerical method for solving differential equations that are decompose multi-dimensional problems into a sum of differential operators.
This method is named after Gilbert Strang.
It is used to speed up {the} calculation for problems involving operators on very different time scales, and to solve the multi-dimensional PEDs by reducing them to a sum of one-dimensional problems.

\subsection{Extensions to two-dimensional case}

We consider the two-dimensional heat conduction equation (\ref{e.11}). As a precursor to Strang splitting, Eq. (\ref{e.11}) can be written as

\begin{equation}\label{e.51}
\bm{\Phi}_t = {\bm{{H_x}\bm{\Phi}}} + {\bm{{H_y}\bm{\Phi}}}
\end{equation}
where ${\bm{H_x}}$ and ${\bm{H_y}}$ are difference operator in the $x$ and $y$ direction. The right side of Eq. (\ref{e.51}) is already split into a sum $a+b$ of relatively simple expressions.
Due to one of the properties of difference operator is the distributive law of multiplication, we obtain the following equations
\begin{equation}\label{e.52}
\bm{\Phi}_t = {\bm{({H_x}+{H_y})\bm{\Phi}}}
\end{equation}

 For Eq. (\ref{e.52}), the analytical solution to the associated initial value problem would be

\begin{equation}\label{equ.86}
  {{\bm{\Phi}}^{k + 1}}(t) = {e^{\bm{({H_x}+{H_y})}t}}{{\bm{\Phi}}^k}
\end{equation}

This section focuses on how to calculate the exponential matrix $e^{\bm{({H_x}+{H_y})}t}$, and the calculation of $e^{\bm{({H_x}+{H_y})}t}$ is too complicated. Thus, we convert it into calculating the product of $e^{\bm{{H_x}}t}$ and $e^{\bm{{H_y}}t}$, but $e^{\bm{{H_x}}t}$ and $e^{\bm{{H_y}}t}$ must satisfy the commutativity of the addition theorem
\begin{equation}\label{}
  e^{\bm{({H_x}+{H_y})}t} =e^{\bm{{H_x}}t}e^{\bm{{H_y}}t} \Leftrightarrow {\bm{{H_x}{H_y}}} = {\bm{{H_y}{H_x}}}
\end{equation}

Besides, the exponentials of ${\bm{H_x}}$ and ${\bm{H_y}}$ are related to that of ${\bm{H_x+H_y}}$ by the Trotter product formula
\begin{equation}\label{}
  {e^{\bm{{H_x}+{H_y}}}} = \mathop {\lim }\limits_{m \to \infty } {\left( {{e^{{{\bm{H_x}} \mathord{\left/
 {\vphantom {{{H_x}} m}} \right.
 \kern-\nulldelimiterspace} m}}}{e^{{{\bm{H_y}} \mathord{\left/
 {\vphantom {{{H_y}} m}} \right.
 \kern-\nulldelimiterspace} m}}}} \right)^m}
\end{equation}

Gottleib et al. \cite{Moler2003} suggested that the Trotter result can be used to approximated $e^{\bm{H}}$ by splitting $\bm{H}$ into $\bm{H_x+H_y}$, because $m=2^{20}$ is already very large that was proposed. Thus, we use the following approximation
\begin{equation}\label{}
  {e^{\bm{H}}} = {\left( {{e^{{{\bm{H_x}} \mathord{\left/
 {\vphantom {{\bm{H_x}} m}} \right.
 \kern-\nulldelimiterspace} m}}}{e^{{{\bm{H_y}} \mathord{\left/
 {\vphantom {{\bm{H_y}} m}} \right.
 \kern-\nulldelimiterspace} m}}}} \right)^m}
\end{equation}

This approach to calculate ${e^{\bm{H}}}$ is of potential interest when the exponentials of ${\bm{H_x}}$ and ${\bm{H_y}}$ can be accurately and efficiently computed. If ${\bm{H_x}}$ and ${\bm{H_y}}$ commute, we rewrite Eq. (\ref{equ.86}) as follows
\begin{equation}\label{}
  {{\bm{\Phi}}^{k + 1}}(t) = {e^{{\bm{H}}t}}{{\bm{\Phi}}^k} = {e^{{{\bm{H}}_x}t}}{e^{{{\bm{H}}_y}t}}{{\bm{\Phi}}^k}
\end{equation}

Thus, the two-dimensional heat conduction equation becomes two one-dimensional problems. For each one-dimensional problem, it can be solved by the CFD-PIM scheme which introduced in Sec. \ref{3.1}.
\subsection{Extensions to three-dimensional case}
For the three-dimensional heat conduction equation, we can also use CFD-PIM based on the SSM(CFD-PIM-SSM) to decompose it into the sum of differential operators of three one-dimensional problems. The CFD-PIM scheme can be extended to three-dimensional case (\ref{e.12}).

As a precursor to Strang splitting, we rewrite Eq. (\ref{e.12}) as follows

\begin{equation}\label{e.58}
\bm{\Phi}_t = {\bm{{H}_x\bm{\Phi}}} + {\bm{{H}_y\bm{\Phi}}}+{\bm{{H}_z\bm{\Phi}}}
\end{equation}
where ${\bm{H_x}}$, ${\bm{H_y}}$ and ${\bm{H_z}}$ are difference operators in the $x$-direction, $y$-direction, and $z$-direction, respectively. The right side of Eq. (\ref{e.58}) is already split, which become a sum  $a+b+c$. We obtain the following equations
\begin{equation}\label{}
   \bm{\Phi}_t = ({\bm{{H}_x}} + {\bm{{H}_y}}+{\bm{{H}_z\bm){\Phi}}}
\end{equation}

 For Eq. (\ref{e.52}), the analytical solution to the associated initial value problem would be

\begin{equation}\label{equ.93}
  {{\bm{\Phi}}^{k + 1}}(t) = {e^{\bm{({H_x}+{H_y}+{H_z})}t}}{{\bm{\Phi}}^k}
\end{equation}

If ${\bm{H_x}}$, ${\bm{H_y}}$ and ${\bm{H_z}}$ commute for Eq. (\ref{equ.93}), the analytical solution to the associated initial value problem would be

\begin{equation}\label{equ.94}
  {{\bm{\Phi}}^{k + 1}}(t) = {e^{{\bm{H}}t}}{{\bm{\Phi}}^k} = {e^{\bm{H_x}t}}{e^{\bm{H_y}t}}{e^{\bm{H_z}t}}{{\bm{\Phi}}^k}
\end{equation}
\par Because we apply SSM to three-dimensional case, we obtain a sum of difference operator of three one-dimensional parabolic problem, the scheme has the same accuracy as one-dimensional cases. For each one-dimensional problem, it can be solved by the CFD-PIM scheme which introduced in Sec. \ref{3.1}.

\section{Numerical Result}
{In this section, we give the six numerical examples to validate the adaptability of the proposed schemes and compare their accuracy with those which are already available in the literature for solving n-dimensional Burgers' system. The accuracy of the schemes is measured in terms of $L_2$ errors, $L_\infty$ errors, computing time and the rate of convergence of the scheme. In our Tables, CPU(s) is computing time. The rate of convergence(ROC) of proposed schemes is defined as
\begin{equation}\label{}
{\log _2}\frac{{{L_\infty }(2h)}}{{{L_\infty }(h)}}
\end{equation}
where ${L_\infty }(2h)$ and ${L_\infty }(h)$ are discrete maximum absolute errors at $2h$ and $h$. All the numerical experiments are conducted on MATLAB R2016a platforms based on an Intel Core i5-6300HQ 2.30 GHz processor.}
\begin{figure}[h!]
\includegraphics[width=5in]{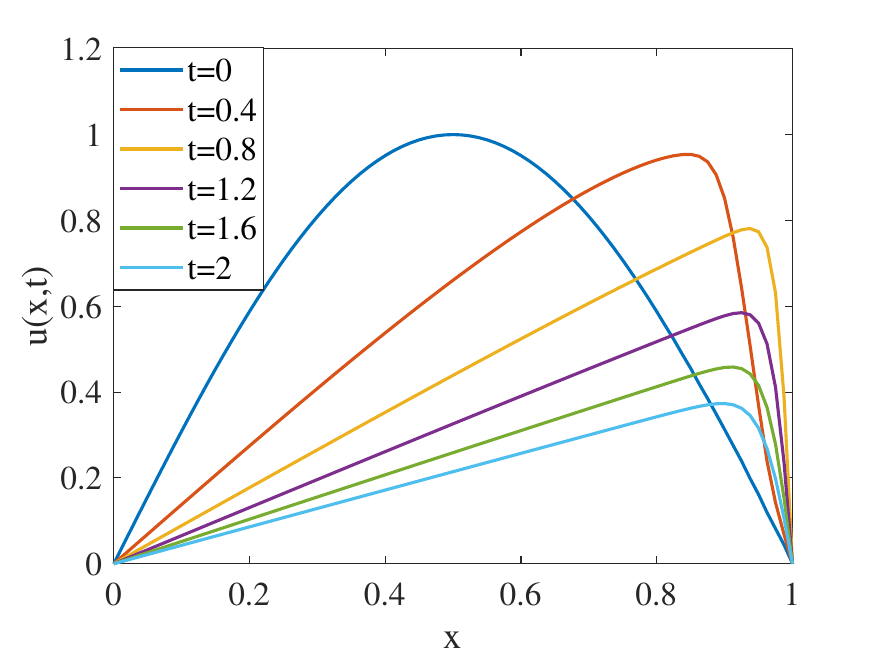}
\caption{ Physical behavior of numerical solutions at different time $t=0,0.4,0.8,1.2,1.6,2.0$ for Example 1}\label{figure.6}
\end{figure}
\\ {\bf{Example 1.}} To verify the effectiveness of the modified algorithm, we test the one-dimensional Burgers' equation (\ref{equ.5}) proposed in Sec. \ref{2.2.1},  over a domain $[0,1]$, with the initial and boundary conditions (\ref{e.60}), according to {the} one-dimensional Hopf-Cole transformation (\ref{ee.26}), the analytical (Fourier series) solution of Eq. (\ref{equ.5}) is
\begin{equation}\label{}
\begin{array}{l}
\\
u(x,t) = 2\pi \omega \frac{{\sum\limits_{\alpha  = 0}^\infty  \alpha  {C_\alpha }\exp ( - {\alpha ^2}{\pi ^2}\omega t)\sin (\alpha \pi x)}}{{\sum\limits_{\alpha  = 0}^\infty  {{C_\alpha }\exp ( - {\alpha ^2}{\pi ^2}\omega t)\cos (\alpha \pi x)} }}
\end{array}
\end{equation}
\par The numerical solutions  are reported in Tables \ref{tt5},\ref{tt6},\ref{tta3} and Fig. \ref{figure.6} for the different values of $t$ {with $Re=100$ and $\tau {\rm{ = }}1 \times {10^{ - 5}}$}.
The results of the numerical solution are compared with those of Refs. \cite{Kutluay1999,Jiwari2015,Jiwari2012, Mittal2012a}, and the numerical results of the proposed scheme are better than their results.
It is evident that the proposed scheme has high accuracy and efficiency than other numerical schemes.
Fig. \ref{figure.6} exhibit that as the increase of time the numerical solution of partial regions becomes steeper and steeper, and the decreasing rate of the approximate solution increases.
This physical phenomenon validates the fact that the numerical solution is capable of describing the shock wave.
For the numerical solutions at different time demonstrated in Fig. \ref{figure.6} are fantastically analogical as depicted in the figures given in Refs. \cite{Jiwari2015,Jiwari2012, Mittal2012a,Li2016a,G.W.Wei1998}. {According to the values of ROC in Table \ref{tta3} ( ROC$>4$ ), the proposed scheme can be verified as a high order scheme}
\begin{table}[h!]
\caption{Comparison with the following numerical and analytical solutions of different values of $x$ and $t$ with $Re=10$ for Example 1.}\label{tt5}
\scriptsize
\centering
\begin{tabular}{llllllll}
  \hline
  $x$ & $t$ & Ref.\cite{Kutluay1999} & Ref.\cite{Kadalbajoo2006} & Ref.\cite{Jiwari2012} & Ref.\cite{Jiwari2015}& \textcolor{green}{Proposed} scheme  & \textcolor{green}{Analytical} solution \\
   \hline
  0.25 & 0.4 & 0.30891 & 0.30881  & 0.30887 & 0.30889 & 0.308894228585555 & 0.308894228585318 \\
      & 0.6 & 0.24075 & 0.24069 & 0.24070 & 0.24075 & 0.240739023291803 & 0.240739023291448 \\
     & 0.8 & 0.19568 & -- & 0.19566 & 0.19569 &0.195675570103972 & 0.195675570103439 \\
    & 1.0 & 0.16257  & 0.16254  & 0.16255 & 0.16258 & 0.162564857111346 & 0.162564857110671\\
      & 3.0 & 0.02720  &0.02720  & 0.02721 & 0.02720 & 0.027202314473410 &0.027202314472951 \\
  0.50 & 0.4 & 0.56964 & 0.56955 & 0.56956  & 0.56956 & 0.569632450695361 & 0.569632450693995 \\
     & 0.6 &0.44721 & 0.44714  & 0.44715  & 0.44724  & 0.447205521200320 & 0.447205521198742 \\
     & 0.8 & 0.35924 & -- &0.35920 & 0.35927  & 0.359236058517410 & 0.359236058515669 \\
       & 1.0 & 0.29192 & 0.29188  &0.29188 & 0.29195  &0.291915957127591 & 0.291915957125836 \\
       & 3.0 & 0.04021 & 0.04021 & 0.04022 & 0.04021 &0.040204924438755 & 0.040204924438046 \\
  0.75 & 0.4 & 0.62542 & 0.62540 & 0.62540 & 0.62537 & 0.625437893711249 & 0.625437893706948 \\
      & 0.6 & 0.48721 & 0.48715 & 0.48716 & 0.48718 & 0.487214974885767 & 0.487214974882163 \\
     & 0.8 & 0.37392 & -- & 0.37389 & 0.37391 &0.373921753212449 & 0.373921753209455 \\
     & 1.0 & 0.28748 & 0.28744 & 0.28743 & 0.28747 & 0.287474405919467 & 0.287474405916976 \\
      & 3.0 & 0.02977 & 0.02978 & 0.02978 & 0.02977 & 0.029772126859293 & 0.029772126858766 \\
  \hline
\end{tabular}
\end{table}
\begin{table}[h!]
\caption{Comparison with the following numerical and analytical solutions of different values of $x$ and $t$ with $Re=100$ for Example 1.}\label{tt6}
\scriptsize
\centering
\begin{tabular}{llllllll}
  \hline
  $x$ & $t$ & Ref.\cite{Kutluay1999} & Ref.\cite{Kadalbajoo2006} & Ref.\cite{Jiwari2012} & Ref.\cite{Jiwari2015}& \textcolor{green}{Proposed} scheme  & \textcolor{green}{Analytical} solution \\
   \hline
  0.25 & 0.4 & 0.34819 & 0.34229  & 0.34184 & 0.34191 & 0.341914932413026 & 0.341914932411983 \\
      & 0.6 & 0.27536 & 0.26902 & 0.26891 & 0.26896 & 0.268964845317425 & 0.268964845316620 \\
     & 0.8 & 0.22752 & -- & 0.22143 & 0.22148 & 0.221481914524793 & 0.221481914524373 \\
    & 1.0 & 0.19375  & 0.18817  & 0.18815 & 0.18820 & 0.188193961397110 & 0.188193961396738 \\
      & 3.0 & 0.07754 &0.07511  & 0.07510 & 0.07511 & 0.075114083887341 &0.075114083887190 \\
  0.50 & 0.4 & 0.66543 & 0.66797 & 0.66060 & 0.66069 & 0.660710972121299 & 0.660710970851541 \\
     & 0.6 & 0.53525 & 0.53211 & 0.52932 & 0.52942 & 0.529418263880240 & 0.529418263729147 \\
     & 0.8 & 0.44526 & -- & 0.43905 & 0.43914 & 0.439138250704212 & 0.374420037644682 \\
       & 1.0 & 0.38047 & 0.37500 & 0.37436 & 0.37443 & 0.374420037662816 & 0.374420037644682 \\
       & 3.0 & 0.15362 & 0.15018 & 0.15017 & 0.15019 & 0.150179005234648 & 0.150179005235832 \\
  0.75 & 0.4 & 0.91201 & 0.93680 & 0.91026 & 0.91023 & 0.910227039712648 & 0.910268136079484 \\
      & 0.6 & 0.77132 & 0.77724 & 0.76719 & 0.76723 & 0.767241968056591 & 0.767243282478514 \\
     & 0.8 & 0.65254 & -- & 0.64745 & 0.64740 & 0.647395126988680 & 0.647395234822761 \\
     & 1.0 & 0.56157 & 0.56157 & 0.55608 & 0.55606 & 0.556050682504419 & 0.556050704468246 \\
      & 3.0 & 0.22874 & 0.22485 & 0.22504 & 0.22486 & 0.224811248098947 & 0.224811248193590 \\
  \hline
\end{tabular}
\end{table}
\begin{table}[h!]
\scriptsize
\centering
{
\caption{\newline Numerical results of $Re=5$ and $Re=10$ with $t=1$ for Example 1.}\label{tta3}
      \begin{tabular}{llllll}
  \toprule
   &$N$ & $11$ & $21$ & $41$ & $81$ \\
              \midrule
      $Re=5$ & ${L_\infty }$ & 1.4625E-06 & 4.2340E-09 & 2.4716E-11 & 3.0104E-13 \\
         & CPU(s)    & 0.604 &0.674  & 0.780 & 1.508 \\
         & ROC    & --     & 8.4322 & 7.4204 & 6.3340 \\
      $Re=10$ & ${L_\infty }$ & 2.8131E-06 & 7.3746E-09 & 3.7254E-11 & 2.1694E-12 \\
         & CPU(s) &1.150  &1.207  & 1.479 & 2.786 \\
         & ROC & -- & 8.5754 & 7.6374 & 4.1719\\
 \bottomrule
\end{tabular}
}
\end{table}
\\
\\{\bf{Example 2.}} To validate the order of convergence of the proposed scheme, a numerical experiment of coupled Burgers' equation (\ref{equ.1}) was carried out in Example 2 with the region $x \in \Omega  = [-\pi,\pi]$ with initial conditions
\begin{equation}
\begin{array}{l}
u(x,0) = \sin x{\kern 1pt} {\kern 1pt} ,{\kern 1pt} {\kern 1pt} {\kern 1pt} x \in \Omega  = [-\pi,\pi]\\
v(x,0) = \sin x{\kern 1pt} {\kern 1pt} ,{\kern 1pt} {\kern 1pt} {\kern 1pt} x \in \Omega  = [-\pi,\pi]
\end{array}
\end{equation}
 extracted from the following exact solution given by Refs. \cite{Bhatt2016,Lai2014} for $\omega_1=\omega_2=1.0$, $\kappa_1=\kappa_2=-2.0$, $\delta_1=\delta_2=1.0$:
\begin{equation}\label{}
  u\left( {x,t} \right) = v\left( {x,t} \right) = \exp ( - t)\sin x,{\kern 1pt} {\kern 1pt} {\kern 1pt} x \in \Omega  = [-\pi,\pi],t>0
\end{equation}
\par
\begin{figure}[h!]
  \includegraphics[width=2.35in]{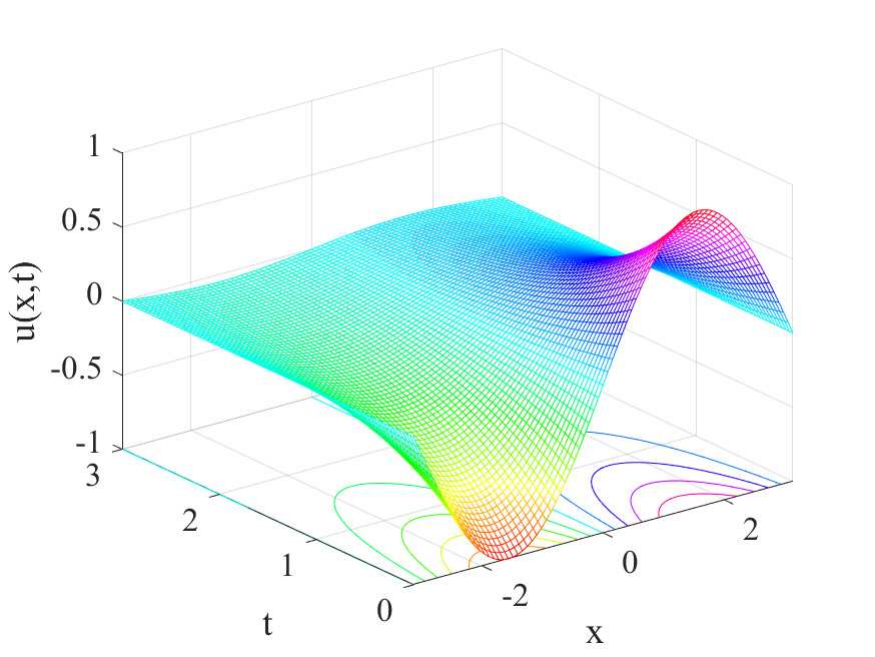}
  \includegraphics[width=2.35in]{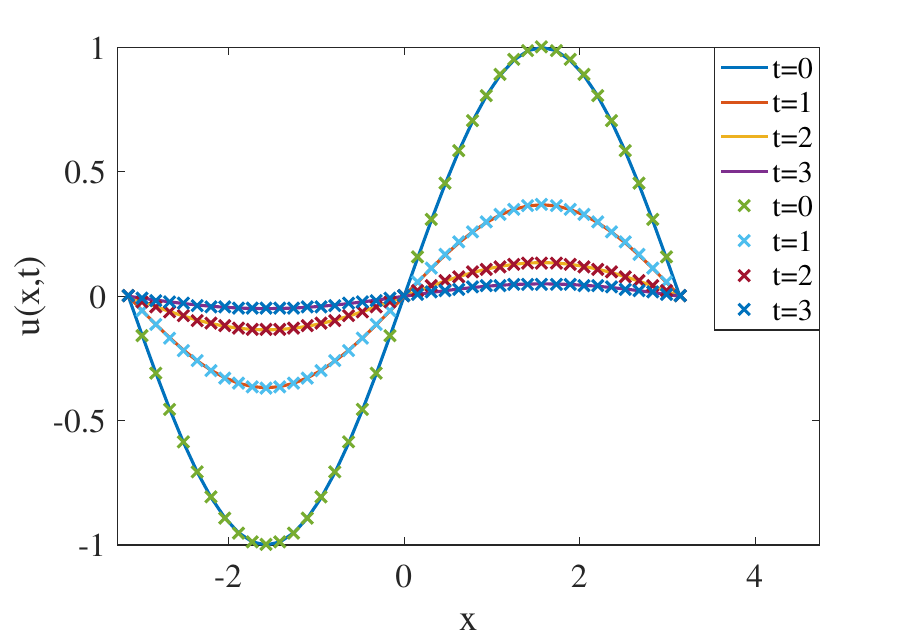}
  \caption{Evolution profile of numerical scheme at time $T =[0,3]$  with $N=41$ (Left), A comparison between numerical and analytical solutions at different time $t$ (Right) with $N=41$ for Example 2. (The solid line represents the analytical solution and the other curve is the numerical solution.)}\label{ff7}
  \end{figure}

 The boundary conditions are extracted from the analytical solution(In this example, $\tau=4\times10^{-4}$).
The spatiotemporal evolution of the numerical solution is shown in the left part of Fig. \ref{ff7}.
To intuitively observe the physical phenomena of the example,
The left figure of Fig. \ref{ff7} is depicted to visually compare the analytical solutions with the numerical solutions of $u ( x , t )$ at different time $t$.
It is observed that the numerical results show great agreement with the analytical solutions.
The numerical results reflect the motion characteristics of wave propagation:
the amplitude of the wave decreases with time while the wavelength remains unchanged.
In addition, the errors of $L_{\infty}$, rates of convergence and CPU running time are listed in Table \ref{tt12}.
As can be observed in Table \ref{tt12}, the accuracy and efficiency of our numerical scheme are much higher than that of Ref. \cite{Bhatt2016} under the same spatial step size.
It is observed that the error of $L_{\infty}$ becomes smaller as the mesh size is refined. The computer operating environment of the algorithm is worse than that of Ref. \cite{Bhatt2016},
and the computation time and convergence order of our numerical scheme are much better than that of Ref. \cite{Bhatt2016},
which shows that the proposed scheme has excellent adaptability.
The proposed scheme presents more accurate and high-efficient solutions in spatial direction than the scheme in Refs .\cite{Bhatt2016,Bak2019}.
\begin{table*}[]
\scriptsize
\centering
\caption{$L_{\infty}$, rates of convergence, and CPU time of the following numerical schemes in Matlab for Example 2.}\label{tt12}
\begin{tabular}{lllllll}
\toprule
\multirow{1}{*}{}
 & \multicolumn{3}{c}{Ref.\cite{Bhatt2016}}
 & \multicolumn{3}{c}{\textcolor{green}{Proposed} scheme}\\
\cmidrule(r){2-4} \cmidrule(r){5-7}
$h$&  $L_\infty $      &  {ROC}      &  CPU(s)
   &  $L_\infty $      &  {ROC}      &  CPU(s)
 \\
\midrule
 $\pi/8$       & {3.660E-05} & -     &0.1492 & 1.149E-05   &  - & 0.071              \\
$\pi/16$         &2.277E-06  &4.0066 & 0.3109     & 2.968E-08   &8.8967       & 0.079                 \\

$\pi/32$         &1.422E-07     &4.0017 & 0.5419      & 1.329E-10
     & 7.8030 & 0.081                 \\
$\pi/64$             &8.882E-09      &4.0004 & 1.0726          & 1.480E-12                                &6.4886   & 0.118\\
\multirow{1}{*}{platform}
 & \multicolumn{3}{c}{Intel Core i7-4510U 2.60 GHz workstation}
 & \multicolumn{3}{c}{Intel Core i5-6300HQ 2.30 GHz processor}\\
\bottomrule
\end{tabular}
\end{table*}
\\
\\{\bf{Example 3.}} To test the applicability of the proposed scheme for two-dimensional problems, we consider the two-dimensional Burgers' equations (\ref{equ.6}) over a square domain $[0,1]\times[0,1]$, with the initial conditions
\begin{equation}\label{}
\begin{array}{l}
u(x,y,0) =  - 2\omega \frac{{ - 2\pi \cos 2\pi x\sin \pi y}}{{2 + \sin 2\pi x\sin \pi y}}\\
v(x,y,0) =  - 2\omega \frac{{ - 2\pi \sin \pi x\cos \pi y}}{{2 + \sin 2\pi x\sin \pi y}}
\end{array}
\end{equation}
for which the analytical solutions \cite{Liao2010A} are
\begin{equation}\label{}
u(x,y,t) = 4\pi \omega \frac{{\exp ( - 5{\pi ^2}\omega t)\cos 2\pi x\sin \pi y}}{{2 + \exp ( - 5{\pi ^2}\omega t)\sin 2\pi x\sin \pi y}}
\end{equation}
\begin{equation}\label{}
v(x,y,t) = 4\pi \omega \frac{{\exp ( - 5{\pi ^2}\omega t)\sin2\pi x\cos \pi y}}{{2 + \exp ( - 5{\pi ^2}\omega t)\sin 2\pi x\sin \pi y}}
\end{equation}
\par
\begin{figure}[h!]
  \includegraphics[width=2.35in]{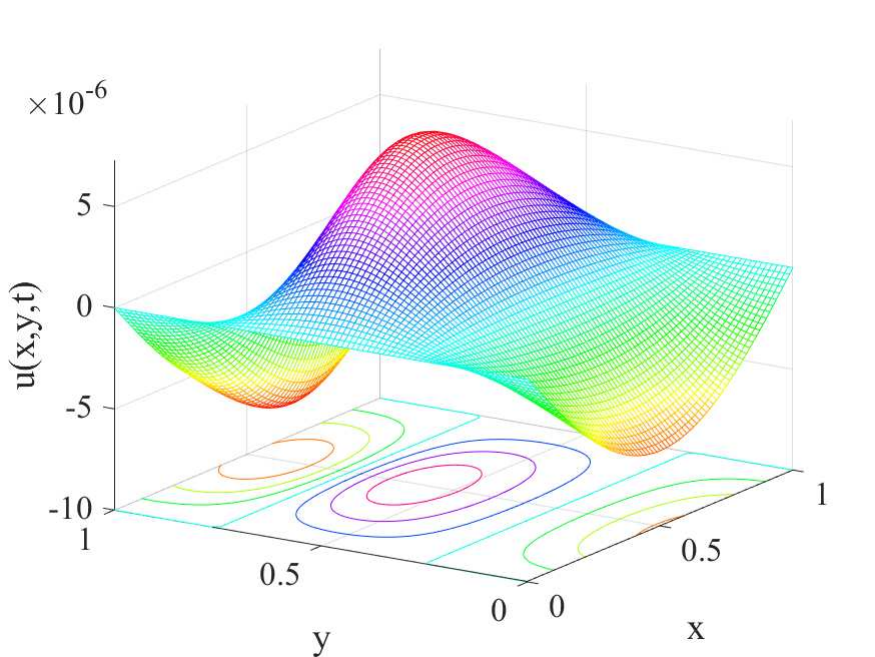}
  \includegraphics[width=2.35in]{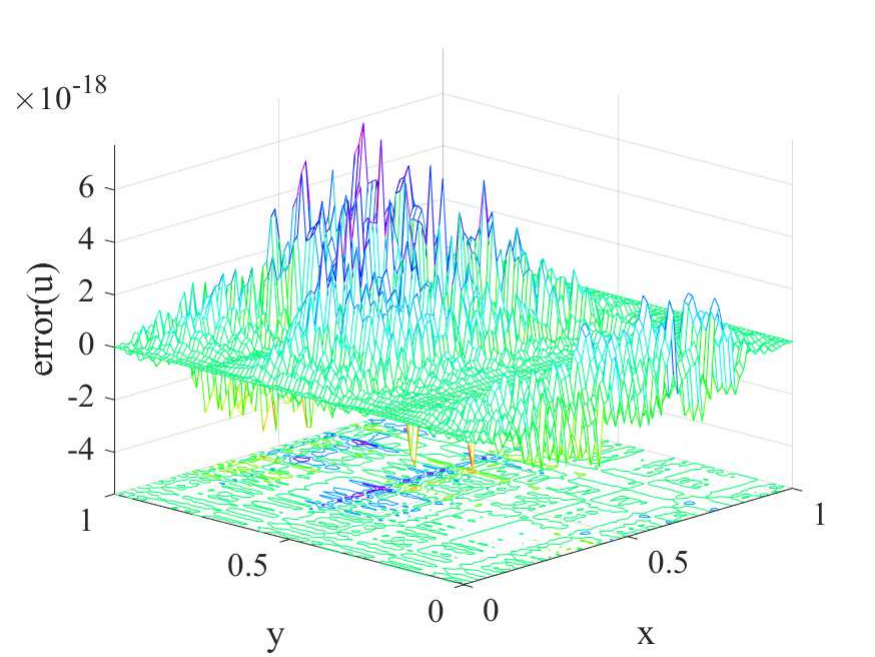}
  \caption{Physical behavior of the numerical solutions $u(x,y,t))$ (Left), and the errors (Right) between the analytical and numerical solutions with $Re=10^6$ and $N \times N= 81\times 81$ at $t=1$ for Example 3. }\label{figure.8}
  \end{figure}

\begin{figure}[h!]
  \includegraphics[width=2.35in]{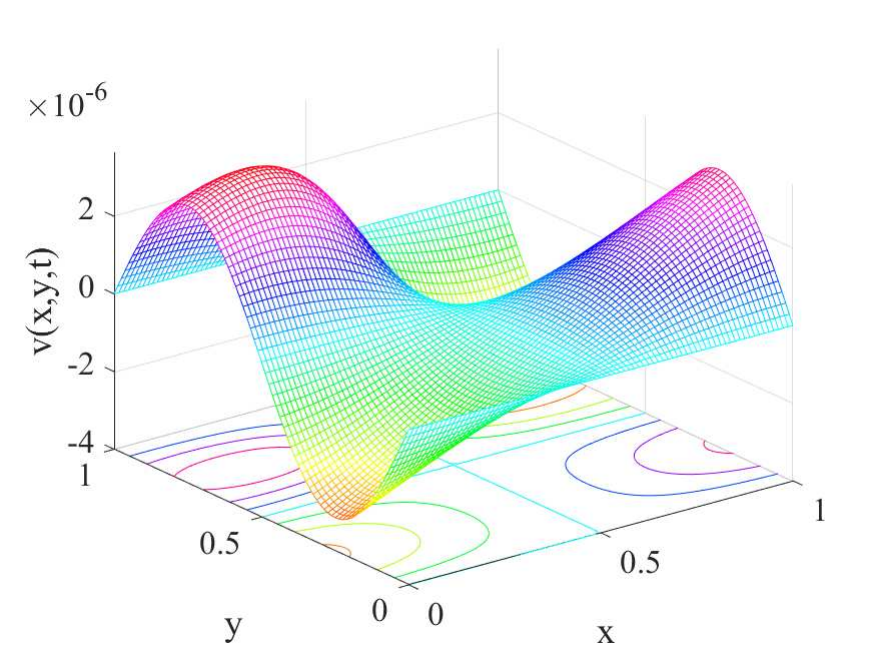}
  \includegraphics[width=2.35in]{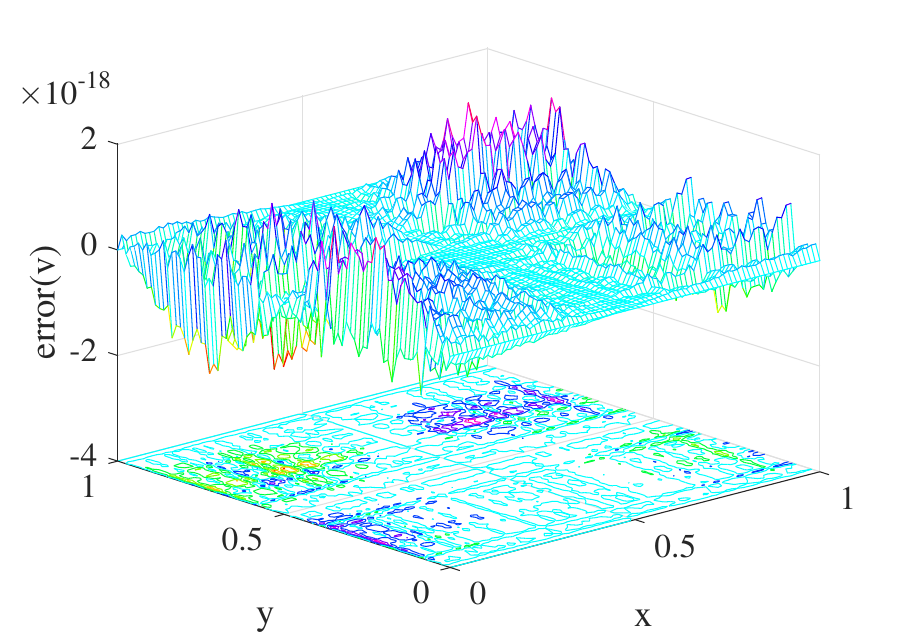}
  \caption{Physical behavior of the numerical solutions $v(x,y,t)$ (Left),and the errors (Right) between the analytical and the numerical solutions with $Re=10^6$ and $N \times N= 81\times 81$ at $t=1$ for Example 3. }\label{figure.9}
  \end{figure}

 The boundary conditions are extracted from the analytical solution.
This example uses the CFD scheme with sixth-order accuracy in space and the PIM with fourth-order accuracy in time.
The numerical simulation results are shown in Table \ref{newta5} and Figs. \ref{figure.8},\ref{figure.9}.
To prove that the method is sixth-order accuracy in space, the time step $\tau$ is fixed to $5\times10^{-4}$, thus the time truncation error can be ignored.
To simplify the demonstration, we use the same grid size in the $x$ and $y$ directions.
 As can be seen from Table \ref{newta5}, when $h$ is reduced by a factor of $2$, the maximal errors for both $u(x,y,t)$ and $v(x,y,t)$ are reduced by a factor of $6$,
which indicates that the method is sixth-order accurate in space.
When $h$ is reduced by a factor of $2$, the $L_{\infty}$ of $u(x,y,t)$ and $v(x,y,t)$ are both reduced by a factor of $2^6$,
which indicates that the method is sixth-order accurate in space.
It is observed that the numerical solutions present excellent agreement with the analytical solutions. {In large Reynolds numbers, the ROC of the proposed scheme has little influence, which means the scheme still is a high order.}
\begin{table}[h!]
\scriptsize
\centering
{
\caption{\newline Comparsion of $L_\infty$ of the different Re at $t=1$ for Example 3.}\label{newta5}
      \begin{tabular}{llllll}
 \toprule
     &$N$ & $11\times11$ & $21\times21$ & $41\times41$ & $81\times81$ \\
             \midrule
     $Re=100$ & $L_\infty(u) $ & 5.2252E-05 & 5.1284E-07 & 2.4176E-09 & 1.4047E-11 \\
         & ROC & -- & 6.6708 & 7.7288 & 7.4272 \\
         & $L_\infty(v) $ & 1.2944E-04 & 1.1137E-06 & 5.0211E-09 & 4.2705E-11 \\
         & ROC & -- & 6.8607 & 7.7932 & 6.8761 \\
    $Re=1000$ & $L_\infty(u) $ & 1.8558E-06 & 4.3507E-08 & 4.0918E-10 & 2.0677E-12 \\
         & ROC & -- & 5.4146 & 6.7323 & 7.6286 \\
         & $L_\infty(v) $ & 2.6879E-06 & 5.8267E-08 & 5.4288E-10 & 2.7184E-12 \\
         & ROC & -- & 5.5277 & 6.7459 & 7.6423 \\
 \bottomrule
\end{tabular}
}
\end{table}
\\
\\{\bf{Example 4.}} In order to verify the effectiveness of the improvement proposed above, we consider the system of {the} two-dimensional Burgers' equations (\ref{equ.6}) proposed in Sec. \ref{2.2.2},  over a square domain $[0,1]\times[0,1]$, with the initial and boundary conditions (\ref{e.70}) and the analytical solution
\begin{equation}\label{}
\begin{array}{l}
u(x,y,t) = 2\pi \omega \frac{{\sum\limits_{\alpha ,\beta  = 0}^\infty  \alpha  {C_{\alpha \beta} }\exp [ - ({\alpha ^2} + {\beta ^2}){\pi ^2}\omega t]\sin (\alpha \pi x)\cos (\beta \pi y)}}{{\sum\limits_{\alpha ,\beta  = 0}^\infty  {{C_{\alpha \beta} }\exp [ - ({\alpha ^2} + {\beta ^2}){\pi ^2}\omega t]\cos (\alpha \pi x)\cos (\beta \pi y)} }}\\
v(x,y,t) = 2\pi \omega \frac{{\sum\limits_{\alpha ,\beta  = 0}^\infty  \beta  {C_{\alpha \beta} }\exp [ - ({\alpha ^2} + {\beta ^2}){\pi ^2}\omega t]cos(\alpha \pi x)\sin (\beta \pi y)}}{{\sum\limits_{\alpha ,\beta  = 0}^\infty  {{C_{\alpha \beta} }\exp [ - ({\alpha ^2} + {\beta ^2}){\pi ^2}\omega t]\cos (\alpha \pi x)\cos (\beta \pi y)} }}
\end{array}
\end{equation}
\par
\begin{figure}[h!]
  \includegraphics[width=2.35in]{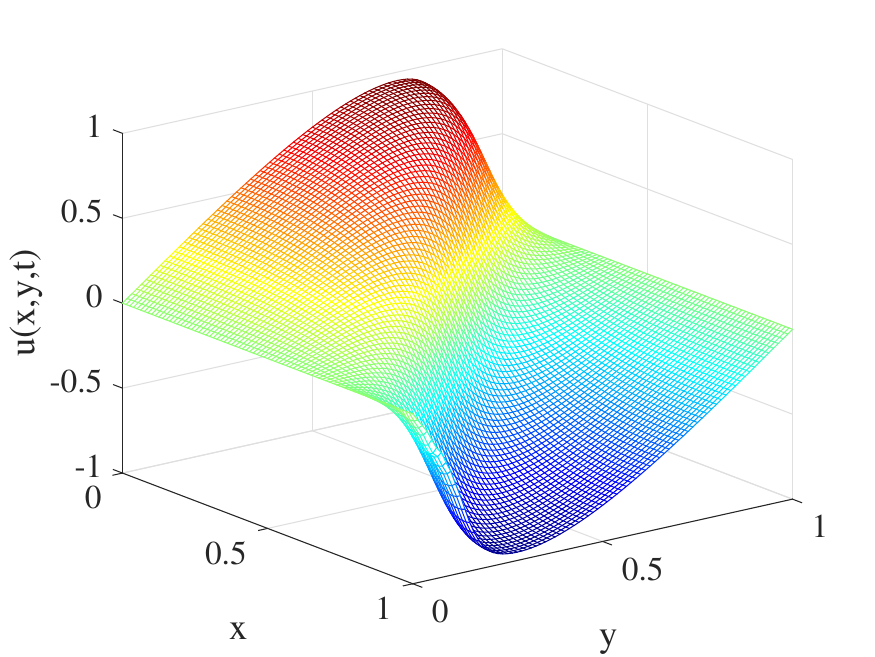}
  \includegraphics[width=2.35in]{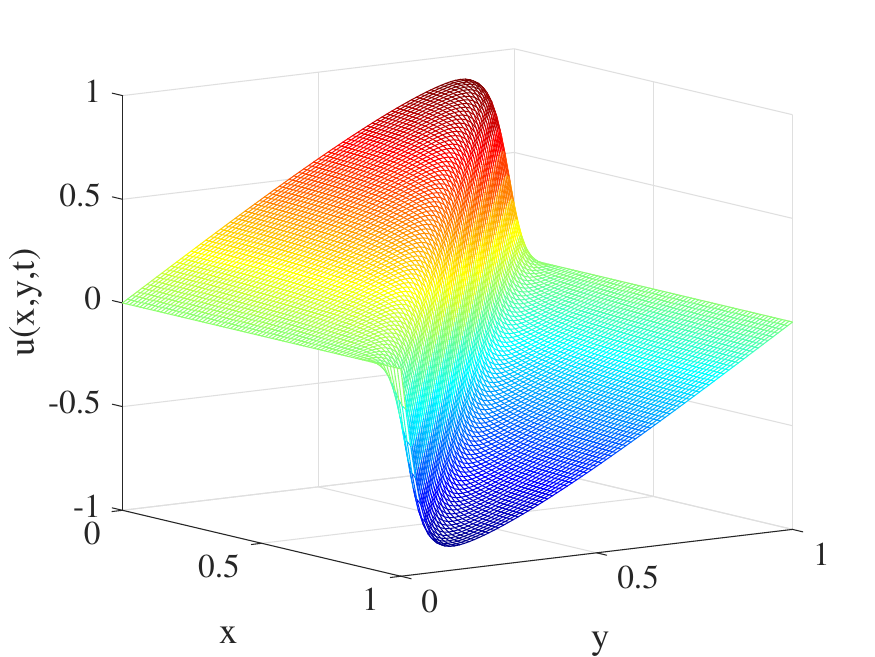}
  \caption{ The numerical solution of the two-dimensional Burgers' equation for $Re = 100$ at $t = 0.25$ (Left) and $t = 0.5$ (Right) with $N \times N= 81\times 81$ for Example 4. }\label{figure.10}
  \end{figure}

\begin{figure}[h!]
  \includegraphics[width=2.35in]{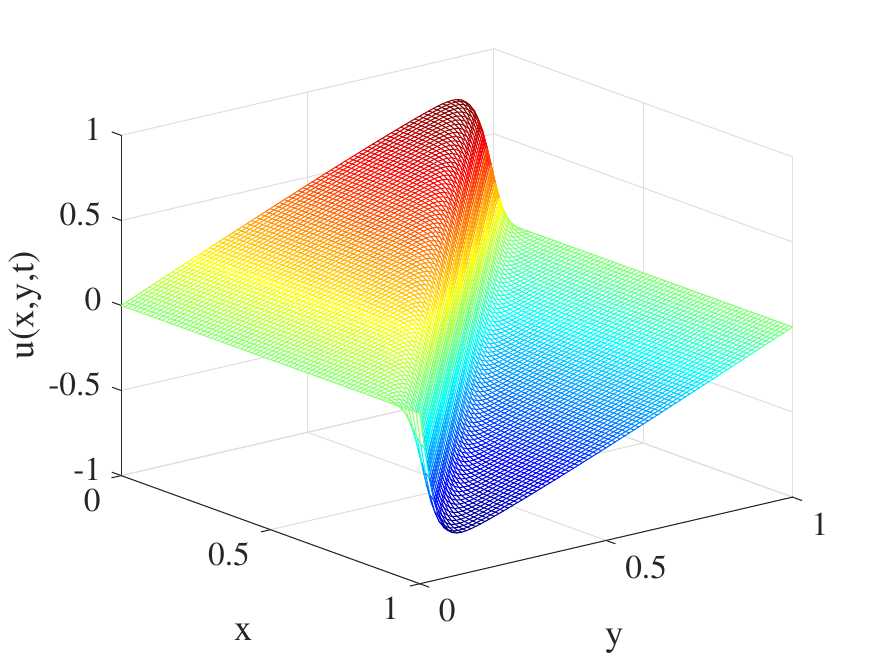}
  \includegraphics[width=2.35in]{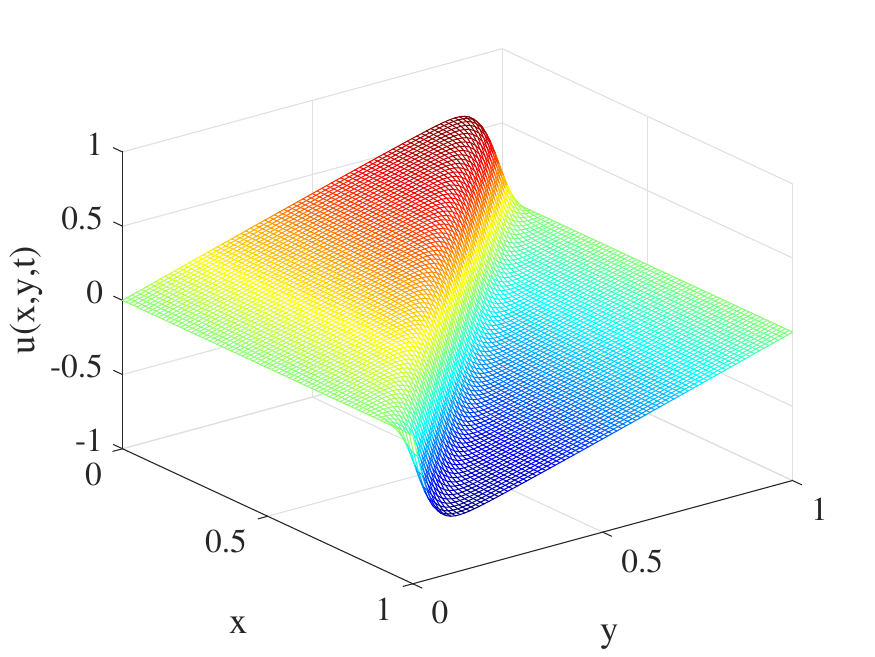}
  \caption{The numerical solution of the two-dimensional Burgers' equation for $Re = 100$ at \textcolor{blue}{$t = 0.75$} (Left) and $t = 1.0$ (Right) with $N \times N= 81\times 81$ for Example 4. }\label{figure.11}
  \end{figure}

\begin{figure}[h!]
  \includegraphics[width=2.35in]{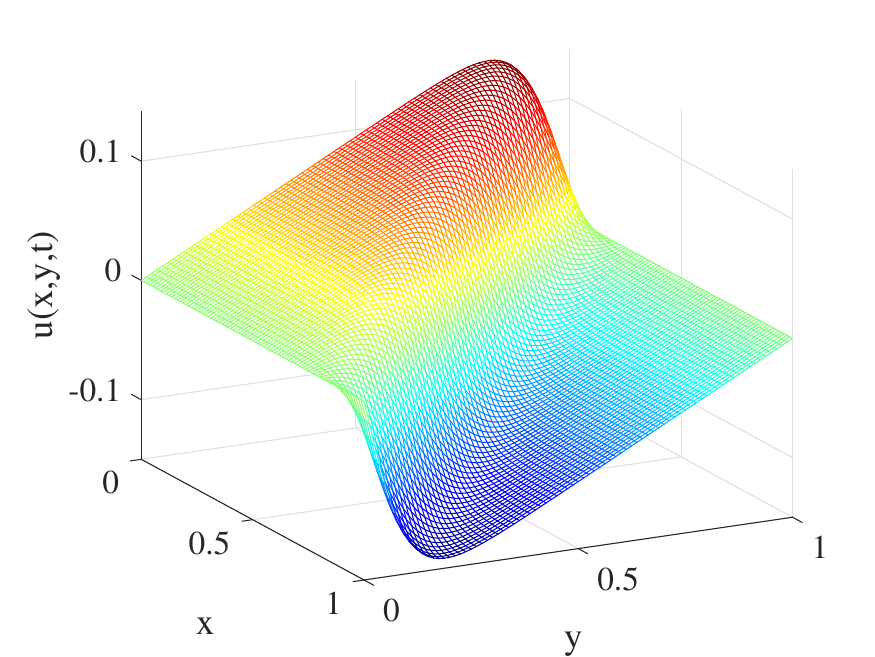}
  \includegraphics[width=2.35in]{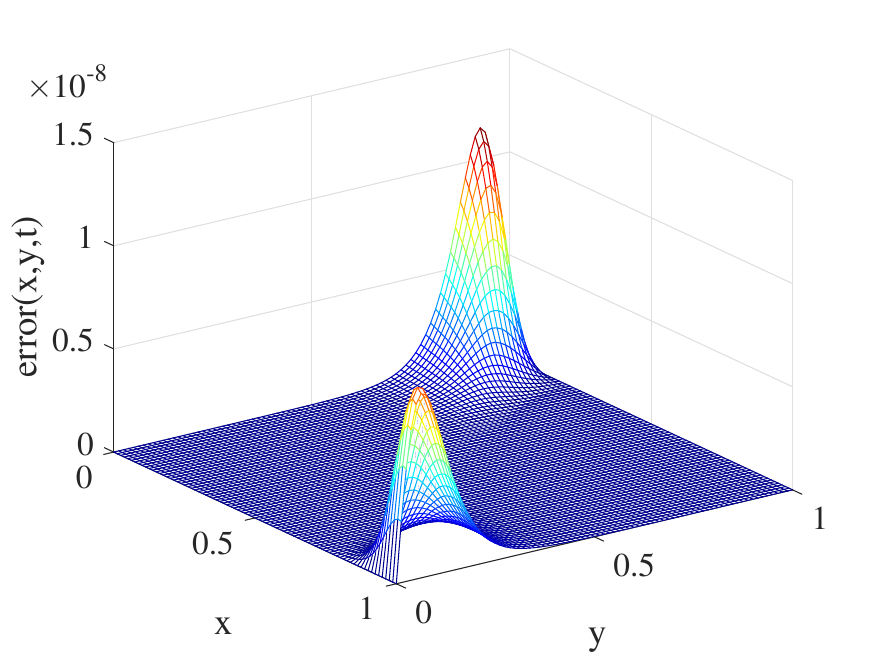}
{
  \caption{ The numerical solution (Left) and the error (Right) of the two-dimensional Burgers' equation for $Re = 200$ at $t = 5.0$ with $N \times N= 81\times 81$ for Example 4. }\label{fi.11}}
  \end{figure}

 The numerical and analytical solutions of two-dimensional examples are present in Tables \ref{ttt},\ref{tttt} with $Re=100$ and $\tau=5\times10^{-5}$. The numerical solutions of the different time are presented in Figs. \ref{figure.10},\ref{figure.11} with $N \times N= 81\times 81$.
From the Tables \ref{ttt},\ref{tttt} and Figs. \ref{fi.11} of numerical simulation results, it can be seen that the proposed numerical scheme has high accuracy under the condition of large Reynolds number ($Re=100,200$). It is observed that the numerical solutions show great agreement with the analytical solutions. The Figs. \ref{figure.10},\ref{figure.11} exhibit that the numerical solution of partial regions becomes steeper and steeper as the increase of time. This physical phenomenon validates the fact that the numerical solution is capable of describing shock wave. Moreover, note that  Tables \ref{ttt},\ref{tttt} and Figs. \ref{figure.10},\ref{figure.11} indicate the property (the boundary condition (\ref{e.72})) of the solution of the two-dimensional Burgers' equation. The physical phenomena depicted in the Figs. \ref{figure.10},\ref{figure.11} are analogical to those in Refs. \cite{Gao2016,Zhang2009,Zhang2010,Siraj-ul-Islam2012}.
\begin{table*}[]
\scriptsize
\centering
\caption{Comparison with the proposed scheme and the analytical solutions of different coordinate positions $(x,y)$ with $Re=100$, at $t=0.25,0.5$ for Example 4.}\label{ttt}
\begin{tabular}{lllllll}
\toprule
\multirow{1}{*}{}
 & \multicolumn{2}{c}{$t=0.25$}
 & \multicolumn{2}{c}{$t=0.5$}\\
\cmidrule(r){2-3} \cmidrule(r){4-5}
$(x,y)$&Analytical solution\cite{Gao2016}    & \textcolor{green}{Proposed} scheme
&Analytical solution\cite{Gao2016}    & \textcolor{green}{Proposed} scheme \\
\midrule
$ (0.25,0.25) $ & 0.3935490117704355 & 0.393549011771103 & 0.2911828920816955 & 0.291182892082807 \\
$ (0.50,0.25) $ & 0.6861822403861596 & 0.686182243447071 & 0.5605467081571704 & 0.560546708625940 \\
$ (0.75,0.25) $ & 0.3935490117704355 & 0.393524523103096 & 0.2911828920816955 & 0.291180266677443 \\
$ (0.25,0.50) $ & 0.2619506158413131 & 0.261950614753384 & 0.2619829318741635 & 0.261982931825796 \\
$ (0.50,0.50) $ & 3.433140255652116E-70 & 0 & 6.648346088881589E-70 & 0 \\
$ (0.75,0.50) $ & -0.2619506158413131 & -0.261950614753467 & -0.2619829318741635 & -0.261982931825921 \\
$ (0.25,0.75) $ & -0.3935490117704355 & -0.393513909760747 & -0.2911828920816955 & -0.291179204683133 \\
$ (0.50,0.75) $ & -0.6861822403861596 & -0.686182243439029 & -0.5605467081571704 & -0.560546708617027 \\
$ (0.75,0.75) $ & -0.3935490117704355 & -0.393549011774452 & -0.2911828920816955 & -0.291182892085496 \\
\bottomrule
\end{tabular}
\end{table*}
\begin{table*}[]
\scriptsize
\centering
\caption{Comparison with the proposed scheme and the exact solutions of different coordinate positions of $(x,y)$ with $Re=100$ and $t=0.75,1.0$ for Example 4.}\label{tttt}
\begin{tabular}{lllllll}
\toprule
\multirow{1}{*}{}
 & \multicolumn{2}{c}{$t=0.75$}
 & \multicolumn{2}{c}{$t=1.0$}\\
\cmidrule(r){2-3} \cmidrule(r){4-5}
$(x,y)$&Analytical solution\cite{Gao2016}    & \textcolor{green}{Proposed} scheme
&Analytical solution\cite{Gao2016}    & \textcolor{green}{Proposed} scheme
 \\
\midrule
$ (0.25,0.25) $ & 0.2273774661403168 & 0.227377466141544 & 0.1858035619888798 & 0.185803561989915 \\
$ (0.50,0.25) $ & 0.4479960634879613 & 0.447996063594949 & 0.3687873026249988 & 0.368787302662550 \\
$ (0.75,0.25) $ & 0.2273774661403168 & 0.227377282040238 & 0.1858035619888798 & 0.185803538891363 \\
$ (0.25,0.50) $ & 0.2180447995061038 & 0.218044799502132 & 0.1818603555820653 & 0.181860355581896 \\
$ (0.50,0.50) $ & 2.301061067978292E-70 & 0 & 7.783304094972026E-71 & 0 \\
$ (0.75,0.50) $ & -0.2180447995061038 & -0.218044799502212 & -0.1818603555820653 & -0.181860355581938 \\
$ (0.25,0.75) $ & -0.2273774661403168 & -0.227377196569641 & -0.1858035619888798 & -0.185803526285714 \\
$ (0.50,0.75) $ & -0.4479960634879613 & -0.447996063584584 & -0.3687873026249988 & -0.368787302652006 \\
$ (0.75,0.75) $ & -0.2273774661403168 & -0.227377466144151 & -0.1858035619888798 & -0.185803561992506 \\
\bottomrule
\end{tabular}
\end{table*}
  \\
\\{\bf{Example 5.}} In order to test the applicability of the proposed scheme for multi-dimensional problems, we consider the three-dimensional Burgers' equations (\ref{equ.7}) over a domain \textcolor{blue}{$[-1,1]\times[-1,1]\times[-1,1]$}, with the initial conditions
\begin{equation}\label{}
\left\{ {\begin{array}{*{20}{l}}
{u\left( {x,y,z,0} \right) = \frac{{ - 2}}{{{\rm{ }}Re}}\frac{{\cos \pi x\sin \pi y\sin \pi z}}{{1 + \sin \pi x\sin \pi y\sin \pi z}}}\\
{v\left( {x,y,z,0} \right) = \frac{{ - 2}}{{\;Re}}\frac{{\sin \pi x\cos \pi y\sin \pi z}}{{1 + \sin \pi x\sin \pi y\sin \pi z}}}\\
{w\left( {x,y,z,0} \right) = \frac{{ - 2}}{{{\rm{ }}Re}}\frac{{\sin \pi x\sin \pi y\cos \pi z}}{{1 + \sin \pi x\sin \pi y\sin \pi z}}}
\end{array}} \right.,(x,y,z) \in \partial \Omega
\end{equation}
The analytical solution for this problem is given by
\begin{equation}\label{}
\left\{ \begin{array}{l}
u = \frac{{ - 2}}{{{\mathop{ Re}\nolimits} }}\frac{{\exp ( - 3\pi^2 \omega t)\cos \pi x\sin \pi y\sin \pi z}}{{1 + \exp ( - 3\pi^2\omega t)\sin \pi x\sin \pi  y\sin \pi z}}\\
v = \frac{{ - 2}}{{{\mathop{\ Re}\nolimits} }}\frac{{\exp ( - 3\pi^2 \omega t)\sin \pi x\cos \pi y\sin \pi z}}{{1 + \exp ( - 3\pi^2\omega t)\sin \pi x\sin \pi y\sin \pi z}}\\
w = \frac{{ - 2}}{{{\mathop{ Re}\nolimits} }}\frac{{\exp ( - 3\pi^2 \omega t)\sin \pi x\sin \pi y\cos \pi z}}{{1 + \exp ( - 3\pi^2\omega t)\sin \pi x\sin \pi y\sin \pi z}}
\end{array} \right.,{\rm{(}}x,y,z{\rm{)}} \in \partial \Omega ,t > 0
\end{equation}
\par
\begin{figure}[h!]
  \includegraphics[width=5in]{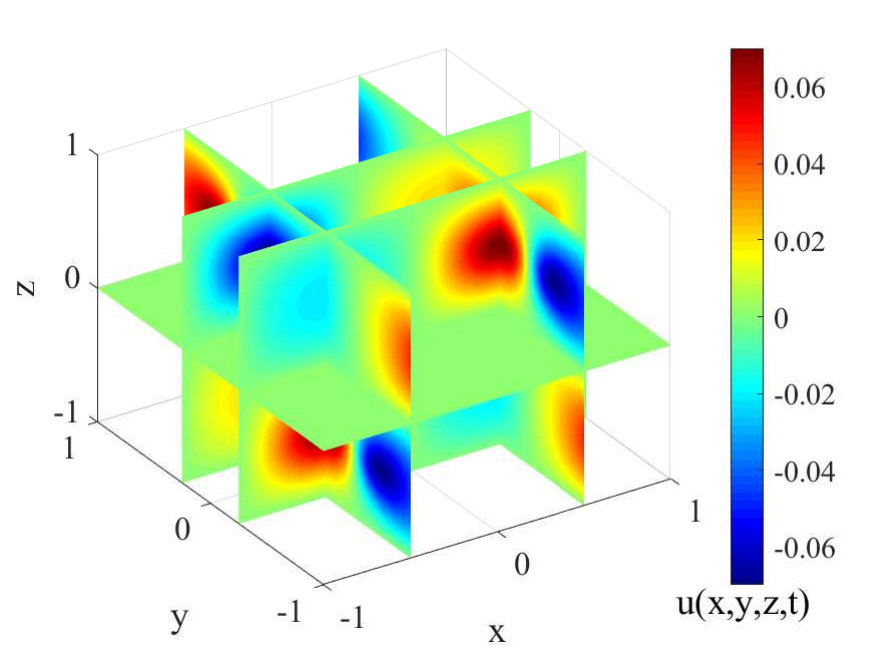}
  \caption{The slices of the four-dimensional figure of the solutions of the three-dimensional
problem for the CFD-PIM-SSM scheme with spatial step size $h_x = h_y = h_z = 1.25 \times 10^{-2}$ and time step size
$\tau = 5 \times10^{-5}$ and $Re=100$ at $t = 1$ for Example 5.}\label{figure.12}
  \end{figure}

\begin{figure}[h!]
  \includegraphics[width=2.35in]{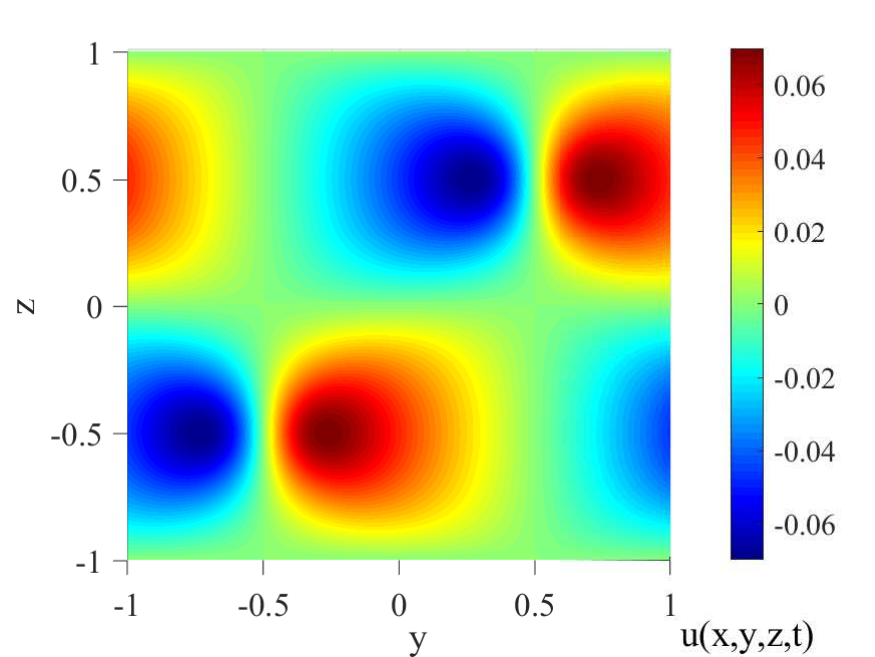}
  \includegraphics[width=2.35in]{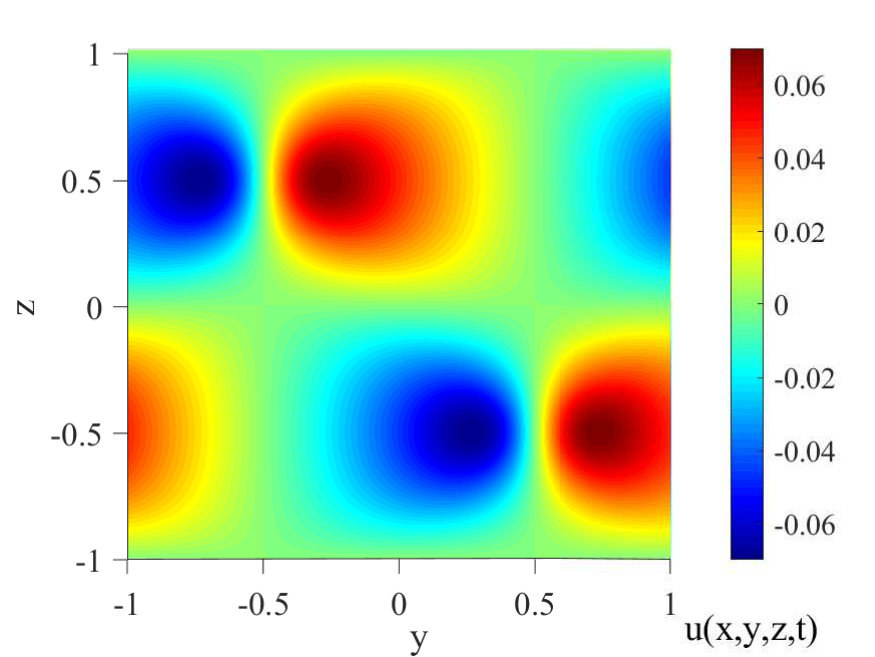}

  \caption{The slices of four-dimensional figures of the solution of the three-dimensional problem
with $x = -0.5$ (Left) and $x = 0.5$ (Right) for Example 5. }\label{figure.13}
  \end{figure}

\begin{figure}[h!]
  \includegraphics[width=2.35in]{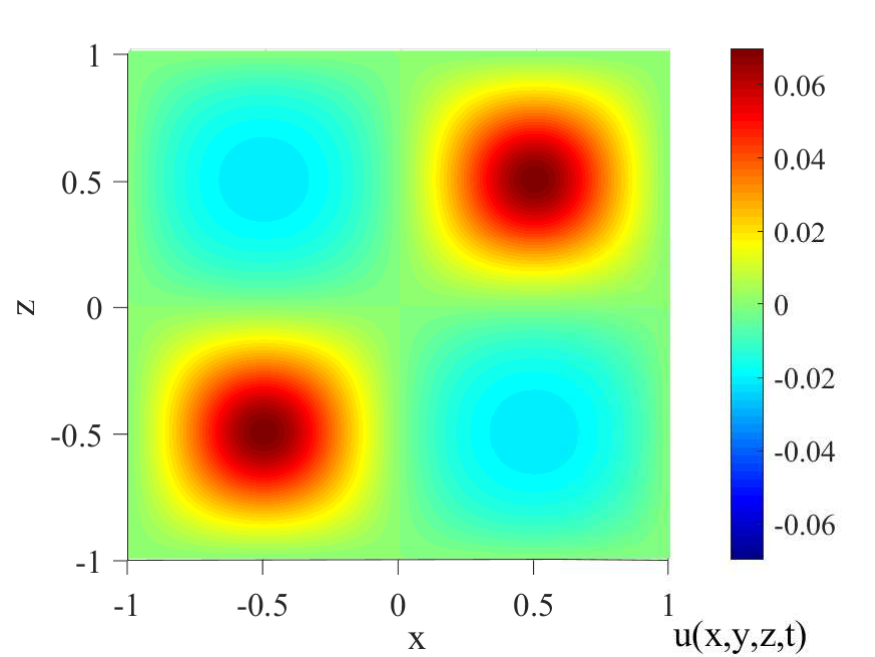}
  \includegraphics[width=2.35in]{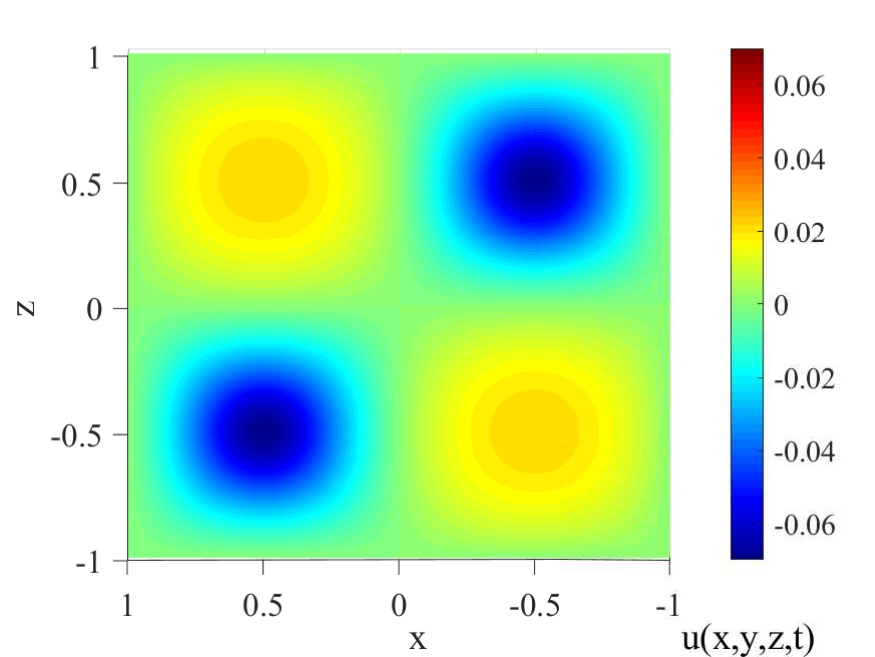}
  \caption{The slices of four-dimensional figures of the solution of the three-dimensional problem
with $y=-0.25$ (Left) and $y=0.25$ (Right) for Example 5. }\label{figure.14}
  \end{figure}

\begin{figure}[h!]
  \includegraphics[width=2.35in]{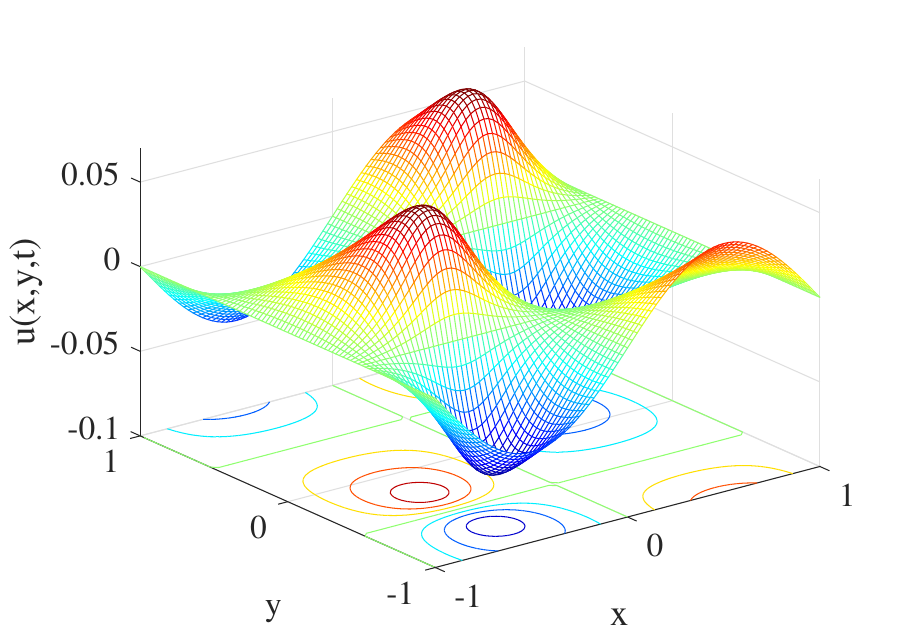}
  \includegraphics[width=2.35in]{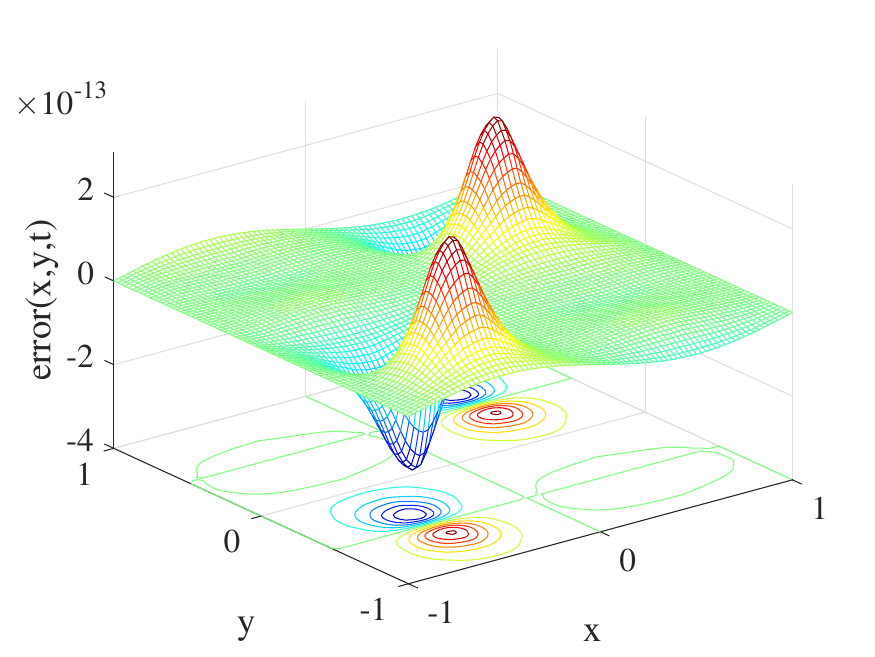}
  \caption{The three-dimensional figures of the solution (Left) and the error (Right) of
the three-dimensional problem with $z=0.25$ and $Re=100$ at $t = 1$ for Example 5. }\label{figure.15}
  \end{figure}

\begin{figure}[h!]
  \includegraphics[width=2.35in]{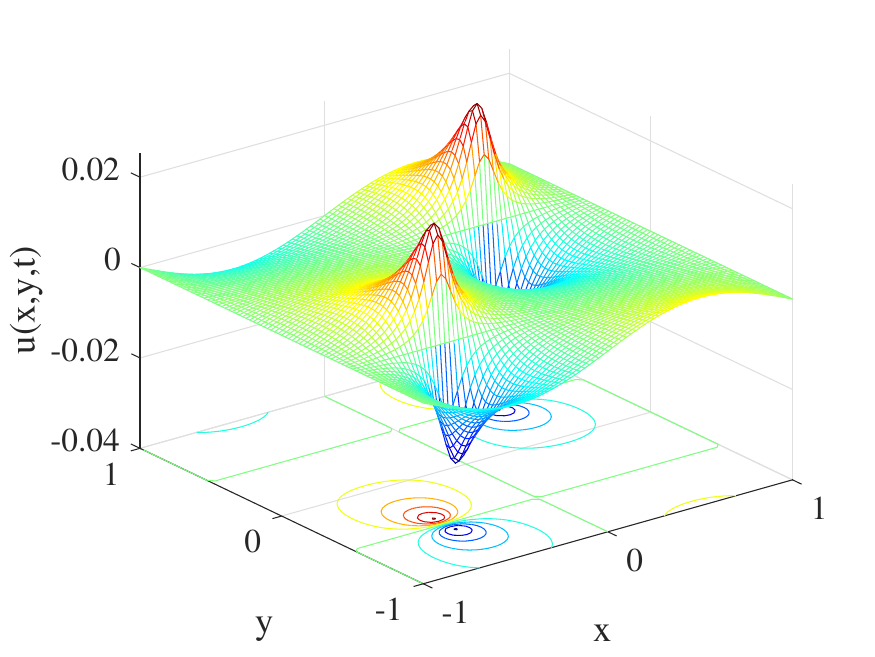}
  \includegraphics[width=2.35in]{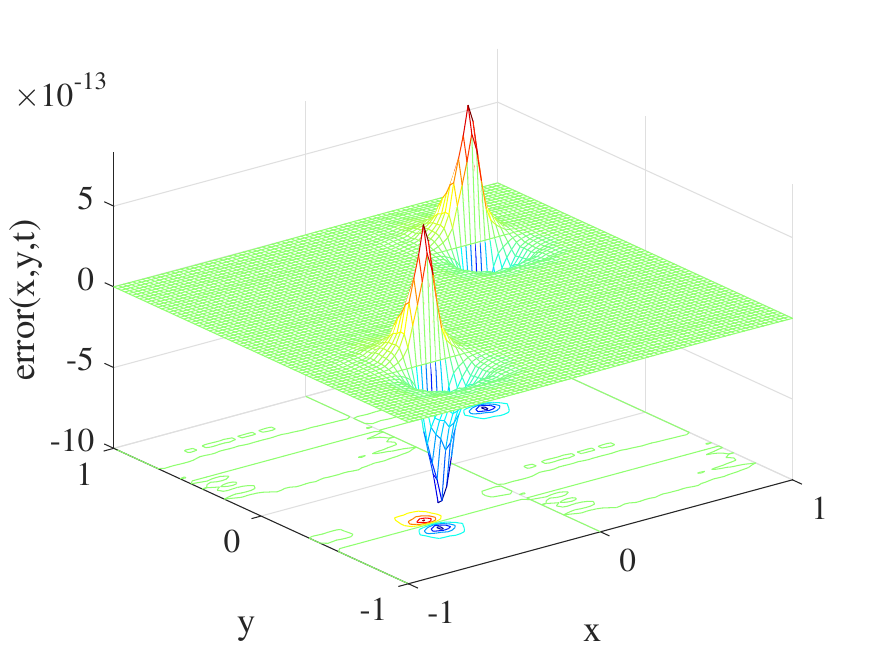}
  \caption{The three-dimensional figures of the solution (Left) and the error (Right) of
the three-dimensional problem with $z=0.25$ and $Re=1000$ at $t = 1$ for Example 5. }\label{figure.16}
  \end{figure}

 The boundary conditions are extracted from the analytical solution.
It is observed that the numerical solutions present great agreement with the analytical solutions.
The slices of the four-dimensional images are used for observing three-dimensional Burgers' equations, which are depicted in Figs. \ref{figure.12},\ref{figure.13},\ref{figure.14}. Figs. \ref{figure.15},\ref{figure.16} present the numerical solutions and the errors between the exact solutions and the numerical solutions at $t = 1$ with $z=0.25$ and $Re=100,1000$. The numerical results of this example show that the CFD-PIM-SMM scheme base on Hopf-Cole transformation is a numerical method with high precision and high efficiency for solving n-dimensional Burgers' system. The comparison is done with solutions obtained by the numerical solutions and the exact solutions for the three-dimensional Burgers' equations, which is presented in on left-hand side Figs. \ref{figure.15},\ref{figure.16}.
All numerical figures of this example show that the proposed scheme has excellent accuracy and efficiency.
Figs. \ref{figure.15},\ref{figure.16} exhibit that the numerical solution of partial regions becomes steeper and steeper as $Re=100\longrightarrow1000$.
This physical phenomenon validates the fact that the numerical solution is capable of describing shock wave.
 { Table \ref{newta8} shows the numerical results of different  $\kappa$ for Burgers' equation. Table \ref{newta8} mainly shows two phenomena: 1. With the increase of $\kappa$, the corresponding error will also increase;  2. Table \ref{newta8} transformation does not cause ROC transformation;  Phenomenon 2 is consistent with n-dimensional Hopf-Cole transformation (Eqs. (\ref{eeqq.20}) to (\ref{eeqq.23})) which has no relationship between $\kappa$ and ROC.
}
\begin{table}[h!]
\scriptsize
\centering
{
\caption{\newline Numerical results of $\kappa=1,2,5,10$ with $t=1$ and $Re=100$ for Example 5.}\label{newta8}
      \begin{tabular}{lllllll}
  \toprule
  $\kappa$ & $N$ & $11\times11$ & $21\times21$ & $41\times41$ & $81\times81$ \\
              \midrule
     1 & ${L_\infty }$ & 1.1483E-07 & 2.7964E-10 & 8.8463E-13 & 7.2616E-15 \\
          & ROC & -- & 8.6817 & 8.3051 & 7.6316 \\
          2 & ${L_\infty }$ & 5.7416E-08 & 1.3982E-10 & 4.4232E-13 & 3.6308E-15 \\
          & ROC & -- & 8.6817 & 8.3051 & 7.6316 \\
          5 & ${L_\infty }$ & 2.2966E-08 & 5.5929E-11 & 1.7693E-13 & 1.4520E-15 \\
           & ROC & -- & 8.6817 & 8.3051 & 7.6316 \\
         10 & ${L_\infty }$ & 1.1483E-08 & 2.7964E-11 & 8.8463E-14 & 7.2598E-16 \\
           & ROC & -- & 8.6817 & 8.3051 & 7.6316 \\
  \bottomrule
\end{tabular}
 }
\end{table}
\begin{figure}[h!]
  \includegraphics[width=5in]{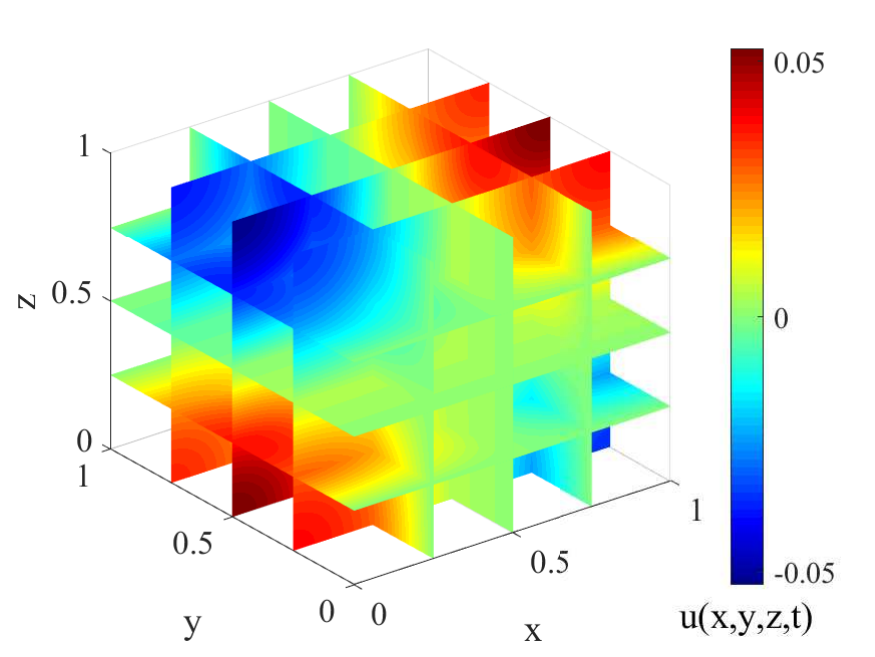}
  \caption{The slices of the four-dimensional figure of the solutions of the three-dimensional
problem for the CFD-PIM-SSM scheme with spatial step size $h_x = h_y = h_z = 1.25\times 10^{-2}$ and time step size
$\tau = 5 \times10^{-5}$ and $Re=10$ at $t = 1$ for Example 6.}\label{figure.18}
  \end{figure}
\begin{figure}[h!]
  \includegraphics[width=2.35in]{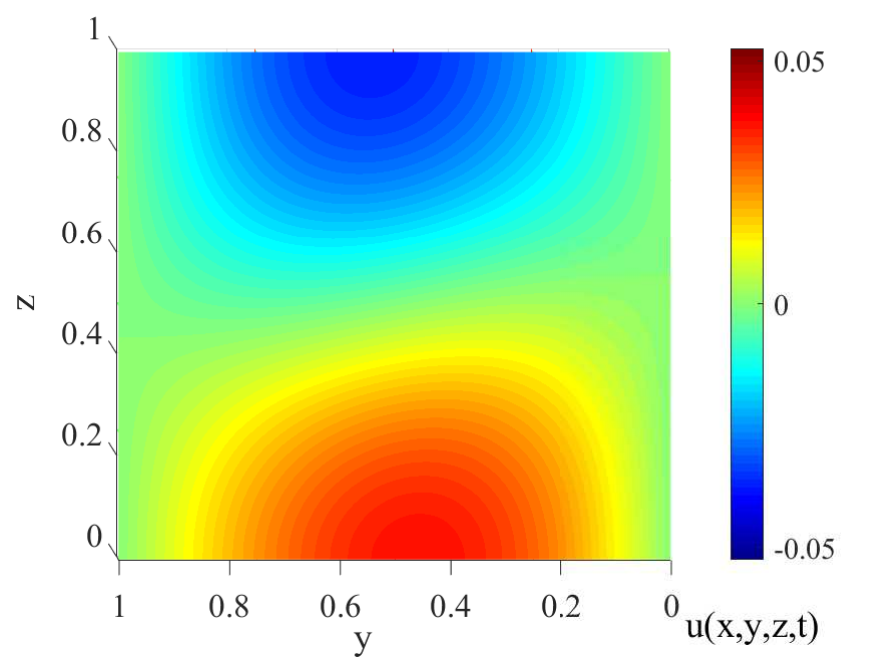}
  \includegraphics[width=2.35in]{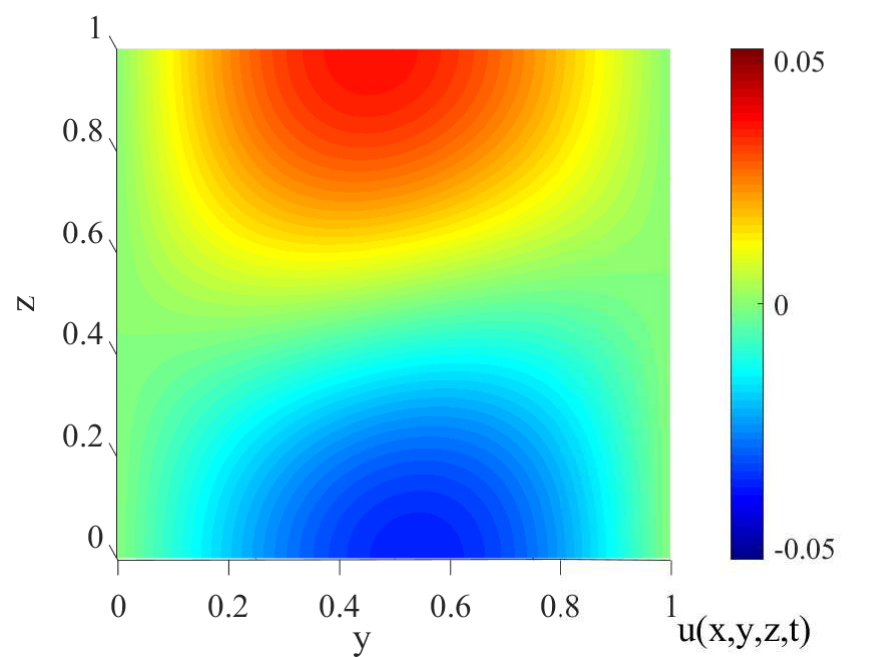}
  \caption{The slices of four-dimensional figures of the solution of the three-dimensional problem
with $x = 0.25$ (Left) and $x=0.75$ (Right) for Example 6. }\label{figure.19}
  \end{figure}
  \begin{figure}[h!]
  \includegraphics[width=2.35in]{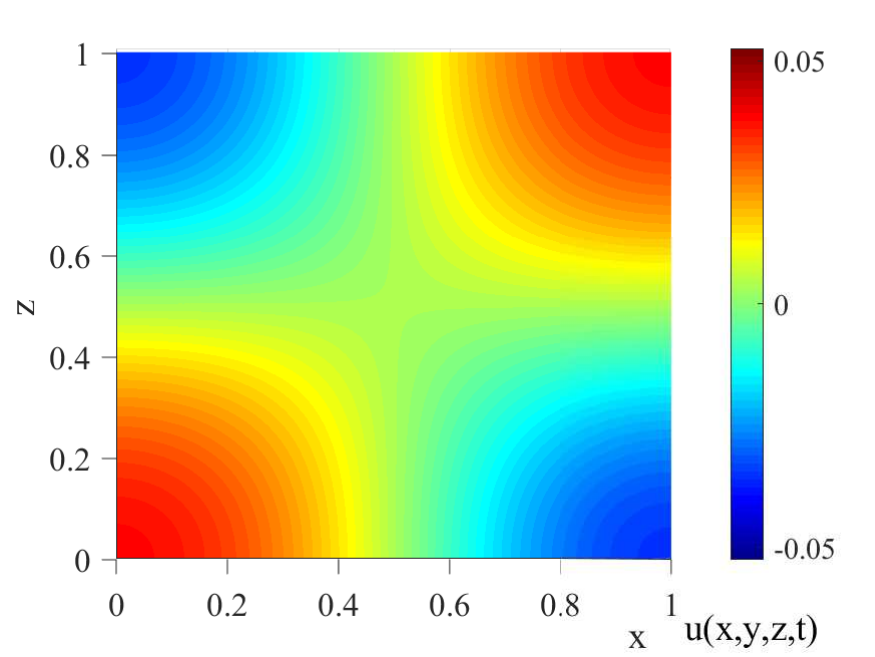}
  \includegraphics[width=2.35in]{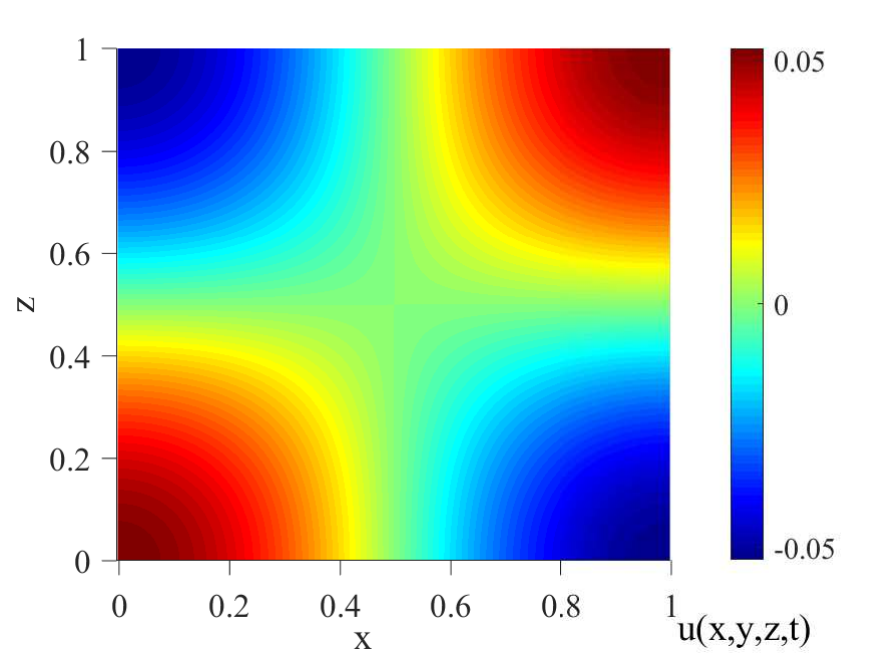}
  \caption{The slices of four-dimensional figures of the solution of the three-dimensional problem
with $y = 0.25$ (Left) and $y = 0.5$ (Right) for Example 6. }\label{figure.20}
  \end{figure}
\begin{figure}[h!]
  \includegraphics[width=2.35in]{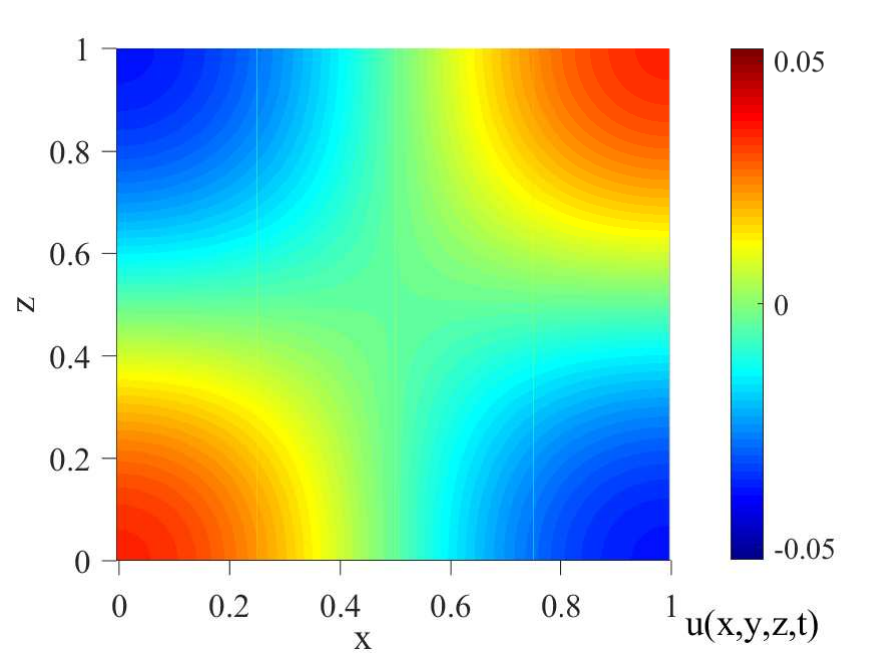}
  \includegraphics[width=2.35in]{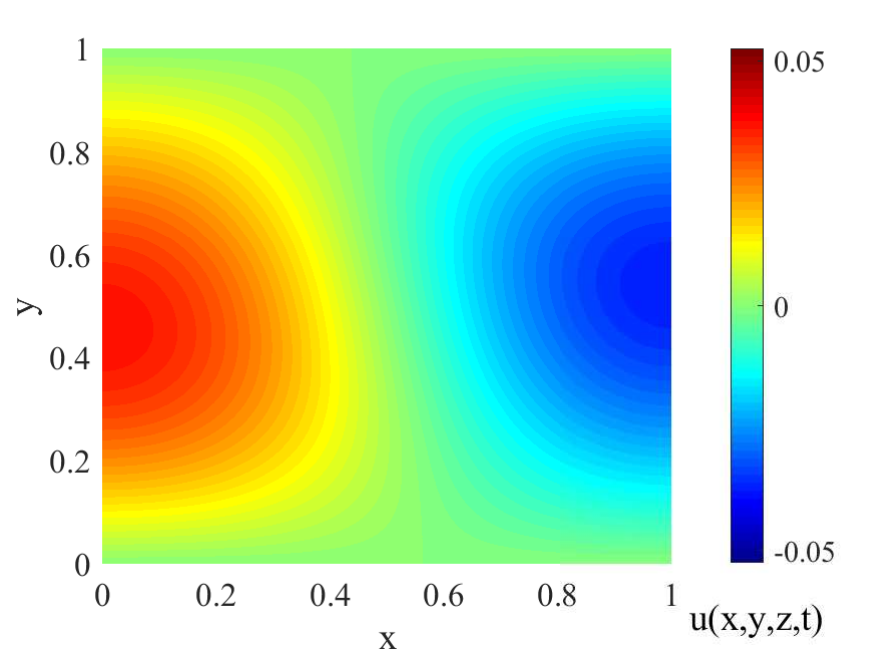}
  \caption{The slices of four-dimensional figures of the solution of the three-dimensional problem
with $y=0.75$ (Left) and $z=0.25$ (Right) for Example 6. }\label{figure.21}
  \end{figure}
\begin{figure}[h!]
  \includegraphics[width=2.35in]{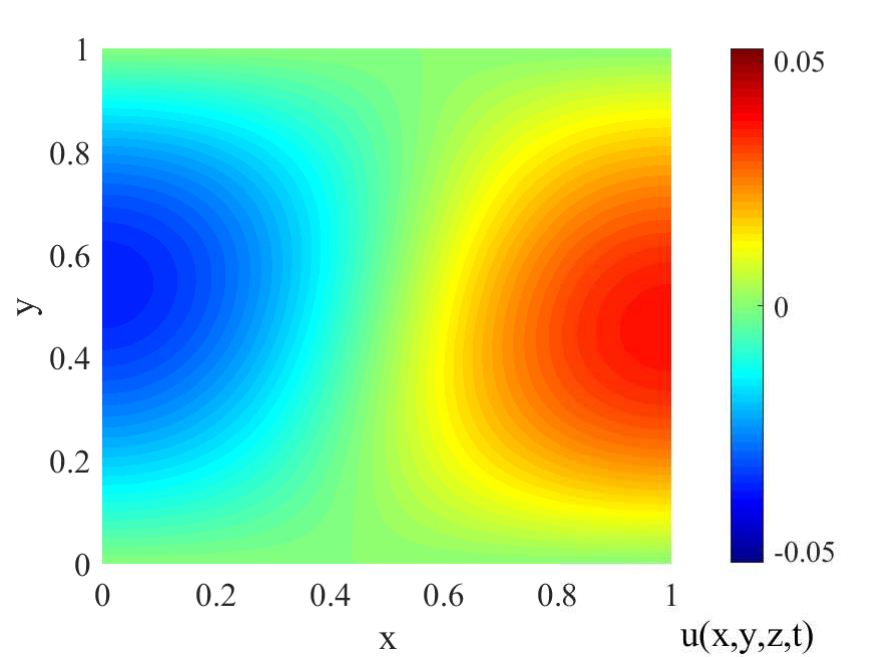}
  \includegraphics[width=2.35in]{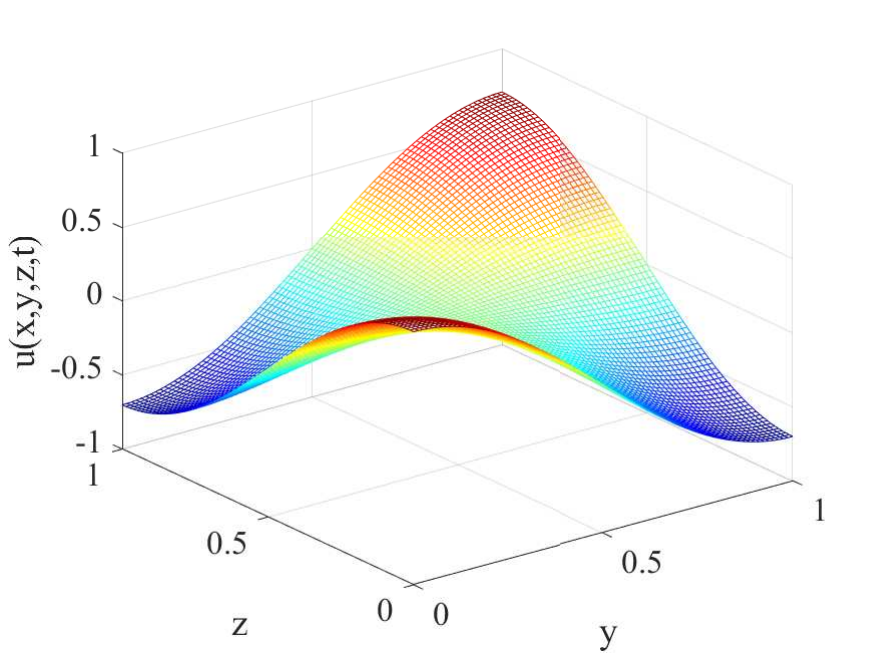}
  \caption{The slices of four-dimensional figures of the solution of the three-dimensional problem
with $z=0.75$ (Left) and the numerical solution of the three-dimensional Burgers' equation at $t = 0$ (Right) with $N \times N= 81\times 81$, $x=0.25$ and $Re = 10$ for Example 6. }\label{figure.22}
  \end{figure}
\begin{figure}[h!]
  \includegraphics[width=2.35in]{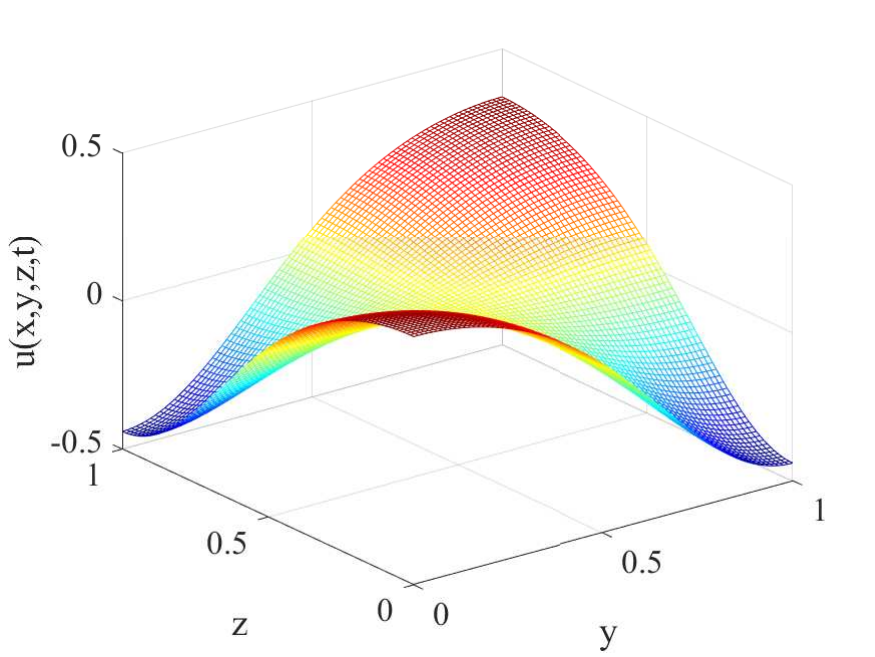}
  \includegraphics[width=2.35in]{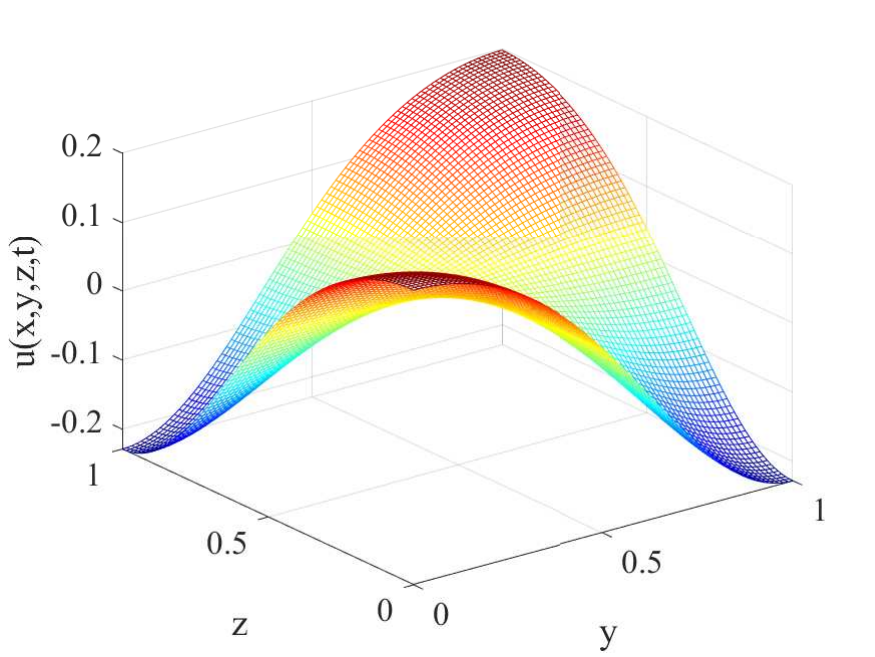}
  \caption{The numerical solution of the three-dimensional Burgers' equation at $t = 0.2$ (Left) and $t = 0.4$ (Right) with $N \times N= 81\times 81$, $x=0.25$ and $Re = 10$ for Example 6. }\label{figure.23}
  \end{figure}
\begin{figure}[h!]
  \includegraphics[width=2.35in]{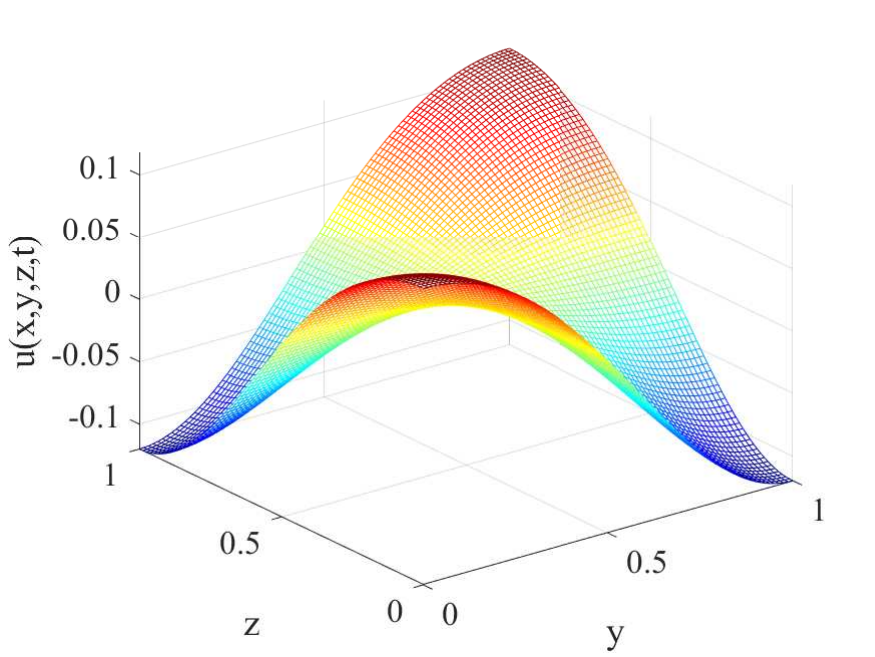}
  \includegraphics[width=2.35in]{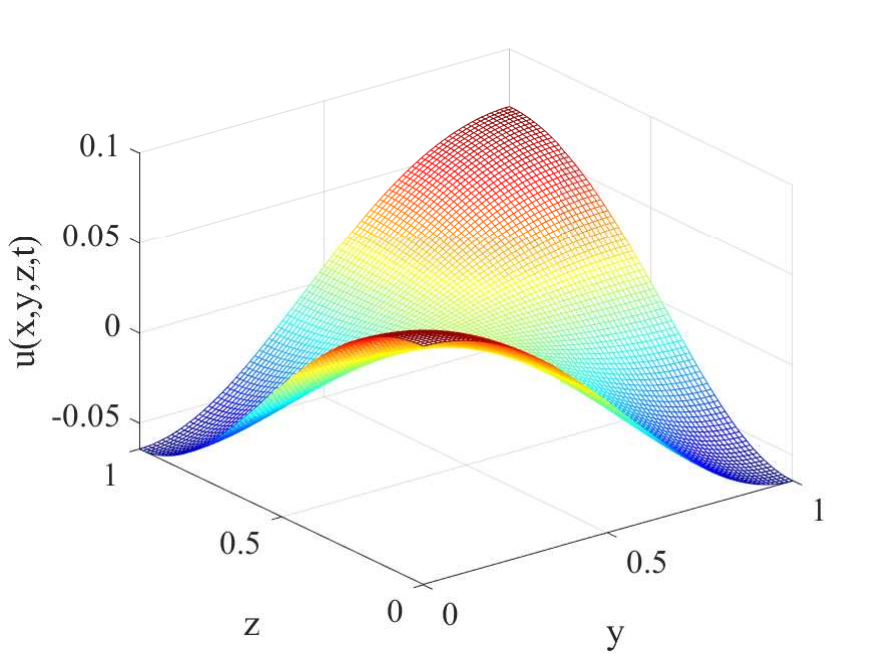}

  \caption{The numerical solution of the three-dimensional Burgers' equation at $t = 0.6$ (Left) and $t = 0.8$ (Right) with $N \times N= 81\times 81$, $x=0.25$ and $Re = 10$ for Example 6. }\label{figure.24}
  \end{figure}
  \begin{figure}[h!]
  \includegraphics[width=2.35in]{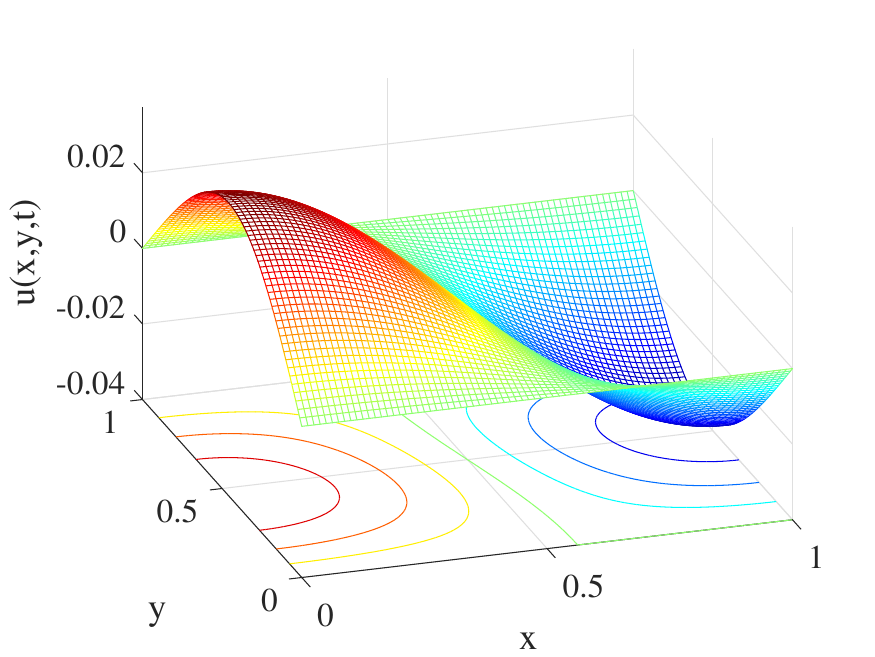}
  \includegraphics[width=2.35in]{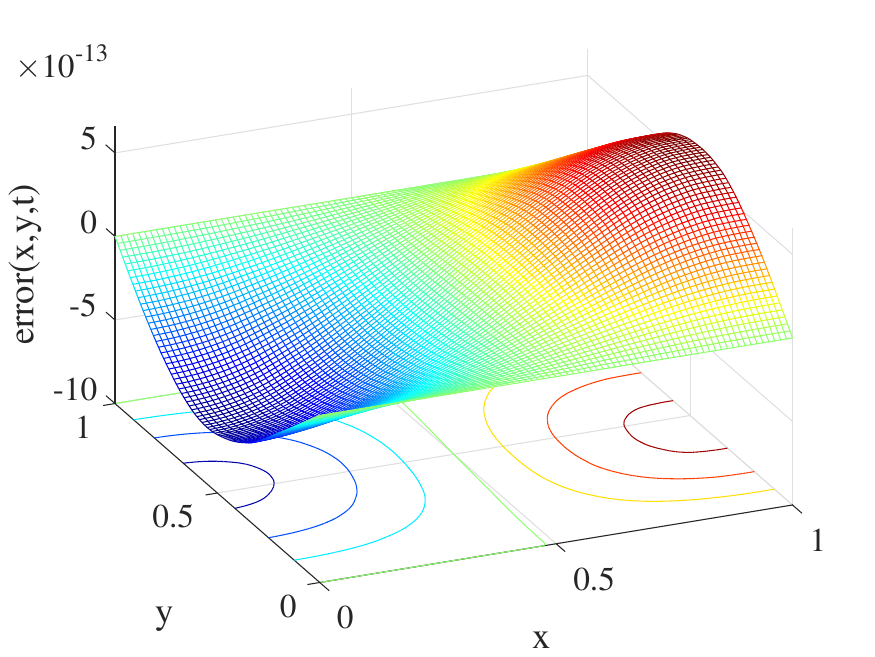}
  \caption{The three-dimensional figures of the solution (Left) and the error (Right) between the numerical and analytical solutions of the three-dimensional problem with $z=0.25$ and $Re=10$ at $t = 1$ for Example 6. }\label{figure.26}
  \end{figure}
\\
\\{\bf{Example 6.}} We consider the system of {the} three-dimensional Burgers' equation (\ref{equ.6}) in Sec.\ref{2.2.3}, over a square domain $[0,1]\times[0,1]\times[0,1]$, with the initial and boundary conditions (\ref{80},\ref{81}), and the analytical solution
\begin{equation}\label{}
\begin{array}{l}
u(x,y,z,t) = 2\pi \omega \frac{{\sum\limits_{\alpha ,\beta ,\gamma  = 0}^\infty  \alpha  {C_{\alpha \beta \gamma} }\exp [ - ({\alpha ^2} + {\beta ^2}+{\gamma^2}){\pi ^2}\omega t]\sin (\alpha \pi x)\cos (\beta \pi y)\cos (\gamma \pi z)}}{{\sum\limits_{\alpha ,\beta,\gamma  = 0}^\infty  {{C_{\alpha \beta \gamma} }\exp [ - ({\alpha ^2} + {\beta ^2} + {\gamma^2}){\pi ^2}\omega t]\cos (\alpha \pi x)\cos (\beta \pi y)\cos (\gamma \pi z)} }}\\
v(x,y,z,t) = 2\pi \omega \frac{{\sum\limits_{\alpha ,\beta ,\gamma = 0}^\infty  \beta  {C_{\alpha \beta \gamma} }\exp [ - ({\alpha ^2} + {\beta ^2}+{\gamma^2}){\pi ^2}\omega t]cos(\alpha \pi x)\sin (\beta \pi y)\cos (\gamma \pi z)}}{{\sum\limits_{\alpha ,\beta ,\gamma = 0}^\infty  {{C_{\alpha \beta \gamma} }\exp [ - ({\alpha ^2} + {\beta ^2}+{\gamma^2}){\pi ^2}\omega t]\cos (\alpha \pi x)\cos (\beta \pi y)\cos (\gamma \pi z)} }}\\
w(x,y,z,t) = 2\pi \omega \frac{{\sum\limits_{\alpha ,\beta,\gamma  = 0}^\infty  \gamma  {C_{\alpha \beta \gamma} }\exp [ - ({\alpha ^2} + {\beta ^2}+{\gamma^2}){\pi ^2}\omega t]cos(\alpha \pi x)\cos (\beta \pi y)\sin (\gamma \pi z)}}{{\sum\limits_{\alpha ,\beta,\gamma  = 0}^\infty  {{C_{\alpha \beta \gamma} }\exp [ - ({\alpha ^2} + {\beta ^2}+{\gamma^2}){\pi ^2}\omega t]\cos (\alpha \pi x)\cos (\beta \pi y)\cos (\gamma \pi z)} }}
\end{array}
\end{equation}
\begin{table}[h!]
\scriptsize
\centering
\caption{The solutions of the three-dimensional heat conduction equation with $Re = 10$ at $t = 0.1$ for Example 6.(The spatial step sizes of the FEM are $h_x=h_y=h_z=10^{-2}$)}\label{t18}
\begin{tabular}{lllllll}
\toprule
\multirow{3}{*}{}
 & \multicolumn{3}{c}{$\phi(x,y,z,t)$}\\
\cmidrule(r){2-4}
$(x,y,z)$ & \multirow{1}{*}{ Analytical solution\cite{Gao2016} }
&\multirow{1}{*}{ FEM\cite{Gao2016}}& \textcolor{green}{Proposed} scheme\\
\midrule
$ (0.25,0.00,0.00) $ & 0.5121139094268445 & 0.5121574025971580 & 0.512112211115770 \\
$ (0.25,0.25,0.00) $ & 0.4021754605524777 & 0.4021832308799626 & 0.402175090853434 \\
$ (0.25,0.25,0.25) $ & 0.3332558679440735 & 0.3332624499504631 & 0.333255848934007 \\
$ (0.25,0.50,0.00) $ & 0.2221335878372145 & 0.2221195126283573 & 0.222134351589450 \\
$ (0.25,0.50,0.25) $ & 0.2145599507332575 & 0.2145578143852998 & 0.214560278902764 \\
$ (0.25,0.75,0.00) $ & 0.1232224673074507 & 0.1232115116009095 & 0.123222177969065 \\
$ (0.50,0.00,0.00) $ & 0.2350307713064107 & 0.2349996700538997 & 0.235032545451199 \\
$ (0.50,0.25,0.00) $ & 0.2221335878372145 & 0.2221195160032975 & 0.222134351729473 \\
$ (0.50,0.25,0.25) $ & 0.2145599507332575 & 0.2145578324085150 & 0.214560278977476 \\
$ (0.75,0.00,0.00) $ & 0.09757645937131027 & 0.09756832939249252 & 0.0975749310522441 \\
$ (0.75,0.25,0.00) $ & 0.1232224673074507 & 0.1232114495876140 & 0.123222177969033 \\
\bottomrule
\end{tabular}
\end{table}
\par In order to verify the validity and practicability of the proposed scheme of three-dimensional Burgers' equation, the numerical solution is compared with the analytical solution.
The numerical solution of the three-dimensional Burgers' equation (\ref{e.12}) is given by the three-dimensional Hopf-Cole transformation (\ref{e.84}).
The correctness of the numerical solution of the heat conduction equation ensures the correctness of the numerical solution of the Burgers' equation.
The higher the accuracy of numerical solution of the heat conduction equation,
the higher the accuracy of the numerical solution of the Burgers' equation will be.
The numerical and analytical solutions of the three-dimensional heat conduction equation are presented in Table \ref{t18} with $Re=10$  and $\tau=5\times10^{-5}$.
From the Table \ref{t18} of numerical simulation results, it is observable that the numerical solutions present excellent consistency with the analytical solutions.
The slices of the four-dimensional images are used for observing three-dimensional Burgers' equation, which are depicted in Figs. \ref{figure.19},\ref{figure.20},\ref{figure.21},\ref{figure.22},\ref{figure.23}.
The numerical solutions of the different time are presented in Figs. \ref{figure.22},\ref{figure.23},\ref{figure.24} with $N \times N \times N= 81\times 81 \times 81$.
It can be observed from Figs. \ref{figure.22},\ref{figure.23},\ref{figure.24} when $Re=10$, the numerical solution of the Burgers' equation will generate shock wave.
This physical phenomenon validates the fact that the numerical solution is capable of describing shock wave.
The physical phenomena depicted in the Figs. \ref{figure.22},\ref{figure.23},\ref{figure.24} are similar to those in Ref.\cite{Gao2016}.
We make a three-dimensional Burgers' equation with numerical solutions and errors in Fig. \ref{figure.26}.
It is observed that the numerical solution is very consistent with the analytical solution.
\section{Conclusion}
This paper proposed the modified Hopf-Cole transformation for n-dimensional Burgers' system. After obtaining the n-dimensional heat conduction equations, we present a new high-order exponential time differencing precise integration method schemes in combination with the sixth-order compact finite difference scheme, which have been developed for the numerical solutions of n-dimensional heat conduction equations.
The proposed scheme is tested on \textcolor{blue}{six} examples, and the results are entirely satisfactory in comparison with the analytical solutions. The findings can be summarized as follows

\begin{enumerate}[(1)]
\item For the examples of the ordinary initial and boundary conditions, the solutions can describe shock wave phenomena for large Reynolds numbers ($Re\geq 1000$), which is characterized by high precision and high efficiency.

\item For the complicated initial and boundary conditions problems, the proposed scheme obtains excellent accuracy and high efficiency in comparison to the other numerical techniques available in the Refs. \cite{Dubansky2017,Bhatt2016,Kadalbajoo2006,Jiwari2015,Jiwari2012,Seydaoglu2018,G.W.Wei1998,Mittal2012a,Sheu2002,Kutluay1999,Lai2014}.
\item After the modification of Hopf-Cole transformation and optimization of computer programming, the present scheme has commendable adaptability and high efficiency in the calculation examples.

\item  It is found that the proposed results are in good agreement with the analytical solutions for two-dimensional and three-dimensional problems. This scheme includes linear problems and nonlinear problems, which can be easily extended to solve model equations of high-dimensional problems.

\end{enumerate}

\section*{Code}
We are happy to share our research results with others. We have uploaded some example codes (Example 1, 4, and 5) on Github. This is Github's link: \url{https://github.com/LzEfreet/CFD-PIM}.

\section*{Acknowledgement}
This work is financially supported by the National Natural Science Foundation of China(No. 11826208).

\section*{References}

\bibliography{mybibfile}

\end{document}